\documentclass[preprint,3p,12pt]{elsarticle}
\usepackage[T1]{fontenc}
\usepackage{lmodern}
\usepackage{amsmath,amssymb,amsthm,mathtools}
\usepackage{microtype}
\usepackage{enumitem}
\usepackage{array,etoolbox}
\usepackage{needspace}
\usepackage{xurl}
\usepackage[hidelinks]{hyperref}
\usepackage{bookmark}
\biboptions{numbers,sort&compress}
\hypersetup{pdftitle={Scott spectra of distributive continuous lattices},
 pdfauthor={Yuxu Chen},
 pdfsubject={An order-theoretic characterization of Scott spectra},
 pdfkeywords={continuous lattice, spectrum, Scott topology, coframe}}
\numberwithin{equation}{section}
\newtheorem{theorem}{Theorem}[section]
\newtheorem{lemma}[theorem]{Lemma}
\newtheorem{proposition}[theorem]{Proposition}
\newtheorem{corollary}[theorem]{Corollary}
\theoremstyle{definition}
\newtheorem{definition}[theorem]{Definition}

\newcommand{\id}{\mathrm{id}}
\newcommand{\Fin}{\mathcal P_{\mathrm{fin}}}
\newcommand{\Nbh}{\mathcal N}
\newcommand{\Spec}{\operatorname{Spec}}
\newcommand{\Prm}{\operatorname{Pr}}
\newcommand{\OO}{\mathcal O}
\newcommand{\GG}{\Gamma}
\newcommand{\NN}{\mathbb N}
\newcommand{\RR}{\mathbb R}
\newcommand{\op}{\mathrm{op}}
\newcommand{\cl}{\operatorname{cl}}
\newcommand{\intr}{\operatorname{int}}
\newcommand{\Pc}{\operatorname{Cop}}
\newcommand{\up}{\mathord\uparrow}
\newcommand{\dn}{\mathord\downarrow}
\newcommand{\wb}{\ll}
\newcommand{\lw}{\lambda}
\newcommand{\hit}[1]{\Diamond #1}
\newcommand{\hata}[1]{\widehat{#1}}
\setlist[enumerate]{leftmargin=*,itemsep=2pt}
\begin{document}
\begin{frontmatter}
\title{Scott spectra of distributive continuous lattices}
\author{Yuxu Chen}
\ead{chenyuxu@scu.edu.cn}
\address{School of Mathematics, Sichuan University, Chengdu 610064, China}

\begin{abstract}
Let \(L\) be a distributive continuous lattice. We prove that the
hull--kernel topology on \(\Spec(L)\) is the Scott topology of the
inclusion order if and only if
\(\sigma(L^{\op})=\upsilon(L^{\op})\). The necessary implication
and the algebraic case are known. We establish sufficiency without
algebraicity or countability assumptions, obtaining the spectrum
characterization requested by Lawson and Mislove. We also show that
a spatial frame is weak hypercontinuous if and only if its dual
Scott space is core-compact.
\end{abstract}
\end{frontmatter}

\section{Introduction}\label{sec:intro}

The spectral representation of distributive continuous lattices connects
domain theory with the topology of locally compact sober spaces. By the
Hofmann--Lawson representation theorem, every distributive continuous
lattice \(L\) is isomorphic to the lattice of open subsets of its
spectrum \(\Spec(L)\) \cite{HL78}. The points of \(\Spec(L)\) are the
proper completely prime filters of \(L\), and its hull--kernel topology
is generated by the sets
\[
 \hata a=\{F\in\Spec(L):a\in F\},\qquad a\in L.
\]
The specialization order of this topology is inclusion. However, the
representation theorem does not imply that the hull--kernel topology
is the Scott topology of this order.

Lawson and Mislove asked which distributive continuous lattices have
spectra that are sober locally compact Scott spaces
\cite[Problem~528]{LM90}. Since local compactness and sobriety already
follow from the representation theorem, the question is whether the
hull--kernel and Scott topologies coincide on the spectrum. An
order-theoretic characterization must therefore distinguish this
additional property from continuity and distributivity of \(L\).

For a poset \(P\), let \(\sigma(P)\) denote its Scott topology and
\(\upsilon(P)\) its upper topology, generated by
\(P\setminus\dn p\), \(p\in P\). The condition
\(\sigma(L^{\op})=\upsilon(L^{\op})\) is necessary for the spectrum
of a distributive continuous lattice \(L\) to be a Scott space
\cite[Theorem~4.15]{CK23}. For distributive algebraic lattices, the
converse is also contained in that theorem. Chen and Kou asked whether
it remains valid for distributive continuous lattices
\cite[Problem~4.12]{CK23}. We prove this implication without algebraicity
or any countability assumption.

\Needspace{7\baselineskip}
\begin{theorem}[Spectrum characterization]\label{thm:main}
Let \(L\) be a distributive continuous lattice. The following conditions
are equivalent.
\begin{enumerate}[label=\textup{(\roman*)}]
\item The hull--kernel topology on \(\Spec(L)\) is the Scott topology
of the inclusion order.
\item \(\sigma(L^{\op})=\upsilon(L^{\op})\).
\end{enumerate}
\end{theorem}

The implication \textup{(i)}\(\Rightarrow\)\textup{(ii)} and the
algebraic case are known. The converse has three steps. First, we
construct a continuous frame-coframe \(T(L)\) and a frame embedding
\(h_L:L\hookrightarrow T(L)\) that preserves the way-below relation.
The latter property is equivalent to Scott continuity of its right
adjoint and is essential for the transfer argument. Second, a
compact-neighbourhood construction shows that the spectrum of every
continuous frame-coframe is a Scott space; equivalently, every locally
compact sober web space carries its Scott topology. Third, we prove
that the Scott property of the spectrum passes along a
way-below-preserving frame embedding.

We also prove that a spatial frame is continuous and satisfies
\(\sigma(L^{\op})=\upsilon(L^{\op})\) whenever its dual Scott space
is core-compact. Applied to \(L=\OO(X)\), this gives
\[
 \Sigma\GG(X)\text{ is core-compact}
 \quad\Longrightarrow\quad
 \begin{cases}
 X\text{ is core-compact},\\
 \sigma(\GG(X))=\upsilon(\GG(X)),
 \end{cases}
\]
where \(\GG(X)\) denotes the lattice of closed subsets of \(X\).
This answers \cite[Problem~4.9]{CK23} affirmatively. The proof of this
implication does not use Theorem~\ref{thm:main}.

Section~\ref{sec:prelim} fixes the notation and records the spectral
representation. Section~\ref{sec:host} constructs the extension,
Section~\ref{sec:coframes} proves the topological criterion, and
Section~\ref{sec:transfer} establishes the transfer theorem and
completes the proof of Theorem~\ref{thm:main}.
Appendix~\ref{app:spectrum-proof} gives a longer alternative proof,
using a countable topological criterion, metrizable lattice quotients,
and explicit approximation nets. Appendix~\ref{app:dual-compactness}
proves the continuity-reflection result and relates lattice properties
to the order structure of spectra. Appendix~\ref{paper:topology}
treats Scott sobrification and sobriety criteria for web spaces.
Appendix~\ref{paper:functions} classifies function lattices and compares
Scott and bi-Scott retracts. Appendix~\ref{paper:auxiliary} gives the
universal property of the join-translation extension and its structural
correspondences.

\section{Preliminaries}\label{sec:prelim}

We use the standard terminology of domain theory \cite{GHK03}.
Directed and filtered subsets are nonempty, and \(\NN=\{0,1,\ldots\}\).
For a set \(S\), put
\[
 \Fin(S)=\{F\subseteq S:|F|<\infty\},\qquad
 \Fin^+(S)=\Fin(S)\setminus\{\varnothing\}.
\]
Both families are ordered by inclusion; \(\Fin(S)\) is always directed,
and \(\Fin^+(S)\) is directed when \(S\ne\varnothing\).
For a poset \(P\) and
\(A\subseteq P\), write
\[
 \up A=\{x\in P:(\exists a\in A)\ a\leq x\},\qquad
 \dn A=\{x\in P:(\exists a\in A)\ x\leq a\}.
\]
A dcpo is a poset in which every directed subset has a supremum.
The way-below relation is defined by \(u\wb x\) if, for every
directed subset \(D\subseteq P\) whose supremum exists,
\[
 x\leq\bigvee D\quad\Longrightarrow\quad
 (\exists d\in D)\ u\leq d.
\]
A dcpo is continuous if
\(\{u:u\wb x\}\) is directed with supremum \(x\), for every \(x\).
For a continuous dcpo, the relation \(\wb\) interpolates and
\(\Uparrow u=\{x:u\wb x\}\) is Scott open. We write
\(\Sigma P=(P,\sigma(P))\).

The least and greatest elements of a complete lattice are denoted by
\(0\) and \(1\). A frame is a complete lattice in which finite meets distribute over
arbitrary joins. A coframe satisfies the dual distributive law.
A frame homomorphism preserves arbitrary joins and finite meets;
a complete lattice homomorphism preserves arbitrary joins and arbitrary
meets. Empty operations are included. Every distributive continuous
lattice is a frame, by meet-continuity and finite distributivity
\cite[Chapter~I]{GHK03}.
A nucleus on a frame is a monotone inflationary idempotent map preserving
finite meets.

\begin{definition}\label{def:weak-hypercontinuous}
A continuous lattice \(L\) is \emph{weak hypercontinuous} if
\(\sigma(L^{\op})=\upsilon(L^{\op})\).
\end{definition}
This convention is used throughout the paper. The same term has a
different meaning in \cite[Definition~4]{LuoXu2020}.

All spaces are \(T_0\). For a topology \(\tau\), write
\(\Nbh_\tau(x)=\{U\in\tau:x\in U\}\).
The specialization order is defined by
\(x\leq y\) if \(x\in\cl\{y\}\). Open subsets are upper, closed
subsets are lower, and saturated subsets are precisely the upper
subsets. We use local compactness in the following form: for
\(x\in U\in\OO(X)\), there are an open set \(V\) and a compact
saturated set \(K\) such that
\begin{equation}\label{eq:sandwich}
 x\in V\subseteq K\subseteq U.
\end{equation}
For a directed family \(\mathcal D\subseteq\OO(X)\), compactness of
\(K\) and directedness of \(\mathcal D\) give
\[
 U\subseteq\bigcup\mathcal D
 \ \Longrightarrow\ K\subseteq\bigcup\mathcal D
 \ \Longrightarrow\ (\exists D\in\mathcal D)\ V\subseteq K\subseteq D.
\]
Thus \(V\wb U\) in \(\OO(X)\). The sets \(V\wb U\) cover \(U\)
and are directed under finite unions, so \(\OO(X)\) is continuous;
equivalently, \(X\) is core compact.
A sober space is well-filtered: if \(\mathcal K\) is a filtered family
of compact saturated sets and \(\bigcap\mathcal K\subseteq U\) is
open, then some \(K\in\mathcal K\) satisfies \(K\subseteq U\)
\cite[Chapter~II]{GHK03}.

A space \(X=(P,\tau)\) is a \emph{\(d\)-space} if \(P\) is a dcpo
and \(\tau\subseteq\sigma(P)\). Every sober space is a \(d\)-space.
For directed \(D\), the irreducible closed set \(\cl D\) is
\(\dn s\) for a unique \(s\). Then
\[
 D\subseteq\dn s,\qquad
 D\subseteq\dn u\Longrightarrow\dn s=\cl D\subseteq\dn u
 \Longrightarrow s\leq u.
\]
Thus \(s=\bigvee D\). For every open \(U\ni s\),
\(s\in\cl D\) gives \(D\cap U\ne\varnothing\), so \(U\) is Scott open.
A neighborhood \(W\) of \(x\) is a \emph{web at \(x\)} if
\[
 (\forall y\in W)(\exists z\in W)\quad z\leq x\text{ and }z\leq y.
\]
A web space is a space in which every neighborhood of every point
contains a web at that point. We use the established characterization
\begin{equation}\label{eq:web-coframe}
 X\text{ is a web space}
 \quad\Longleftrightarrow\quad
 \OO(X)\text{ is a coframe},
\end{equation}
due to Ern\'e \cite{Erne2009}; see also \cite[p.~23:4]{Erne2019}.

Let \(\GG(X)\) be the lattice of closed subsets of \(X\), including
\(\varnothing\), ordered by inclusion. Its meets are intersections
and its joins are closures of unions. The lower Vietoris topology is
generated by
\[
 \hit U=\{C\in\GG(X):C\cap U\ne\varnothing\},\qquad U\in\OO(X).
\]
It is the upper topology, since
\begin{equation}\label{eq:lower-vietoris}
 \hit U=\GG(X)\setminus\dn_{\GG(X)}(X\setminus U).
\end{equation}
The map \(\eta_X:X\to\GG(X)\), \(\eta_X(x)=\dn x\), satisfies
\(\eta_X^{-1}(\hit U)=U\). Hence it is a topological embedding into
the lower Vietoris topology.

A lattice filter is an upper set containing \(1\) and closed under
finite meets. For a frame \(A\), a proper completely prime filter
is a lattice filter \(F\) such that \(0\notin F\) and
\[
 \bigvee S\in F\quad\Longrightarrow\quad S\cap F\ne\varnothing
 \qquad(S\subseteq A).
\]
The hull--kernel topology \(\tau_A\) on \(\Spec(A)\) consists of
the sets \(\hata a\), \(a\in A\), because
\[
 \hata a\cap\hata b=\widehat{a\wedge b},\qquad
 \bigcup_{a\in S}\hata a=\widehat{\bigvee S}.
\]
A directed union of proper completely prime filters is a proper
completely prime filter and is their supremum in the inclusion order.
Consequently \(\Spec(A)\) is a dcpo and
\(\tau_A\subseteq\sigma(\Spec(A))\).
The Hofmann--Lawson theorem states that, when \(A\) is continuous,
\(\Spec(A)\) is locally compact and sober and
\(a\mapsto\hata a\) is a frame isomorphism
\(A\cong\OO(\Spec(A))\) \cite[Theorem~9.6]{HL78}.

We will also use the prime-element description of the same spectrum.
For a proper completely prime filter \(F\), put
\(p_F=\bigvee(A\setminus F)\). Complete primality gives
\(p_F\notin F\) and \(A\setminus F=\dn p_F\). Moreover,
\[
 a\wedge b\leq p_F\quad\Longrightarrow\quad
 a\leq p_F\text{ or }b\leq p_F.
\]
Conversely, if \(p<1\) has this property, then
\(A\setminus\dn p\) is a proper completely prime filter. These
constructions are inverse, and
we denote the set of such prime elements by \(\Prm(A)\). Thus
\begin{equation}\label{eq:prime-spectrum}
 A\setminus\dn p\subseteq A\setminus\dn q
 \quad\Longleftrightarrow\quad q\leq p.
\end{equation}
In particular, \(\Spec(A)\) is order-isomorphic to
\(\Prm(A)^{\op}\), and \(\hata a\) corresponds to
\(\{p\in\Prm(A):a\not\leq p\}\).
For every complete lattice \(Q\), put
\begin{equation}\label{eq:cop-definition}
 \Pc(Q)=\left\{p\in Q\setminus\{0\}:
 \begin{array}{l}
 (\forall a,b\in Q)\quad p\leq a\vee b\\[-2pt]
 \hfill\Longrightarrow p\leq a\text{ or }p\leq b
 \end{array}\right\}.
\end{equation}
If \(Q^{\op}\) is a frame, its spectrum is order-isomorphic to
\(\Pc(Q)\), and the hull--kernel opens are
\(\{p\in\Pc(Q):p\not\leq a\}\), \(a\in Q\).

For a continuous lattice \(L\), the Lawson topology
\(\lw(L)=\sigma(L)\vee\upsilon(L^{\op})\) is compact Hausdorff,
the order is Lawson closed, and binary meet is jointly Lawson
continuous \cite[Chapter~III]{GHK03}. Products of continuous lattices
are continuous, with Lawson topology equal to the product of the
factor Lawson topologies. A space is coherent if the intersection of
every two compact saturated subsets is compact. A stably compact space
is a compact, locally compact, sober, coherent space. Its patch topology is generated by its opens and the
complements of its compact saturated subsets. This topology is compact
Hausdorff, the specialization order is patch closed, and the original
opens are precisely the patch-open upper sets
\cite[Section~VI-6]{GHK03}.

We record the net-convergence facts used below. Let \((K,\tau)\) be
compact Hausdorff. If a net \((x_i)_{i\in I}\) does not converge to
\(x\), then
\[
 (\exists U\in\Nbh_\tau(x))(\forall i_0\in I)(\exists i\geq i_0)
 \quad x_i\notin U.
\]
There is a subnet in the compact set \(K\setminus U\), and a further
subnet converges to a point of \(K\setminus U\). Consequently,
\begin{equation}\label{eq:unique-subnet-limit}
 \left[\text{every convergent subnet of }(x_i)\text{ has limit }x\right]
 \quad\Longrightarrow\quad x_i\xrightarrow{\tau}x.
\end{equation}

Suppose also that the order on \(K\) is closed. For an increasing net
\((x_i)_{i\in I}\), let \(x_{i(j)}\to x\) be a convergent subnet.
For \(i_0\in I\) and any upper bound \(u\) of \(\{x_i:i\in I\}\),
\[
 x_{i_0}\leq x_{i(j)}\leq u\quad\text{eventually}
 \quad\Longrightarrow\quad x_{i_0}\leq x\leq u.
\]
Thus \(x=\bigvee_{i\in I}x_i\). Compactness and
\eqref{eq:unique-subnet-limit}, and the same argument for the reverse
order, give
\begin{align}
 x_i\leq x_j\ (i\leq j)
 &\quad\Longrightarrow\quad
 x_i\xrightarrow{\tau}\bigvee_{i\in I}x_i,
 \label{eq:increasing-net}\\
 x_j\leq x_i\ (i\leq j)
 &\quad\Longrightarrow\quad
 x_i\xrightarrow{\tau}\bigwedge_{i\in I}x_i.
 \label{eq:decreasing-net}
\end{align}
In particular, the suprema and infima in these formulas exist.
If \(a_i\leq b_i\leq c_i\), \(a_i\to a\), and \(c_i\to a\),
then every convergent subnet \(b_{i(j)}\to b\) satisfies
\[
 (a_{i(j)},b_{i(j)})\to(a,b),\quad
 (b_{i(j)},c_{i(j)})\to(b,a)
 \quad\Longrightarrow\quad a\leq b\leq a.
\]
Equation~\eqref{eq:unique-subnet-limit} therefore yields
\begin{equation}\label{eq:squeeze}
 a_i\leq b_i\leq c_i,\quad a_i\to a,\quad c_i\to a
 \quad\Longrightarrow\quad b_i\to a.
\end{equation}

\section{Embeddings into continuous coframes}\label{sec:host}

We construct a frame embedding into a continuous frame-coframe that
preserves the way-below relation. Write
\(\omega(A)=\upsilon(A^{\op})\) for the lower topology of a lattice
\(A\). Let \(L\) be a continuous frame. For \(a\in L\), define
\(c_a:L\to L\) by \(c_a(x)=a\vee x\), and put
\begin{equation}\label{eq:hostdefinition}
 T=\cl_{L^L}\{c_a:a\in L\},\qquad h(a)=c_a,
\end{equation}
where \(L^L\) carries the product of the Lawson topologies and the
pointwise order. The basic properties of \(T\) require only continuity
and finite distributivity of \(L\). To prove the coframe law, we will
then assume
\begin{equation}\label{eq:dualcondition}
 \sigma(L^{\op})=\upsilon(L^{\op}).
\end{equation}

\subsection{The Lawson closure of join translations}\label{subsec:host-base}

\begin{lemma}\label{lem:hostbase}
The ordered set \(T\) is a continuous lattice. Its meets and
directed joins are pointwise, and its Lawson topology is the subspace
topology inherited from \(L^L\). The map \(h:L\to T\) is injective,
preserves arbitrary joins and finite meets, and satisfies
\[
 a\wb_L b\quad\Longrightarrow\quad h(a)\wb_T h(b).
\]
\end{lemma}
\begin{proof}
Put \(Z=L^L\), with the pointwise order and
\(\lambda(Z)=\prod_{x\in L}\lambda(L)\), and put
\(H=\{c_a:a\in L\}\). Then \(T=\cl_ZH\).
Finite distributivity in \(L\) gives
\[
 (c_a\wedge_Zc_b)(x)
 =(a\vee x)\wedge(b\vee x)=(a\wedge b)\vee x=c_{a\wedge b}(x).
\]
The map \(\mu:Z^2\to Z\), \(\mu(f,g)(x)=f(x)\wedge_Lg(x)\),
is continuous, and \(T\times T=\cl_{Z^2}(H\times H)\). Hence
\[
 \mu[T\times T]
 \subseteq\cl_Z\mu[H\times H]\subseteq\cl_ZH=T,
 \qquad c_1\in T.
\]

For \(S\subseteq T\), index by \((\Fin(S),\subseteq)\) and define
\[
 m_F(x)=\bigwedge_{t\in F}t(x),\qquad
 m_S(x)=\bigwedge_{t\in S}t(x)\quad(x\in L).
\]
Finite-meet closure gives \(m_F\in T\), including
\(m_\varnothing=c_1\). For every \(x\in L\),
\[
 F\subseteq G\ \Longrightarrow\ m_G(x)\leq m_F(x),\qquad
 \bigwedge_{F\in\Fin(S)}m_F(x)=m_S(x).
\]
By \eqref{eq:decreasing-net},
\[
 (\forall x\in L)\quad m_F(x)\xrightarrow{\lambda(L)}m_S(x),
 \qquad m_F\xrightarrow{\lambda(Z)}m_S.
\]
Since \(m_F\in T=\cl_ZT\),
\(m_S\in\cl_Z\{m_F:F\in\Fin(S)\}\subseteq T\).
The pointwise order now gives \(\bigwedge_TS=m_S\).

Let \(D\subseteq T\) be directed, indexed by \((D,\leq)\), and define
\(v_D(x)=\bigvee_{d\in D}d(x)\). For every \(x\in L\),
\[
 d\leq d'\ \Longrightarrow\ d(x)\leq d'(x),\qquad
 d(x)\xrightarrow{\lambda(L)}v_D(x)
\]
by \eqref{eq:increasing-net}. Thus \(d\xrightarrow{\lambda(Z)}v_D\),
so \(v_D\in\cl_ZD\subseteq T\), and \(\bigvee_TD=v_D\).
In particular, \(T\) has all meets, including the empty meet \(c_1\),
and hence is a complete lattice.

The product lattice \(Z\) is continuous. For inclusion \(j:T\to Z\),
define
\[
 p(f)=\bigwedge_T\{t\in T:f\leq_Zj(t)\}.
\]
The defining family contains \(c_1\), and pointwise meets give
\[
 p(f)\leq_Tt\iff f\leq_Zj(t),\qquad
 p\dashv j,\quad pj=\id_T.
\]
For \(t\in T\), put \(D_t=\{u\in Z:u\wb_Zj(t)\}\). Then
\[
 p[D_t]\text{ directed},\qquad
 \bigvee_Tp[D_t]=p\left(\bigvee_ZD_t\right)=pj(t)=t.
\]
If \(u\in D_t\), \(D\subseteq T\) is directed, and
\(t\leq\bigvee_TD\), preservation of directed joins by \(j\) gives
\[
 u\wb_Zj(t)\leq j\left(\bigvee_TD\right)=\bigvee_Zj[D]
 \ \Longrightarrow\ (\exists d\in D)\ u\leq j(d)
 \ \Longrightarrow\ p(u)\leq d.
\]
Hence \(p[D_t]\subseteq\{v:v\wb_Tt\}\).
Conversely, \(v\wb_Tt=\bigvee_Tp[D_t]\) implies
\(v\leq p(u)\) for some \(u\in D_t\).
Thus \(\{v:v\wb_Tt\}\) is directed with supremum \(t\), proving
continuity of \(T\).

Let \(t\in T\), and choose \((a_i)_{i\in I}\subseteq L\) with
\(c_{a_i}\xrightarrow{\lambda(Z)}t\).
Closedness of the order in \(L\) and continuity of meet give
\begin{align*}
 x\leq y&\ \Longrightarrow\ c_{a_i}(x)\leq c_{a_i}(y)
                     \ \Longrightarrow\ t(x)\leq t(y),\\
 x&\leq c_{a_i}(x)\quad(i\in I)
                     \ \Longrightarrow\ x\leq t(x),\\
 t(x\wedge y)
 &=\lim_i c_{a_i}(x\wedge y)
  =\lim_i\bigl(c_{a_i}(x)\wedge c_{a_i}(y)\bigr)
  =t(x)\wedge t(y),\\
 t(1)&=\lim_i c_{a_i}(1)=1.
\end{align*}
For \(r:T\to L\), \(r(t)=t(0)\),
\begin{align*}
 h(a)\leq t&\ \Longrightarrow\ a=c_a(0)\leq t(0)=r(t),\\
 a\leq r(t)&\ \Longrightarrow\
 a\vee x\leq t(0)\vee x\leq t(x)\quad(x\in L).
\end{align*}
Consequently
\begin{equation}\label{eq:host-adjunction}
 h(a)\leq t\iff a\leq r(t),\qquad h\dashv r,\qquad rh=\id_L.
\end{equation}
Thus \(h\) is injective and, as a left adjoint, preserves arbitrary
joins. Moreover,
\[
 h(a\wedge b)=h(a)\wedge_Th(b),\qquad h(1)=c_1=1_T,
\]
so \(h\) preserves finite meets. For directed \(D\subseteq T\),
\[
 r\left(\bigvee_TD\right)
 =\left(\bigvee_TD\right)(0)=\bigvee_{d\in D}d(0)=\bigvee_Lr[D].
\]
If \(a\wb_Lb\) and \(h(b)\leq\bigvee_TD\), then
\[
 b=rh(b)\leq\bigvee_Lr[D]
 \ \Longrightarrow\ (\exists d\in D)\ a\leq r(d)
 \ \Longrightarrow\ h(a)\leq d.
\]
Hence \(h(a)\wb_Th(b)\).

For \(x\in L\), evaluation \(\operatorname{ev}_x(t)=t(x)\) satisfies
\[
 \operatorname{ev}_x\left(\bigwedge_TS\right)
 =\bigwedge_L\operatorname{ev}_x[S],\qquad
 \operatorname{ev}_x\left(\bigvee_TD\right)
 =\bigvee_L\operatorname{ev}_x[D]
\]
for arbitrary \(S\subseteq T\) and directed \(D\subseteq T\).
Thus \(\operatorname{ev}_x\) is Scott continuous and has a left
adjoint \(\ell_x:L\to T\). The equality
\[
 \operatorname{ev}_x^{-1}(L\setminus\up_La)
 =T\setminus\up_T\ell_x(a)
\]
proves Lawson continuity of \(\operatorname{ev}_x\).
Hence the identity
\[
 (T,\lambda(T))\longrightarrow(T,\lambda(Z)|_T)
\]
is a continuous bijection from a compact space to a Hausdorff space,
and therefore a homeomorphism.
\end{proof}

Arbitrary joins in \(T\) are its intrinsic lattice joins; only directed
joins have been shown to be pointwise. Membership in \(T\) does not
assert Scott continuity of \(t:L\to L\). The proof also gives that
every \(t\in T\) is monotone, inflationary and finite-meet-preserving,
and that \(r(t)=t(0)\) is the Scott-continuous right adjoint of \(h\).
In particular, \(0_T=c_0=\id_L\), whereas \(0_{L^L}\) is the
constant-zero map.

\subsection{The coframe law}

For the rest of this section, assume \eqref{eq:dualcondition}.
Every upper subset of \(L\) closed under filtered infima is then
Lawson closed. Put \(Z=L^L\), as in the proof of
Lemma~\ref{lem:hostbase}.

\begin{lemma}\label{ap:sequence}
If \(b_n\in L\), \(n\geq1\), and \(v\leq b_n\vee b_m\) for all \(n\ne m\), then
\begin{equation}\label{ap:eq:liminf}
 v\leq\bigvee_{N\geq1}\bigwedge_{n\geq N} b_n.
\end{equation}
Consequently, if a net \((z_i)_{i\in I}\) converges to \(z\) in \(\lambda(L)\) and, for each \(i\), \(v\leq z_i\vee z_j\) eventually in \(j\), then \(v\leq z\).
\end{lemma}
\begin{proof}
Suppose \eqref{ap:eq:liminf} fails. Choose \(u\ll v\) with \(u\nleq\bigvee_N\bigwedge_{n\geq N}b_n\). Every tail contains an index \(n\) with \(u\nleq b_n\), for otherwise \(u\leq\bigwedge_{n\geq N}b_n\) for some \(N\). Passing to a subsequence and relabelling, we may therefore assume \(u\nleq b_n\) for all \(n\). For \(N\geq1\), let
\[
 E_N=\uparrow v\cup\bigcup_{n\geq N}\uparrow b_n.
\]
Let \(D\subseteq E_N\) be filtered. If \(v\nleq\bigwedge D\), choose \(d_0\in D\) with \(v\nleq d_0\). There is a unique \(n\geq N\) with \(b_n\leq d_0\), by the hypothesis on pairwise joins. For \(d\in D\), take \(e\in D\) with \(e\leq d,d_0\). Then \(v\nleq e\); the same uniqueness gives \(b_n\leq e\leq d\). Thus \(b_n\leq\bigwedge D\). In either case \(\bigwedge D\in E_N\). By \eqref{eq:dualcondition}, \(E_N\) is Lawson closed.

By compactness take \(b\in\bigcap_N\cl_{\lambda(L)}\{b_n:n\geq N\}\). Then \(b\in\bigcap_NE_N\). If \(v\nleq b\), infinitely many \(b_n\) lie below \(b\), which contradicts the condition on pairwise joins. Therefore \(v\leq b\). Now \(u\ll b\), whereas every \(b_n\) lies outside the Scott-open set \(\Uparrow u\). This contradicts the choice of \(b\).

For the net assertion, suppose \(v\nleq z\) and choose \(u\ll v\) with \(u\nleq z\). Eventually \(u\nleq z_i\), since \(\uparrow u\) is Lawson closed. Recursively choose indices \(i_n\) from this tail such that \(v\leq z_{i_m}\vee z_{i_n}\) whenever \(m<n\); each step uses finitely many eventual conditions. By \eqref{ap:eq:liminf}, \(v\leq\bigvee_N\bigwedge_{n\geq N}z_{i_n}\). This forces \(u\leq z_{i_n}\) for all sufficiently large \(n\), a contradiction.
\end{proof}

\begin{lemma}\label{ap:limits}
Under \eqref{eq:dualcondition}, the following hold.
\begin{enumerate}
\item For each \(b\in L\), \(c_b\) is Lawson continuous at \(b\).
\item Every \(t\in T\) is a nucleus on \(L\).
\item If \(c_{a_i}\to t\) in \(T\), then \(t\circ c_{a_i}\to t\) in \(\lambda(Z)\).
\end{enumerate}
\end{lemma}
\begin{proof}
(1) Suppose \(z_i\to b\) but \(b\vee z_i\not\to b\). By compactness, pass to a subnet with \(b\vee z_i\to q\ne b\); then \(b\leq q\). Choose \(v\ll q\) with \(v\nleq b\). On a tail, \(v\ll b\vee z_i\). For each fixed \(i\) in this tail, Scott continuity of \(c_{z_i}\) and convergence \(z_j\to b\) imply \(v\leq z_i\vee z_j\) eventually in \(j\). The net assertion of Lemma~\ref{ap:sequence} gives \(v\leq b\), a contradiction.

(2) For \(c_{a_i}\to t\), fix \(x\in L\) and put \(b=t(x)\). Since \(x\leq b\), (1) gives
\[
 t(b)=\lim_i(a_i\vee b)=\lim_i\big((a_i\vee x)\vee b\big)=b.
\]
Together with Lemma~\ref{lem:hostbase}, this proves that \(t\) is a nucleus.

(3) Fix \(x\in L\), put \(z_i=a_i\vee x\) and \(b=t(x)\). If \(t(z_i)\not\to b\), pass to a subnet with \(t(z_i)\to q\ne b\). Since \(z_i\to b\) and \(z_i\leq t(z_i)\), we have \(b\leq q\). Take \(v\ll q\) with \(v\nleq b\). On a tail, \(v\ll t(z_i)\). For each fixed \(i\) in that tail,
\[
 z_i\vee z_j=c_{a_j}(z_i)\longrightarrow t(z_i),
\]
so \(v\leq z_i\vee z_j\) eventually. Again Lemma~\ref{ap:sequence} gives a contradiction. Therefore \(t(z_i)\to b\) for every \(x\), as required.
\end{proof}

\begin{lemma}\label{ap:join}
For \(a\in L\) and \(t\in T\),
\begin{equation}\label{ap:eq:translation}
 c_a\vee_T t=t\circ c_a.
\end{equation}
For each \(t\in T\), the map \(s\mapsto t\vee_T s\) is Lawson continuous.
\end{lemma}
\begin{proof}
Take \(c_{b_i}\to t\). Then \(c_{a\vee b_i}\to t\circ c_a\) pointwise in the Lawson topology, so \(t\circ c_a\in T\). It is an upper bound of \(t,c_a\). If \(w\in T\) is another upper bound, then by idempotence,
\[
 t(a\vee x)\leq w(a\vee x)\leq w(w(x))=w(x).
\]
This proves \eqref{ap:eq:translation}. Consequently joining with a fixed \(c_a\) is Lawson continuous, by coordinate evaluation, and preserves arbitrary meets.

By Lemma~\ref{ap:limits}(3), whenever \(c_{a_i}\to t\),
\begin{equation}\label{ap:eq:diagonal}
 t\vee_T c_{a_i}\longrightarrow t.
\end{equation}
We extend this to arbitrary translations. Fix \(u,v\in T\), take \(c_{a_i}\to u\), and put \(w=u\vee_T v\). If \(v\vee_T c_{a_i}\not\to w\), there is a subnet on which \(v\vee_T c_{a_i}\to q\ne w\). Closedness of the order gives \(w\leq q\). Choose \(x\in L\) and \(d\ll q(x)\) with \(d\nleq w(x)\). On a tail,
\begin{equation}\label{ap:eq:badjoin}
 d\leq (v\vee_T c_{a_i})(x).
\end{equation}
Take \(c_{b_j}\to w\). Since \(c_{b_j}\leq u\vee_T c_{b_j}\leq w\vee_T c_{b_j}\), \eqref{ap:eq:diagonal} gives \(u\vee_T c_{b_j}\to w\). For any Lawson neighbourhood \(V\) of \(w\), first choose \(j\) such that \(u\vee_T c_{b_j}\in V\), and then choose \(i\) on the tail in \eqref{ap:eq:badjoin} with \(c_{a_i\vee b_j}\in V\), using continuity of joining with \(c_{b_j}\). For \(a=a_i\vee b_j\),
\[
 c_a\in V,\qquad d\leq(v\vee_T c_{a_i})(x)\leq(w\vee_T c_a)(x).
\]
Thus \(w\) is in the closure of \(\{c_a:d\leq(w\vee_T c_a)(x)\}\). A net in that set converging to \(w\), together with \eqref{ap:eq:diagonal} and closedness of \(\uparrow d\), implies \(d\leq w(x)\), a contradiction.

We have proved \(v\vee_T c_{a_i}\to v\vee_T u\) for every net of join translations converging to \(u\). The closure in \(T^2\) of \(\{(c_a,v\vee_T c_a):a\in L\}\) is therefore precisely the graph of \(u\mapsto v\vee_T u\): density gives one inclusion and the convergence just proved gives the other. The projection of this graph onto the first factor is a continuous bijection from a compact space to the Hausdorff space \(T\), hence a homeomorphism. Composing its inverse with the second projection proves continuity.
\end{proof}

\begin{theorem}\label{thm:embedding}
Let \(L\) be a continuous frame with
\(\sigma(L^{\op})=\upsilon(L^{\op})\). There exist a continuous
lattice \(T\) that is both a frame and a coframe and an injective
frame homomorphism \(h:L\to T\) preserving the way-below relation.
\end{theorem}
\begin{proof}
By Lemma~\ref{ap:join}, for \(a\in L\) and \(s,t\in T\),
\[
 c_a\vee_T(s\wedge t)=(c_a\vee_T s)\wedge(c_a\vee_T t).
\]
By continuity of the three join translations and of meet, this equality extends from \(h[L]\) to \(T\). Hence \(T\) is distributive. Since \(T\) is continuous, it is a frame. For \(u\in T\) and a filtered subset \(D\subseteq T\), continuity of joining with \(u\) and convergence of decreasing nets give
\[
 u\vee_T\bigwedge D=\bigwedge_{d\in D}(u\vee_T d).
\]
Every arbitrary meet is the filtered infimum of the finite meets. Finite distributivity therefore yields the coframe law. The assertions about \(h\) and \(r\) follow from Lemma~\ref{lem:hostbase}.
\end{proof}

\section{A topological criterion for spectra}\label{sec:coframes}

We prove a Scott-topology criterion by constructing nested compact
saturated neighbourhoods. It applies to the spectrum of every
continuous frame-coframe and hence to every locally compact sober
web space.

\subsection{Compact saturated neighbourhoods}\label{ap:topology}

For a space \(X\) and \(U\in\OO(X)\), put
\begin{equation}\label{ap:eq:R}
 R_U=\{(y,z)\in X^2:U\cap\downarrow y\cap\downarrow z\ne\varnothing\},
 \qquad H_U(y)=\{z:(y,z)\in R_U\}.
\end{equation}
Assume that each \(R_U\) is open in \(X^2\). The next two lemmas give
first-countable points below any given point. A triangular construction,
related to the compact selection argument in \cite[Section~5]{He2026},
then proves equality with the Scott topology.

\begin{lemma}\label{ap:principal}
Let \(X\) be a locally compact sober space satisfying the condition on \(R_U\). Suppose \((K_n)_{n\geq1}\) is a sequence of nonempty compact saturated sets with
\begin{equation}\label{ap:eq:principal}
 K_{n+1}\subseteq\intr K_n,\qquad K_{n+1}^2\subseteq R_{\intr K_n}.
\end{equation}
Then \(\bigcap_n K_n=\uparrow m\) for a first-countable point \(m\) of \(X\).
\end{lemma}
\begin{proof}
By well-filteredness, \(K=\bigcap_nK_n\) is nonempty and compact saturated. For \(y,z\in K\), each \(K_n\) meets the closed set \(\downarrow y\cap\downarrow z\), by the second inclusion in \eqref{ap:eq:principal}. Well-filteredness, applied to its open complement, implies that \(K\) meets this closed set. Thus \(K\) is filtered. The closed subsets \(K\cap\downarrow y\), \(y\in K\), have the finite intersection property in \(K\), so \(K\) has a least element \(m\). Saturation gives \(K=\uparrow m\). Since \(m\in\intr K_n\) for every \(n\), and well-filteredness gives \(K_n\subseteq U\) for some \(n\) whenever \(m\in U\in\OO(X)\), the sets \(\intr K_n\) form a countable local base at \(m\).
\end{proof}

\begin{lemma}\label{ap:firstcount}
Let \(X\) satisfy the assumptions of Lemma~\ref{ap:principal}. Let \(G\) be its set of first-countable points. For \(x\in X\) and any countable family \((V_n)_{n\geq1}\) of open neighbourhoods of \(x\), there is an \(m\in G\) with
\[
 m\leq x,\qquad m\in\bigcap_n V_n.
\]
In particular, \(G\cap\downarrow x\) is directed and has supremum \(x\).
\end{lemma}
\begin{proof}
Choose a compact saturated neighbourhood \(K_1\subseteq V_1\) of \(x\). Given \(K_n\), the relation \(R_{\intr K_n}\) is open and contains \((x,x)\). Hence there is an open neighbourhood \(W\) of \(x\) with \(W\subseteq V_{n+1}\cap\intr K_n\) and \(W^2\subseteq R_{\intr K_n}\). Choose a compact saturated neighbourhood \(K_{n+1}\) of \(x\) contained in \(W\). By Lemma~\ref{ap:principal}, \(\bigcap_nK_n=\uparrow m\) for some \(m\in G\). The construction gives the desired inequalities and memberships.

For \(a,b\in G\cap\downarrow x\), apply the assertion to the combined countable local bases of \(a,b\), all of whose members contain \(x\). The resulting \(m\) satisfies \(a,b\leq m\leq x\). The set \(G\cap\downarrow x\) is nonempty and meets each neighbourhood of \(x\), by the assertion. Since \(X\) is a \(d\)-space, it has a supremum \(s\leq x\). Every neighbourhood of \(x\) contains \(s\), so \(s=x\).
\end{proof}

\begin{proposition}\label{ap:topcriterion}
Let \(X=(P,\tau)\) be a locally compact sober space. If \(R_U\) is open in \(X^2\) for every \(U\in\tau\), then \(\tau=\sigma(P)\).
\end{proposition}
\begin{proof}
Let \(B\) be Scott closed in \(P\). For \(x\in\cl_X B\) and any open neighbourhood \(U\) of \(x\), the set \(H_U(x)\) is an open neighbourhood of \(x\). Choose \(b\in B\cap H_U(x)\) and then \(c\in U\cap\downarrow b\cap\downarrow x\). By lowerness of \(B\),
\begin{equation}\label{ap:eq:lowerhit}
 B\cap\downarrow x\cap U\ne\varnothing.
\end{equation}
First suppose that \(x\) has a decreasing countable local base \((U_n)_{n\geq1}\). We construct points \(c_{n,k}\in B\cap\downarrow x\) and compact saturated sets \(K_{n,k}\) for \(1\leq n\leq k\), with the following properties whenever the indices are defined:
\begin{align}
 c_{n,k}&\in\intr K_{n,k}\subseteq K_{n,k}\subseteq U_n,\label{ap:eq:array1}\\
 c_{n,k}&\leq c_{n+1,k},\qquad c_{n,k+1}\leq c_{n,k},\label{ap:eq:array2}\\
 K_{n+1,k}&\subseteq\intr K_{n,k}\cap H_{\intr K_{n,k}}(c_{n,k}),\label{ap:eq:array3}\\
 K_{n,k+1}&\subseteq\intr K_{n,k},\qquad K_{n,k+1}^2\subseteq R_{\intr K_{n,k}}.\label{ap:eq:array4}
\end{align}
Choose \(c_{1,1}\) by \eqref{ap:eq:lowerhit}, and choose its compact saturated neighbourhood \(K_{1,1}\subseteq U_1\). Suppose column \(k\) has been constructed. We choose the points of column \(k+1\) in reverse row order. Since \(c_{k,k}\leq x\), the set \(H_{\intr K_{k,k}}(c_{k,k})\) is an open neighbourhood of \(x\). By \eqref{ap:eq:lowerhit}, choose
\[
 c_{k+1,k+1}\in B\cap\downarrow x\cap U_{k+1}\cap H_{\intr K_{k,k}}(c_{k,k}).
\]
For \(n=k,k-1,\ldots,1\), choose
\begin{equation}\label{ap:eq:choosepoint}
 c_{n,k+1}\in\intr K_{n,k}\cap\downarrow c_{n,k}\cap\downarrow c_{n+1,k+1}.
\end{equation}
When \(n=k\), the intersection is nonempty by the choice of \(c_{k+1,k+1}\). For \(n<k\), the already chosen point \(c_{n+1,k+1}\) lies in \(\intr K_{n+1,k}\), which by \eqref{ap:eq:array3} is contained in \(H_{\intr K_{n,k}}(c_{n,k})\). This proves that \eqref{ap:eq:choosepoint} is possible. The new points lie in \(B\cap\downarrow x\) and satisfy \eqref{ap:eq:array2}.

Next choose their compact saturated neighbourhoods in the row order \(n=1,\ldots,k+1\). If \(n\leq k\), choose an open neighbourhood \(W\) of \(c_{n,k+1}\) in \(\intr K_{n,k}\) such that \(W^2\subseteq R_{\intr K_{n,k}}\); such a neighbourhood exists by openness of the relation and \eqref{ap:eq:choosepoint}. If \(n=k+1\), start with \(W=U_{k+1}\). In either case, for \(n>1\), intersect \(W\) with
\[
 \intr K_{n-1,k+1}\cap H_{\intr K_{n-1,k+1}}(c_{n-1,k+1}).
\]
This remains a neighbourhood of \(c_{n,k+1}\), by \(c_{n-1,k+1}\leq c_{n,k+1}\). Choose \(K_{n,k+1}\subseteq W\) with \(c_{n,k+1}\in\intr K_{n,k+1}\). All of \eqref{ap:eq:array1}--\eqref{ap:eq:array4} now hold. This completes the recursion.

For a fixed row \(n\), Lemma~\ref{ap:principal} applied to \eqref{ap:eq:array4} gives
\[
 \bigcap_{k\geq n} K_{n,k}=\uparrow a_n
\]
for some \(a_n\in X\). By \eqref{ap:eq:array2}, \(c_{n,k}\leq c_{n,n}\) for every \(k\geq n\). As the sets \(K_{n,k}\) are upper, \(c_{n,n}\in\bigcap_{k\geq n}K_{n,k}\); hence \(a_n\leq c_{n,n}\in B\), and \(a_n\in B\). The vertical inclusions \eqref{ap:eq:array3} give \(a_n\leq a_{n+1}\). Moreover, each \(K_{n,k}\) contains \(x\) since \(c_{n,k}\leq x\). By \eqref{ap:eq:array1}, therefore,
\[
 a_n\leq x,\qquad a_n\in U_n.
\]
Consequently \(\bigvee_na_n=x\), because \(\bigcap_nU_n=\uparrow x\). Scott closedness of \(B\) gives \(x\in B\).

Finally let \(x\in\cl_X B\) be arbitrary. As \(\cl_X B\) is lower, every first-countable point below \(x\) belongs to \(\cl_X B\) and hence to \(B\), by the first part of the proof. By Lemma~\ref{ap:firstcount}, these points form a directed set with supremum \(x\). Thus \(x\in B\). This proves that every Scott-closed set is \(\tau\)-closed. The reverse inclusion follows from sobriety, which implies that \(X\) is a \(d\)-space.
\end{proof}

\subsection{Spectra of continuous frame-coframes}

We use the following consequence of \cite[Corollary~VII-2.4]{GHK03}.

\begin{lemma}\label{ap:joint}
Binary join is jointly Lawson continuous on every continuous coframe.
\end{lemma}

We also use \cite[Proposition~III-3.18(ii)]{GHK03} in the following form.

\begin{lemma}\label{un:lawson-lower}
On a continuous lattice, the Lawson-open lower sets are exactly the
lower-open sets.
\end{lemma}

Use the prime-element representation \eqref{eq:prime-spectrum}.
The correspondence \(p\mapsto A\setminus\downarrow p\)
identifies \(\Spec(A)\) with \(\Prm(A)^{\op}\), and
\(\widehat a\) is represented by
\[
 U_a=\{p\in\Prm(A):a\nleq p\},\qquad a\in A.
\]
Directed suprema in this representation are filtered infima in \(A\).
When \(A\) is continuous, the Hofmann--Lawson representation recalled
in Section~\ref{sec:prelim} gives
\begin{equation}\label{ap:eq:spatial}
 a\nleq b\quad\Longrightarrow\quad
 \exists p\in\Prm(A),\quad b\leq p,\quad a\nleq p.
\end{equation}
Indeed, \(U_a\nsubseteq U_b\) under the isomorphism
\(A\cong\OO(\Spec(A))\).

\begin{corollary}\label{cor:coframe}
If \(A\) is a continuous frame and a coframe, then the hull--kernel
topology on \(\Spec(A)\) is the Scott topology of the inclusion order.
\end{corollary}
\begin{proof}
Let \(M\) be a continuous lattice that is both a frame and a coframe.
Lemma~\ref{ap:joint} shows that join is jointly Lawson continuous on \(M\). By Lemma~\ref{un:lawson-lower}, the lower topology
\(\omega(M)\) consists exactly of the Lawson-open lower sets. The same assertion holds for \(M^2\), whose lower topology is the product
of the lower topologies of \(M\), by \cite[Lemma~III-1.3]{GHK03}.

For \(a\in M\), the set
\[
 \{(p,q)\in M^2:a\nleq p\vee q\}
\]
is Lawson open and lower, so it is open for the product lower topology. Put \(X=\Spec(M)\) and identify its points with \(\Prm(M)^{\op}\). Its hull--kernel topology is the restriction of \(\omega(M)\). By \eqref{ap:eq:spatial},
\[
 \begin{split}
 (p,q)\in R_{U_a}&\Longleftrightarrow
 \exists s\in\Prm(M),\quad p\vee q\leq s,\ a\nleq s\\
 &\Longleftrightarrow a\nleq p\vee q.
 \end{split}
\]
Thus \(R_{U_a}\) is open in \(X^2\). Every open of \(X\) has the form \(U_a\). By the Hofmann--Lawson theorem, \(X\) is locally compact and sober, so Proposition~\ref{ap:topcriterion} applies.
\end{proof}

\begin{theorem}\label{thm:web}
Every locally compact sober web space is the Scott space of its
specialization order.
\end{theorem}
\begin{proof}
Let \(X\) be locally compact, sober and a web space. Then
\(L=\OO(X)\) is a continuous frame and, by
\eqref{eq:web-coframe}, a coframe. Sobriety identifies \(X\) with
\(\Spec(L)\). Corollary~\ref{cor:coframe} therefore gives the
asserted equality of topologies.
\end{proof}

\begin{lemma}\label{lem:fc-core}
Let \(X\) be a locally compact sober web space. For every point \(x\)
and every countable family of neighborhoods \((U_n)_{n\geq1}\) of
\(x\), there is a first-countable point \(m\leq x\) with
\(m\in\bigcap_nU_n\). Moreover, if \(G\) denotes the set of
first-countable points of \(X\), then \(G\cap\dn x\) is directed
and has supremum \(x\).
\end{lemma}
\begin{proof}
Put \(L=\OO(X)\). As in the proof of Corollary~\ref{cor:coframe},
the representation \(X\cong\Spec(L)\) makes \(R_U\) open in
\(X^2\) for every \(U\in\OO(X)\). Apply
Lemma~\ref{ap:firstcount} to the interiors of the prescribed
neighbourhoods.
\end{proof}

\section{The spectrum characterization}\label{sec:transfer}

We prove that equality of the hull--kernel and Scott topologies passes
along an injective frame homomorphism preserving the way-below relation.
Applying this to the embedding in Theorem~\ref{thm:embedding} completes
the characterization.

We use the Scott-open filter separation theorem
\cite[Lemma~1.5]{JM08}: every Scott-open lattice filter \(H\)
of a frame \(A\) is the intersection of the completely prime
filters containing it. Equivalently,
\begin{equation}\label{eq:filter-separation}
 v\notin H
 \quad\Longrightarrow\quad
 (\exists G\in\Spec(A))\ H\subseteq G,\quad v\notin G.
\end{equation}

The adjoint criterion \cite[Theorem~IV-1.4]{GHK03} gives the following
equivalent formulation of way-below preservation.

\begin{lemma}\label{ap:adjoint-waybelow}
Let \(L\) be a continuous dcpo, let \(M\) be a dcpo, and let
\(h:L\to M\) and \(r:M\to L\) satisfy \(h\dashv r\).
Then \(h\) preserves the way-below relation if and only if
\(r\) is Scott continuous.
\end{lemma}

We use the prime-element form of \eqref{eq:filter-separation}.

\begin{lemma}\label{ap:prime}
Let \(h:L\to M\) be a frame homomorphism with Scott-continuous right adjoint \(r\). If \(p\in\Prm(L)\) and \(t\in M\), then
\begin{equation}\label{ap:eq:prime}
 r(t)\leq p\quad\Longleftrightarrow\quad
 \exists q\in\Prm(M),\quad t\leq q,\ r(q)\leq p.
\end{equation}
\end{lemma}
\begin{proof}
Assume \(r(t)\leq p\). Since \(r\) preserves finite meets and directed
joins,
\[
 F=r^{-1}(L\setminus\downarrow p)
\]
is a Scott-open filter of \(M\), and \(t\notin F\). In the
prime-element representation, \eqref{eq:filter-separation} gives
\[
 M\setminus F=\downarrow_M\bigl(\Prm(M)\setminus F\bigr).
\]
Hence there is a \(q\in\Prm(M)\setminus F\) with \(t\leq q\), which
means \(r(q)\leq p\). The converse follows by monotonicity of \(r\).
\end{proof}

\begin{theorem}\label{thm:transfer}
Let \(e:L\to A\) be an injective frame homomorphism, where \(L\)
is continuous, and suppose that \(e\) preserves the way-below
relation. If the hull--kernel topology on \(\Spec(A)\) is its
Scott topology, then the same holds for \(\Spec(L)\).
\end{theorem}
\begin{proof}
Let \(r:A\to L\) be the right adjoint of \(e\). The hypotheses and Lemma~\ref{ap:adjoint-waybelow}
make \(r\) Scott continuous. The induced spectral map is
\[
 f:\Prm(A)^{\op}\longrightarrow\Prm(L)^{\op},\qquad f(q)=r(q).
\]
Since \(e\) is a frame homomorphism, \(f\) is well defined. Preservation
of arbitrary meets by \(r\) makes \(f\) Scott continuous. Also
\(re=\id_L\). For \(p\in\Prm(L)\), Lemma~\ref{ap:prime}, applied to
\(t=e(p)\), gives a prime \(q\) with \(e(p)\leq q\) and \(r(q)\leq p\).
The adjunction forces \(r(q)=p\), so \(f\) is onto.

Let \(B\) be a Scott-closed subset of \(\Prm(L)^{\op}\). By Scott continuity of \(f\) and the hypothesis on \(\Spec(A)\), for some \(t\in A\) we have
\[
 f^{-1}(B)=\Prm(A)\cap\uparrow_A t.
\]
Since \(B\) is upper in the original lattice order of \(L\), surjectivity and \eqref{ap:eq:prime} give, for \(p\in\Prm(L)\),
\[
 \begin{split}
 p\in B&\Longleftrightarrow
 \exists q\in\Prm(A),\quad t\leq q,\ r(q)\leq p\\
 &\Longleftrightarrow r(t)\leq p.
 \end{split}
\]
Hence \(B=\Prm(L)\cap\uparrow_L r(t)\), which is hull--kernel closed. The opposite inclusion of topologies is valid for every frame.
\end{proof}

\begin{proof}[Proof of Theorem~\ref{thm:main}]
The implication \textup{(i)}\(\Rightarrow\)\textup{(ii)} is
\cite[Theorem~4.15, \textup{(4)}\(\Rightarrow\)\textup{(2)}]{CK23}.
Assume \textup{(ii)}. If \(L\) has one element, then
\(\Spec(L)=\varnothing\) and both conditions hold.
Otherwise, \(L\) is a continuous frame. Theorem~\ref{thm:embedding}
gives a continuous frame and coframe \(T\), and an injective frame
homomorphism \(h:L\to T\) with
\[
 a\wb_L b\quad\Longrightarrow\quad h(a)\wb_T h(b).
\]
Corollary~\ref{cor:coframe} gives
\(\tau_T=\sigma(\Spec(T))\).
Applying Theorem~\ref{thm:transfer} to \(h\) yields
\(\tau_L=\sigma(\Spec(L))\), proving \textup{(i)}.
\end{proof}

\appendix
\renewcommand{\thesection}{\Alph{section}}
\makeatletter
\renewcommand{\@seccntformat}[1]{%
  \ifstrequal{#1}{section}{\appendixname~}{}%
  \csname the#1\endcsname.\hskip 0.5em}
\makeatother

\section{An alternative proof of the spectrum characterization}\label{app:spectrum-proof}
\label{ap:short-proof}

We retain a second proof of the three steps in
Theorem~\ref{thm:main}. The topological criterion follows from a
countably based argument and metrizable lattice quotients; the coframe
law follows from simultaneous approximation by join translations;
and spectral transfer follows directly from completely prime filters.
This proof uses the preliminaries and Lemma~\ref{lem:hostbase}, but
not the sequence or triangular-neighbourhood arguments of
Sections~\ref{sec:host} and~\ref{sec:coframes}.

\subsection{Upper and patch topologies}\label{sec:upper-topology}

\Needspace{7\baselineskip}
\begin{lemma}\label{lem:upper-topology}
For every complete lattice \(M\), the space \((M,\upsilon(M))\)
is compact, sober, and coherent.
\end{lemma}
\begin{proof}
Put \(\tau=\upsilon(M)\). Every member of \(\tau\) is upper, and
\[
 x\not\leq y\quad\Longrightarrow\quad
 x\in M\setminus\dn y\in\tau,\qquad y\notin M\setminus\dn y.
\]
Hence the specialization order of \(\tau\) is \(\leq_M\).
For \(U\in\tau\), \(0\in U\) implies \(U=M\). Every open cover
of \(M\) therefore contains \(M\), so \((M,\tau)\) is compact.
Sobriety follows from \cite[Lemma~4.17]{CK23}.

Define \(v:M^2\to M\) by \(v(x,y)=x\vee y\). For every \(b\in M\),
\[
 v^{-1}(M\setminus\dn b)
 =((M\setminus\dn b)\times M)\cup(M\times(M\setminus\dn b))
 \in\tau\times\tau.
\]
Thus \(v\) is continuous. If \(K,H\subseteq M\) are compact saturated,
then \(K=\up K\), \(H=\up H\), and
\[
 z\in v[K\times H]\ \Longrightarrow\ z\in K\cap H,
 \qquad
 z\in K\cap H\ \Longrightarrow\ z=v(z,z)\in v[K\times H].
\]
Consequently \(K\cap H=v[K\times H]\) is compact. Thus intersections
of two compact saturated subsets are compact, and \((M,\tau)\) is coherent.
\end{proof}

We use two standard results on compact Hausdorff semilattices.
A compact Hausdorff semitopological semilattice is a topological
semilattice \cite{Lawson74}; see also the statement preceding
\cite[Theorem~7.1]{BB2019}. For a compact Hausdorff topological meet-semilattice,
a lower subset is closed if and only if it is closed under directed
suprema \cite[Theorem~13]{Lawson73}. In the latter reference,
\emph{ideal} means a lower subset, without a directedness requirement.

\Needspace{7\baselineskip}
\begin{proposition}\label{prop:upper-patch}
Let \(M\) be a frame for which \(\tau=\upsilon(M)\) is locally
compact. Then \(\tau=\sigma(M)\), and meet is jointly continuous
in the patch topology \(\pi\) of \((M,\tau)\). If \(M\) is
also a coframe, then \((M,\pi)\) is a compact Hausdorff topological
lattice, and both lattice operations are jointly \(\tau\)-continuous.
\end{proposition}
\begin{proof}
Lemma~\ref{lem:upper-topology} and local compactness make \((M,\tau)\)
stably compact. The patch correspondence gives
\[
 (M,\pi)\text{ compact Hausdorff},\qquad
 \{(x,y):x\leq y\}\text{ closed in }(M,\pi)^2,
\]
and
\begin{equation}\label{eq:patch-upper}
 \tau=\{U\in\pi:U=\up U\}.
\end{equation}
Let \(\mathcal K\) denote the compact saturated subsets of \((M,\tau)\).
The topology \(\pi\) is generated by
\(\tau\cup\{M\setminus K:K\in\mathcal K\}\).

Fix \(a\in M\). The frame law gives a right adjoint \(r_a\) to
\(m_a(x)=a\wedge x\). Thus
\[
 a\wedge x\leq b\iff x\leq r_a(b),\qquad
 m_a^{-1}(M\setminus\dn b)=M\setminus\dn r_a(b)\in\tau.
\]
For \(K\in\mathcal K\), upperness of \(K\) gives
\begin{align*}
 x\in m_a^{-1}(K)
 &\iff a\wedge x\in K\\
 &\iff (\exists k\in K)\ k\leq a,\ k\leq x\\
 &\iff x\in\up(K\cap\dn a).
\end{align*}
Hence
\begin{equation}\label{eq:meet-compact-preimage}
 m_a^{-1}(K)=\up(K\cap\dn a).
\end{equation}
The set \(\dn a\) is \(\tau\)-closed, so \(K\cap\dn a\) is compact.
Saturation preserves compactness; therefore
\[
 m_a^{-1}(K)\in\mathcal K,\qquad
 m_a^{-1}(M\setminus K)=M\setminus\up(K\cap\dn a)\in\pi.
\]
Thus \(m_a:(M,\pi)\to(M,\pi)\) is continuous for every \(a\).
The compact semilattice joint-continuity theorem \cite{Lawson74}
implies that \(\wedge:(M,\pi)^2\to(M,\pi)\) is continuous.

For \(C\subseteq M\) with \(C=\dn C\), Lawson's theorem
\cite[Theorem~13]{Lawson73} gives
\[
 C\text{ is \(\pi\)-closed}
 \iff
 (\forall D\subseteq C\text{ directed})\ \bigvee D\in C.
\]
By \eqref{eq:patch-upper}, for arbitrary \(C\subseteq M\),
\begin{align*}
 C\text{ is \(\tau\)-closed}
 &\iff C=\dn C\text{ and }C\text{ is \(\pi\)-closed}\\
 &\iff C=\dn C\text{ and }
       (\forall D\subseteq C\text{ directed})\ \bigvee D\in C\\
 &\iff C\text{ is Scott closed}.
\end{align*}
Therefore \(\tau=\sigma(M)\).

Suppose that \(M\) is also a coframe. For fixed \(a\), the map
\(j_a(x)=a\vee x\) preserves arbitrary meets, including the empty
meet. Define
\[
 \ell_a(y)=\bigwedge\{x\in M:y\leq a\vee x\}.
\]
The adjunction and the subbasic inverse images are
\begin{gather*}
 \ell_a(y)\leq x\iff y\leq j_a(x),\\
 j_a^{-1}(M\setminus\dn b)=
 \begin{cases}
 M\setminus\dn b,&a\leq b,\\
 M,&a\not\leq b,
 \end{cases}
 \qquad
 \ell_a^{-1}(M\setminus\dn b)=M\setminus\dn(a\vee b).
\end{gather*}
Both \(j_a\) and \(\ell_a\) are \(\tau\)-continuous. For
\(K\in\mathcal K\),
\begin{align*}
 x\in j_a^{-1}(K)
 &\iff j_a(x)\in K\\
 &\iff (\exists k\in K)\ k\leq j_a(x)\\
 &\iff (\exists k\in K)\ \ell_a(k)\leq x.
\end{align*}
Consequently
\begin{equation}\label{eq:join-compact-preimage}
 j_a^{-1}(K)=\up\ell_a[K]\in\mathcal K,
\end{equation}
since \(\ell_a[K]\) is compact. The inverse images under \(j_a\)
of both types of patch generators belong to \(\pi\); hence \(j_a\)
is \(\pi\)-continuous. The joint-continuity theorem, applied to
\((M,\pi,\vee)\), gives joint \(\pi\)-continuity of \(\vee\).
The natural meet-semilattice order associated with \(\vee\) is
\(\geq_M\), which is also closed.

We verify
\begin{equation}\label{eq:patch-product}
 \tau\times\tau
 =\{O\in\pi\times\pi:O=\up_{M^2}O\}.
\end{equation}
The inclusion from left to right follows from
\eqref{eq:patch-upper}. Conversely, let \(O\in\pi\times\pi\) be
upper and \((x,y)\in O\). Then
\(\up x\times\up y\subseteq O\). Both principal upper sets are
compact. The tube lemma gives \(A,B\in\pi\) with
\[
 \up x\subseteq A,\qquad \up y\subseteq B,\qquad A\times B\subseteq O.
\]
For every \(\pi\)-closed \(F\),
\[
 \dn F=\operatorname{pr}_1
 \{(z,f)\in M\times F:z\leq f\}
\]
is compact and hence \(\pi\)-closed. Therefore
\[
 U=M\setminus\dn(M\setminus A),\qquad
 V=M\setminus\dn(M\setminus B)
\]
are \(\pi\)-open upper sets satisfying
\(x\in U\subseteq A\), \(y\in V\subseteq B\).
By \eqref{eq:patch-upper}, \(U,V\in\tau\), and
\((x,y)\in U\times V\subseteq O\), proving
\eqref{eq:patch-product}.
For \(*\in\{\wedge,\vee\}\) and \(W\in\tau\), joint
\(\pi\)-continuity and monotonicity give
\[
 *^{-1}(W)\in\pi\times\pi,\qquad
 \up_{M^2}*^{-1}(W)=*^{-1}(W).
\]
Equation~\eqref{eq:patch-product} gives \(*^{-1}(W)\in\tau\times\tau\),
proving joint \(\tau\)-continuity of both operations.
\end{proof}

\Needspace{7\baselineskip}
\begin{lemma}\label{lem:lower-vietoris-locally-compact}
If \(X\) is locally compact, then the lower Vietoris topology on
\(\GG(X)\) is locally compact.
\end{lemma}
\begin{proof}
Put \(M=\GG(X)\) and \(\tau=\upsilon(M)\). For compact
\(K\subseteq X\), define
\[
 H(K)=\{C\in M:C\cap K\ne\varnothing\}.
\]
Since \(C=\dn C\) for \(C\in M\),
\[
 C\in H(K)
 \iff (\exists x\in K)\ \dn x\subseteq C
 \iff C\in\up_M\eta_X[K].
\]
The map \(\eta_X:X\to(M,\tau)\) is continuous; hence
\(H(K)=\up_M\eta_X[K]\) is compact saturated. By
Lemma~\ref{lem:upper-topology}, \(\bigcap_{i=1}^nH(K_i)\) is compact
saturated for every finite family of compact subsets \(K_i\subseteq X\).

Let \(C\in\bigcap_{i=1}^n\hit{U_i}\), where \(n\geq1\) and
\(U_i\in\OO(X)\). Choose
\[
 x_i\in C\cap U_i,\qquad
 x_i\in V_i\subseteq K_i\subseteq U_i,
\]
with \(V_i\) open and \(K_i\) compact saturated. Then
\[
 C\in\bigcap_{i=1}^n\hit{V_i}
 \subseteq\bigcap_{i=1}^nH(K_i)
 \subseteq\bigcap_{i=1}^n\hit{U_i}.
\]
The first intersection is \(\tau\)-open and the middle intersection
is compact saturated. For \(n=0\), take both neighborhoods to be
\(M\), which is compact by Lemma~\ref{lem:upper-topology}.
Thus every point and basic open neighborhood admit a compact
neighborhood as in \eqref{eq:sandwich}; \((M,\tau)\) is locally compact.
\end{proof}

\Needspace{7\baselineskip}
\begin{proposition}\label{prop:canonical}
Let \(X\) be a locally compact sober web space and put
\(M=\GG(X)\). Then
\[
 \tau_M=\sigma(M)=\upsilon(M)
\]
is stably compact. Its patch topology \(\pi_M\) makes \(M\) a
compact Hausdorff topological lattice with operations union and
intersection, and both operations are jointly \(\tau_M\)-continuous.
The map \(\eta_X:X\to\Sigma M\) is a topological and order
embedding, and
\begin{equation}\label{eq:eta-saturation}
 \up_M\eta_X(U)=\hit U\qquad(U\in\OO(X)).
\end{equation}
Moreover, \(\OO(X)\) is a continuous frame and a coframe, with
\(\sigma(\OO(X)^{\op})=\upsilon(\OO(X)^{\op})\).
\end{proposition}
\begin{proof}
By \eqref{eq:web-coframe}, \(\OO(X)\) is a coframe. The map
\[
 \OO(X)^{\op}\longrightarrow M,\qquad U\longmapsto X\setminus U
\]
is an order isomorphism, so \(M\) is a frame. For \(C\in M\) and
\(\mathcal A\subseteq M\),
\[
 C\cup\bigcap_{A\in\mathcal A}A
 =\bigcap_{A\in\mathcal A}(C\cup A),
\]
where both sides equal \(X\) when \(\mathcal A=\varnothing\).
Thus \(M\) is also a coframe.
Lemma~\ref{lem:lower-vietoris-locally-compact} gives local compactness
of \(\upsilon(M)\). Lemma~\ref{lem:upper-topology} and
Proposition~\ref{prop:upper-patch} now give
\[
 \tau_M=\sigma(M)=\upsilon(M),\qquad
 (M,\tau_M)\text{ stably compact},
\]
and joint continuity of union and intersection for both \(\tau_M\)
and its compact Hausdorff patch topology \(\pi_M\).

For \(x,y\in X\) and \(U\in\OO(X)\),
\[
 \dn x\subseteq\dn y\iff x\leq y,\qquad
 \eta_X^{-1}(\hit U)=U.
\]
Since the sets \(\hit U\) generate \(\tau_M\), the map
\(\eta_X:X\to\Sigma M\) is a topological and order embedding.
Moreover,
\begin{align*}
 C\in\up_M\eta_X[U]
 &\iff (\exists x\in U)\ \dn x\subseteq C\\
 &\iff C\cap U\ne\varnothing
 \iff C\in\hit U,
\end{align*}
which proves \eqref{eq:eta-saturation}.
Local compactness of \(X\) gives continuity of \(\OO(X)\).
Finally, the order isomorphism \(\OO(X)^{\op}\cong M\) and
\(\sigma(M)=\upsilon(M)\) imply
\(\sigma(\OO(X)^{\op})=\upsilon(\OO(X)^{\op})\).
\end{proof}

The topology \(\pi_M\) is the patch topology of
\((M,\sigma(M))\). No continuity of \(M\), or identification of
\(\pi_M\) with a Lawson topology, is required.

\Needspace{10\baselineskip}
\subsection{Countably based spaces}\label{sec:countable}

\Needspace{7\baselineskip}
\begin{lemma}\label{lem:fusion}
Let \(X\) be core compact, and suppose that
\(\eta_X:X\to\Sigma\GG(X)\) is continuous. For every decreasing
sequence \((A_n)_{n\geq0}\) of lower subsets of \(X\),
\begin{equation}\label{eq:fusion}
 \cl_X\!\left(\bigcap_{n\geq0}A_n\right)=\bigcap_{n\geq0}\cl_X A_n.
\end{equation}
\end{lemma}
\begin{proof}
Suppose that
\[
 x\in\bigcap_{n\geq0}\cl A_n,\qquad
 U\in\Nbh_{\OO(X)}(x),\qquad U\cap\bigcap_{n\geq0}A_n=\varnothing.
\]
Core compactness gives \(V\in\Nbh_{\OO(X)}(x)\) with \(V\wb U\)
in \(\OO(X)\). Choose \(b_n\in V\cap A_n\) for each \(n\).
Put \(F=X\setminus U\), \(M=\GG(X)\), and
\[
 E_N=\dn_M F\ \cup\ \bigcup_{n\geq N}\dn_M(\dn_X b_n)
 \qquad(N\geq0).
\]
Fix \(N\). The set \(E_N\) is lower in \(M\). Let
\(\mathcal D\subseteq E_N\) be directed, and put
\(C=\bigvee_M\mathcal D=\cl\bigcup\mathcal D\).
If \(C\subseteq F\), then \(C\in E_N\).
If \(C\not\subseteq F\), openness of \(U\) gives
\[
 U\cap\cl\bigcup\mathcal D\ne\varnothing
 \ \Longrightarrow\ U\cap\bigcup\mathcal D\ne\varnothing.
\]
Choose \(C_0\in\mathcal D\), \(e\in C_0\cap U\), and \(m\) with
\(e\notin A_m\). For \(C'\in\mathcal D\) with \(C_0\subseteq C'\),
\[
 C'\not\subseteq F,\qquad
 (\exists n\geq N)\ C'\subseteq\dn_X b_n.
\]
Every such \(n\) satisfies \(n<m\): otherwise
\(e\leq b_n\in A_n\subseteq A_m=\dn_X A_m\), a contradiction.
Thus \(m>N\), and the directed cofinal subset
\(\mathcal D_0=\mathcal D\cap\up_M C_0\) satisfies
\[
 \bigvee_M\mathcal D_0=C,\qquad
 \mathcal D_0\subseteq\bigcup_{N\leq n<m}\dn_M(\dn_X b_n).
\]
If \(\mathcal D_0\not\subseteq\dn_M(\dn_X b_n)\) for every
\(N\leq n<m\), choose
\(C_n\in\mathcal D_0\setminus\dn_M(\dn_X b_n)\).
Directedness gives \(C_*\in\mathcal D_0\) with
\(\bigcup_{N\leq n<m}C_n\subseteq C_*\). Then
\[
 C_*\not\subseteq\dn_X b_n\quad(N\leq n<m),
\]
contrary to the displayed inclusion for \(\mathcal D_0\).
Hence, for some \(N\leq n<m\),
\[
 \mathcal D_0\subseteq\dn_M(\dn_X b_n)
 \quad\Longrightarrow\quad
 C=\bigvee_M\mathcal D_0\subseteq\dn_X b_n.
\]
Therefore \(C\in E_N\), and \(E_N\) is Scott closed.

By continuity of \(\eta_X\),
\[
 C_N:=\eta_X^{-1}(E_N)=F\cup\bigcup_{n\geq N}\dn_X b_n
\]
is closed in \(X\). Clearly \(F\subseteq\bigcap_N C_N\).
For \(y\in U\), choose \(m\) with \(y\notin A_m\). For \(n\geq m\),
\[
 y\leq b_n\quad\Longrightarrow\quad
 y\in\dn_X A_n\subseteq A_m,
\]
a contradiction. Thus \(y\notin C_m\), and \(\bigcap_N C_N=F\).
Set \(O_N=X\setminus C_N\). Then
\[
 O_N\subseteq O_{N+1},\qquad \bigcup_NO_N=U,\qquad b_N\in V\setminus O_N.
\]
But \(V\wb U\) implies \(V\subseteq O_N\) for some \(N\), a
contradiction. Hence
\(\bigcap_n\cl A_n\subseteq\cl\bigcap_nA_n\).
The reverse inclusion follows from
\(\bigcap_nA_n\subseteq A_m\) for every \(m\).
\end{proof}

We use the following consequence of de Brecht's results
\cite[Theorems~23 and 44]{deBrecht2013}. If \(X\) is countably
based, locally compact, and sober, and \((B_i)_{i\geq0}\) is an
open basis, then
\begin{equation}\label{eq:borel-coding}
 e:X\longrightarrow\mathcal P(\NN),\qquad
 e(x)=\{i:x\in B_i\},
\end{equation}
is a topological embedding into \(\Sigma\mathcal P(\NN)\), and
\(e(X)\) is Borel in Cantor space \(2^{\NN}\).
More precisely, there are \(U_k,V_k\in\sigma(\mathcal P(\NN))\)
such that
\[
 e(X)=\bigcap_{k\geq0}
       \bigl((\mathcal P(\NN)\setminus U_k)\cup V_k\bigr).
\]
Scott-open subsets of \(\mathcal P(\NN)\) are Cantor open, so
\(e(X)\) is Cantor Borel. Also,
\(x\leq y\iff e(x)\subseteq e(y)\), and every open subset of
\(X\) is the inverse image of a Cantor-open set. We use the
Cantor topology only for these Borel statements.

\Needspace{7\baselineskip}
\begin{theorem}\label{thm:countable}
Let \(X=(P,\tau)\) be countably based, locally compact, and sober.
If \(\eta_X:X\to\Sigma\GG(X)\) is continuous, then
\(\tau=\sigma(P)\).
\end{theorem}
\begin{proof}
Suppose that \(B\subseteq P\) is Scott closed and
\(x\in\cl_\tau B\setminus B\). Choose a decreasing open local basis
\((V_n)_{n\geq0}\) at \(x\), and choose \(b_n\in B\cap V_n\).
Identify \(P\) with its image under \eqref{eq:borel-coding}; then
\(a\leq_P b\iff a\subseteq b\). Put
\[
 S=\bigcup_{n\geq0}\dn_P b_n,\qquad S_n=S\cap V_n.
\]
Since \(B=\dn_P B\), \(S\subseteq B\). The set \(P\) is Cantor Borel,
and \(\dn_{2^{\NN}}b_n\) is Cantor closed. For each \(n\), choose a
Cantor-open set \(O_n\) with \(V_n=P\cap O_n\). Then
\[
 S=P\cap\bigcup_{m\geq0}\dn_{2^{\NN}}b_m,\qquad S_n=S\cap O_n
\]
are Borel in Cantor space. Moreover, for \(U\in\Nbh_\tau(x)\) and
\(n\geq0\), choose \(m\geq n\) with \(V_m\subseteq U\cap V_n\).
Then \(b_m\in S_n\cap U\). Consequently,
\[
 S_{n+1}\subseteq S_n,\qquad x\in\cl_\tau S_n\quad(n\geq0).
\]

On \(\NN\times2^{\NN}\), define
\[
 (n+1,b)R(n,a)
 \iff a\in S_n,\ b\in S_{n+1},\ a\subseteq b.
\]
The sets \(S_n\) are Borel and inclusion is Cantor closed, so \(R\)
is Borel. If \(z_{m+1}Rz_m\) for every \(m\geq0\), then for some
\(n\geq0\),
\[
 z_m=(n+m,a_m),\qquad
 a_m\leq_Pa_{m+1},\qquad a_m\in B\cap V_{n+m}.
\]
Sobriety makes \(P\) a dcpo. Scott closedness gives
\(a:=\bigvee_{m\geq0}a_m\in B\), and upperness of the \(V_k\) gives
\[
 (\forall k)(\exists m)\quad n+m\geq k,\quad
 a_m\in V_{n+m}\subseteq V_k,\quad a_m\leq a.
\]
Hence \(a\in\bigcap_kV_k=\up_Px\), so \(x\leq a\in B\), contrary
to \(x\notin B=\dn_PB\). Thus \(R\) is well-founded. By the
Kunen--Martin boundedness theorem \cite[Theorem~1, p.~37]{Dell77},
\[
 r(z)=\sup\{r(w)+1:wRz\},\qquad
 \sup_{z\in\NN\times2^{\NN}}r(z)<\omega_1.
\]
Write \(r_n(a)=r(n,a)\). Since \(S_{n+1}\subseteq S_n\),
\begin{equation}\label{eq:rankstep}
 a\in S_{n+1}\quad\Longrightarrow\quad
 (n+1,a)R(n,a)\quad\Longrightarrow\quad r_n(a)\geq r_{n+1}(a)+1.
\end{equation}

For \(U\in\Nbh_\tau(x)\), put
\[
 \rho_n(U)=\sup\{r_n(a):a\in S_n\cap U\},\qquad
 \beta_n=\min\{\rho_n(U):U\in\Nbh_\tau(x)\}.
\]
Choose \(U_n\in\Nbh_\tau(x)\) with \(\rho_n(U_n)=\beta_n\).
For \(U\in\Nbh_\tau(x)\) with \(U\subseteq U_n\),
\[
 \beta_n\leq\rho_n(U)\leq\rho_n(U_n)=\beta_n.
\]
Therefore
\begin{align*}
 \beta_n
 &=\sup\{r_n(a):a\in S_n\cap U_n\cap U_{n+1}\}\\
 &\geq\sup\{r_{n+1}(a)+1:a\in S_{n+1}\cap U_n\cap U_{n+1}\}\\
 &\geq\sup\{r_{n+1}(a):a\in S_{n+1}\cap U_n\cap U_{n+1}\}
 =\beta_{n+1}.
\end{align*}
Choose \(n_0\) with
\(\lambda:=\beta_{n_0}=\min\{\beta_n:n\geq0\}\). Then
\[
 \lambda\leq\beta_n\leq\beta_{n_0}=\lambda\quad(n\geq n_0),
 \qquad\lambda<\omega_1.
\]
Fix \(n\geq n_0\). If \(\lambda=0\), choose
\(a\in S_{n+1}\cap U_n\cap U_{n+1}\); then
\(1\leq r_{n+1}(a)+1\leq r_n(a)\leq\lambda=0\).
If \(\lambda=\alpha+1\), the equality
\[
 \sup\{r_{n+1}(a):a\in S_{n+1}\cap U_n\cap U_{n+1}\}=\alpha+1
\]
gives \(a\) in the displayed set with \(r_{n+1}(a)=\alpha+1\).
Then \(\alpha+2\leq r_n(a)\leq\lambda=\alpha+1\), a contradiction.
Thus \(\lambda\) is a nonzero countable limit ordinal.
For \(n>n_0\), \eqref{eq:rankstep} gives
\begin{equation}\label{eq:strict-rank}
 a\in S_n\cap U_{n-1}\quad\Longrightarrow\quad
 r_n(a)+1\leq r_{n-1}(a)\leq\lambda,\qquad r_n(a)<\lambda.
\end{equation}

Fix \(n>n_0\), and choose
\(\gamma_k<\gamma_{k+1}<\lambda\) with \(\sup_k\gamma_k=\lambda\).
Define
\[
 T_k=\{a\in S_{n+1}\cap U_{n+1}:r_{n+1}(a)\geq\gamma_k\},
 \qquad A_k=\dn_PT_k.
\]
Then \(A_{k+1}\subseteq A_k=\dn_PA_k\subseteq S\).
No Borel assumption on \(T_k\) or \(A_k\) is used.
For every \(U\in\Nbh_\tau(x)\),
\[
 \rho_{n+1}(U\cap U_{n+1})=\lambda>\gamma_k
 \quad\Longrightarrow\quad U\cap T_k\ne\varnothing.
\]
Hence, by Lemma~\ref{lem:fusion},
\[
 x\in\bigcap_k\cl_\tau A_k=\cl_\tau\bigcap_kA_k.
\]
Choose \(b\in U_{n-1}\cap V_n\cap\bigcap_kA_k\).
Since \(b\in S_n\cap U_{n-1}\), \eqref{eq:strict-rank} gives
\(r_n(b)<\lambda\). For each \(k\), choose \(a_k\in T_k\) with
\(b\leq_Pa_k\). Then
\[
 (n+1,a_k)R(n,b),\qquad
 r_n(b)\geq\sup_k\bigl(r_{n+1}(a_k)+1\bigr)
 \geq\sup_k(\gamma_k+1)=\lambda,
\]
a contradiction. Every Scott-closed subset of \(P\) is therefore
\(\tau\)-closed, giving \(\sigma(P)\subseteq\tau\).
Sobriety gives \(\tau\subseteq\sigma(P)\), so \(\tau=\sigma(P)\).
\end{proof}

\Needspace{10\baselineskip}
\subsection{Metrizable lattice quotients}\label{sec:generalweb}

The patch lattice \(\GG(X)\) allows us to apply the countably based
result without imposing a countability assumption on \(X\).

\Needspace{7\baselineskip}
\begin{proof}[Alternative proof of Lemma~\ref{lem:fc-core}]
For \(A\subseteq X\), define
\[
 R_A=\{(a,b)\in X^2:(\exists c\in A)\ c\leq a,b\}.
\]
Proposition~\ref{prop:canonical} makes
\(\mu(a,b)=\eta_X(a)\cap\eta_X(b)\) a continuous map
\(X^2\to\Sigma\GG(X)\). Thus
\[
 U\in\OO(X)\quad\Longrightarrow\quad
 R_U=\mu^{-1}(\hit U)\in\OO(X^2).
\]
Since \((x,x)\in R_U\) for \(x\in U\),
\[
 U\in\Nbh_{\OO(X)}(x)
 \quad\Longrightarrow\quad
 (\exists V\in\Nbh_{\OO(X)}(x))\ V^2\subseteq R_U.
\]
Choose a compact saturated neighborhood \(K_1\subseteq\intr U_1\)
of \(x\). Given \(K_n\), choose \(V\in\Nbh_{\OO(X)}(x)\) with
\(V^2\subseteq R_{\intr K_n}\), and choose a compact saturated
neighborhood \(K_{n+1}\) of \(x\) such that
\[
 K_{n+1}\subseteq V\cap\intr K_n\cap\intr U_{n+1}.
\]
For all \(n\geq1\),
\begin{equation}\label{eq:compact-neighborhood-sequence}
 x\in\intr K_n,\quad K_n\subseteq U_n,\quad
 K_{n+1}\subseteq\intr K_n,\quad K_{n+1}^2\subseteq R_{K_n}.
\end{equation}

Put \(K=\bigcap_{n\geq1}K_n\). Then \(x\in K=\up K\).
For an open cover \(\mathcal U\subseteq\OO(X)\) of \(K\),
well-filteredness and compactness give
\[
 K\subseteq\bigcup\mathcal U
 \ \Longrightarrow\ (\exists n)\ K_n\subseteq\bigcup\mathcal U
 \ \Longrightarrow\ (\exists\mathcal V\in\Fin(\mathcal U))\
 K\subseteq K_n\subseteq\bigcup\mathcal V.
\]
Hence \(K\) is compact. Let \(a,b\in K\), and put
\(F=\dn_Xa\cap\dn_Xb\). By \eqref{eq:compact-neighborhood-sequence},
\[
 a,b\in K_{n+1}\quad\Longrightarrow\quad K_n\cap F\ne\varnothing
 \qquad(n\geq1).
\]
If \(K\cap F=\varnothing\), well-filteredness applied to the open
set \(X\setminus F\) gives \(K_n\subseteq X\setminus F\) for some
\(n\), a contradiction. Thus \(K\) is filtered.
For every \(E\in\Fin(K)\), filteredness and nonemptiness give
\[
 K\cap\bigcap_{k\in E}\dn_Xk\ne\varnothing.
\]
Compactness of \(K\) yields
\[
 m\in\bigcap_{k\in K}(K\cap\dn_Xk),\qquad
 m\leq k\ (k\in K),\qquad K=\up m.
\]
In particular, \(m\leq x\) and \(m\in\bigcap_nU_n\).
For every \(n\), \(m\in K_{n+1}\subseteq\intr K_n\).
For \(O\in\Nbh_{\OO(X)}(m)\), upperness gives \(K=\up m\subseteq O\);
well-filteredness then gives \(K_n\subseteq O\) for some \(n\).
Therefore
\[
 (\forall O\in\Nbh_{\OO(X)}(m))(\exists n)\quad
 m\in\intr K_n\subseteq K_n\subseteq O,
\]
so \((\intr K_n)_{n\geq1}\) is a countable local basis at \(m\).

Let \(D_x=G\cap\dn x\). Taking \(U_n=X\) gives \(D_x\ne\varnothing\).
For \(a,b\in D_x\), let \(\mathcal B_a,\mathcal B_b\) be countable
open local bases. Since \(a,b\leq x\),
\(\mathcal B_a\cup\mathcal B_b\subseteq\Nbh_{\OO(X)}(x)\).
The first part of the proof gives \(g\in D_x\) with
\[
 g\in\bigcap\mathcal B_a\cap\bigcap\mathcal B_b
 =\up a\cap\up b.
\]
Thus \(D_x\) is directed. Taking \(U_1=U\) and \(U_n=X\) for
\(n\geq2\) shows
\(D_x\cap U\ne\varnothing\) for every \(U\in\Nbh_{\OO(X)}(x)\).
Hence \(x\in\cl_X D_x\). Sobriety gives \(s:=\bigvee D_x\), and
\[
 s\leq x,\qquad D_x\subseteq\dn s,\qquad
 \cl_XD_x\subseteq\dn s.
\]
Therefore \(x\leq s\leq x\), and \(\bigvee D_x=x\).
\end{proof}

\Needspace{7\baselineskip}
\begin{lemma}\label{lem:metric-quotient}
Let \(M\) be a compact Hausdorff topological complete lattice.
Every countable family of continuous real-valued functions on \(M\)
factors through a continuous surjective complete lattice homomorphism
\(q:M\to Q\), where \(Q\) is a compact metrizable topological lattice.
\end{lemma}
\begin{proof}
Adjoin the zero function to the given family and enumerate the
functions as \((f_n)_{n\geq1}\), allowing repetitions.
Use the discrete algebraic signature
\[
 E_0=\{\mathbf0,\mathbf1\},\qquad E_2=\{\wedge,\vee\},\qquad
 E_k=\varnothing\quad(k\notin\{0,2\}),
\]
where \(\mathbf0,\mathbf1\) are interpreted as the bounds of \(M\).
The spaces \(E_k\times M^k\) are compact and hence Lindel\"of.
By \cite[Theorem~14]{Reznichenko2024}, for each \(n\) there are a
separable metrizable topological algebra \(A_n\), a continuous
homomorphism \(q_n:M\to A_n\), and a continuous
\(g_n:A_n\to\RR\) such that \(f_n=g_n\circ q_n\).
Replace \(A_n\) by \(q_n[M]\) and \(g_n\) by its restriction.
Then \(A_n\) is compact metrizable, and the bounded lattice identities
hold on \(A_n\) by surjectivity of \(q_n\). Define
\[
 q(x)=(q_n(x))_{n\geq1},\qquad
 Q=q[M]\subseteq\prod_{n\geq1}A_n.
\]
The set \(Q\) is a compact metrizable bounded topological lattice;
\(q:M\to Q\) is a continuous surjective lattice homomorphism.
For the coordinate projection \(\operatorname{pr}_n\),
\[
 f_n=(g_n\circ\operatorname{pr}_n|_Q)\circ q.
\]

The order on \(Q\) is closed because
\[
 \{(a,b):a\leq b\}=\{(a,b):a\wedge b=a\}.
\]
For nonempty \(S\subseteq Q\), index the nets
\(v_F=\bigvee_QF\), \(m_F=\bigwedge_QF\) by
\((\Fin^+(S),\subseteq)\). For \(F\subseteq G\),
\[
 v_F\leq v_G,\qquad m_G\leq m_F.
\]
Equations~\eqref{eq:increasing-net}--\eqref{eq:decreasing-net} give
\(v=\sup_Fv_F\) and \(m=\inf_Fm_F\). For \(z\in Q\),
\[
 (\forall s\in S)\ s\leq z
 \iff(\forall F\in\Fin^+(S))\ v_F\leq z,
\]
and the dual equivalence holds for lower bounds. Thus
\(v=\bigvee_QS\) and \(m=\bigwedge_QS\). The constants \(0_Q,1_Q\)
give the empty supremum and infimum; hence \(Q\) is complete.

For nonempty \(S\subseteq M\), put \(I=\Fin^+(S)\). In the given
topologies of \(M\) and \(Q\), monotone convergence gives
\begin{gather*}
 \bigvee_MF\longrightarrow\bigvee_MS,\qquad
 \bigwedge_MF\longrightarrow\bigwedge_MS,\\
 \bigvee_Qq[F]\longrightarrow\bigvee_Qq[S],\qquad
 \bigwedge_Qq[F]\longrightarrow\bigwedge_Qq[S]
 \qquad(F\in I).
\end{gather*}
By continuity and finite-operation preservation,
\begin{align*}
 q\left(\bigvee_MS\right)
 &=\lim_{F\in I}q\left(\bigvee_MF\right)
  =\lim_{F\in I}\bigvee_Qq[F]=\bigvee_Qq[S],\\
 q\left(\bigwedge_MS\right)
 &=\lim_{F\in I}q\left(\bigwedge_MF\right)
  =\lim_{F\in I}\bigwedge_Qq[F]=\bigwedge_Qq[S].
\end{align*}
For \(S=\varnothing\), use \(q(0_M)=0_Q\) and \(q(1_M)=1_Q\).
Therefore \(q\) is a complete lattice homomorphism.
\end{proof}

\Needspace{7\baselineskip}
\begin{lemma}\label{lem:quotient}
Let \(A\) be a continuous frame with
\(\sigma(A^{\op})=\upsilon(A^{\op})\), and let
\(q:A\to B\) be a surjective complete lattice homomorphism.
Then \(B\) is a continuous frame and
\(\sigma(B^{\op})=\upsilon(B^{\op})\).
\end{lemma}
\begin{proof}
For \(b\in B\) and \((b_j)_{j\in J}\subseteq B\), choose
\(a,a_j\in A\) with \(q(a)=b\) and \(q(a_j)=b_j\). Then
\[
 b\wedge\bigvee_{j\in J}b_j
 =q\left(a\wedge\bigvee_{j\in J}a_j\right)
 =q\left(\bigvee_{j\in J}(a\wedge a_j)\right)
 =\bigvee_{j\in J}(b\wedge b_j).
\]
Thus \(B\) is a frame. Complete join and meet preservation give
\[
 \ell\dashv q\dashv i,\qquad
 \ell(b)=\bigwedge\{a\in A:b\leq q(a)\},\qquad
 i(b)=\bigvee\{a\in A:q(a)\leq b\}.
\]
For \(q(a)=b\), the adjunctions give
\[
 b\leq q\ell(b)\leq q(a)=b,
 \qquad b=q(a)\leq qi(b)\leq b.
\]
Hence \(q\ell=qi=\id_B\).
Fix \(b\in B\), and put \(D_b=\{a\in A:a\wb_A\ell(b)\}\).
For directed \(D\subseteq B\) with \(b\leq\bigvee_BD\),
\[
 a\in D_b
 \ \Longrightarrow\ a\wb_A\ell(b)\leq\ell\left(\bigvee_BD\right)
   =\bigvee_A\ell[D]
 \ \Longrightarrow\ (\exists d\in D)\ q(a)\leq d.
\]
Thus \(q[D_b]\subseteq\{v\in B:v\wb_Bb\}\). Continuity of \(A\)
and join preservation give
\[
 q[D_b]\text{ directed},\qquad
 \bigvee_Bq[D_b]=q\left(\bigvee_AD_b\right)=q\ell(b)=b.
\]
For \(v\wb_Bb\), the last equality yields
\(v\leq q(a)\) for some \(a\in D_b\). Therefore the full set
\(\{v:v\wb_Bb\}\) is directed with supremum \(b\), and \(B\) is continuous.

The map \(q:A^{\op}\to B^{\op}\) preserves directed suprema.
For \(a\in A\), the adjunction \(q\dashv i\) gives
\[
 i^{-1}(A\setminus\up_Aa)=B\setminus\up_Bq(a).
\]
Hence \(i:(B^{\op},\upsilon(B^{\op}))\to
(A^{\op},\upsilon(A^{\op}))\) is continuous. If
\(W\in\sigma(B^{\op})\), then
\[
 q^{-1}(W)\in\sigma(A^{\op})=\upsilon(A^{\op}),\qquad
 W=i^{-1}q^{-1}(W)\in\upsilon(B^{\op}).
\]
Since \(\upsilon(B^{\op})\subseteq\sigma(B^{\op})\), the topologies coincide.
\end{proof}

\Needspace{7\baselineskip}
\begin{proof}[Alternative proof of Theorem~\ref{thm:web}]
Let \(X=(P,\tau)\) satisfy the hypotheses. Suppose that
\(B\subseteq P\) is Scott closed and \(x\in\cl_\tau B\setminus B\).
By Lemma~\ref{lem:fc-core},
\[
 D_x=\{z\leq x:z\text{ is first countable in }X\}
 \text{ is directed},\qquad \bigvee_PD_x=x.
\]
If \(D_x\subseteq B\), then \(x=\bigvee_PD_x\in B\), a contradiction.
Choose \(g\in D_x\setminus B\). Since \(g\leq x\) and
\(\cl_\tau B\) is lower, \(g\in\cl_\tau B\).
Let \((U_n)_{n\geq1}\) be a decreasing open local basis at \(g\).
Choose \(b_n\in B\cap U_n\). Applying Lemma~\ref{lem:fc-core} at
\(b_n\), with the first prescribed neighborhood \(U_n\), gives a
first-countable \(x_n\leq b_n\) in \(U_n\). Thus
\begin{equation}\label{eq:sequence-outside}
 x_n\in B\cap U_n,\qquad x_n\xrightarrow{\tau}g,\qquad g\notin B.
\end{equation}

Put \(M=\GG(X)\) and \(Y=\eta_X(X)\), with the topology
\(\sigma(M)|_Y\). By Proposition~\ref{prop:canonical}, the space
\((M,\pi_M)\) is a compact Hausdorff topological lattice,
\[
 \sigma(M)=\{U\in\pi_M:U=\up_MU\},
\]
and \(\eta_X:X\to Y\) is a homeomorphism and an order isomorphism.
If \((V_n)\) is an open local basis at \(z\in X\), then
\[
 \dn z\in\bigcap_{i=1}^k\hit{W_i}
 \ \Longrightarrow\ z\in\bigcap_{i=1}^kW_i
 \ \Longrightarrow\ (\exists n)\
 \dn z\in\hit{V_n}\subseteq\bigcap_{i=1}^k\hit{W_i}.
\]
Thus \((\hit{V_n})\) is a local basis at \(\eta_X(z)\) in \(\Sigma M\).
Each \(t\) in the countable set
\[
 E=\{\eta_X(g)\}\cup\{\eta_X(x_n):n\geq1\}
\]
therefore has a countable Scott-open local basis \((W_{t,k})_{k\geq1}\)
in \(M\). Since \(W_{t,k}\in\sigma(M)\subseteq\pi_M\), normality
of \((M,\pi_M)\) gives continuous \(f_{t,k}:M\to[0,1]\) with
\[
 f_{t,k}(t)=1,\qquad f_{t,k}[M\setminus W_{t,k}]\subseteq\{0\}.
\]
Apply Lemma~\ref{lem:metric-quotient} to \(\{f_{t,k}:t\in E,k\geq1\}\).
Let \(q:M\to Q\) be the resulting continuous surjective complete
lattice homomorphism. Complete meet preservation gives a left adjoint
\(\ell:Q\to M\), and surjectivity gives
\begin{equation}\label{eq:least-section}
 \ell(a)\leq m\iff a\leq q(m),\qquad
 q\ell=\id_Q,\qquad \ell q\leq\id_M.
\end{equation}
For each \(t,k\), write \(f_{t,k}=\bar f_{t,k}\circ q\). Then
\begin{align*}
 f_{t,k}(\ell q(t))
 &=\bar f_{t,k}(q\ell q(t))
  =\bar f_{t,k}(q(t))=1,\\
 \ell q(t)&\in\bigcap_{k\geq1}W_{t,k}=\up_Mt.
\end{align*}
Combining \(\ell q(t)\in\up_Mt\) with \(\ell q(t)\leq t\) yields
\begin{equation}\label{eq:fixed-points}
 \ell q(t)=t\qquad(t\in E).
\end{equation}

Let \(\pi_Q\) be the compact metrizable topology of \(Q\).
The continuous surjection \(q:(M,\pi_M)\to(Q,\pi_Q)\) is a quotient
map. Since \(q\) preserves directed joins,
\[
 V\in\sigma(Q)\quad\Longrightarrow\quad
 q^{-1}(V)\in\sigma(M)\subseteq\pi_M
 \quad\Longrightarrow\quad V\in\pi_Q.
\]
Conversely, let \(V\in\pi_Q\) be upper, and let \(D\subseteq Q\)
be directed with \(\bigvee_QD\in V\). Equation~\eqref{eq:increasing-net}
gives \(d\xrightarrow{\pi_Q}\bigvee_QD\), indexed by \((D,\leq)\).
Thus \(D\cap V\ne\varnothing\), proving
\begin{equation}\label{eq:quotient-scott-topology}
 \sigma(Q)=\{V\in\pi_Q:V=\up_QV\}.
\end{equation}
For \(O\in\pi_Q\) and \(b\in\up_QO\), choose \(a\in O\) with
\(a\leq b\). The set \(V=\{z:a\wedge z\in O\}\) is \(\pi_Q\)-open,
and
\[
 b\in V,\qquad z\in V\ \Longrightarrow\ a\wedge z\in O,\quad
 a\wedge z\leq z\ \Longrightarrow\ z\in\up_QO.
\]
Hence \(\up_QO\in\pi_Q\). If \(\mathcal B\) is a countable
\(\pi_Q\)-basis, then \(\{\up_QO:O\in\mathcal B\}\) is a
\(\sigma(Q)\)-basis: for \(b\in U\in\sigma(Q)\), choose
\(O\in\mathcal B\) with \(b\in O\subseteq U\); then
\(b\in\up_QO\subseteq\up_QU=U\).

For any complete lattice \(R\), let \(D\subseteq\Pc(R)\) be directed
and put \(p=\bigvee_RD\). Then \(p\ne0\). If \(p\leq a\vee b\)
and \(p\not\leq a,b\), choose \(d_1,d_2\in D\) with
\(d_1\not\leq a\), \(d_2\not\leq b\), and choose \(d\in D\)
with \(d_1,d_2\leq d\). Then
\[
 d\leq p\leq a\vee b,\qquad d\not\leq a,\quad d\not\leq b,
\]
contrary to \(d\in\Pc(R)\). Consequently
\begin{equation}\label{eq:cop-directed-closure}
 D\subseteq\Pc(R)\text{ directed}
 \quad\Longrightarrow\quad\bigvee_RD\in\Pc(R).
\end{equation}
For \(M=\GG(X)\), the elements of \(\Pc(M)\) are exactly the nonempty
irreducible closed subsets of \(X\). Sobriety gives \(Y=\Pc(M)\).

The adjunction \eqref{eq:least-section} and join preservation imply
\begin{equation}\label{eq:coprime}
 a\in\Pc(Q)\quad\Longleftrightarrow\quad\ell(a)\in\Pc(M).
\end{equation}
Indeed, for \(a\in\Pc(Q)\),
\begin{align*}
 \ell(a)\leq m\vee n
 &\Longrightarrow a=q\ell(a)\leq q(m)\vee q(n)\\
 &\Longrightarrow a\leq q(m)\text{ or }a\leq q(n)\\
 &\Longrightarrow \ell(a)\leq m\text{ or }\ell(a)\leq n.
\end{align*}
For \(\ell(a)\in\Pc(M)\),
\begin{align*}
 a\leq b\vee c
 &\Longrightarrow \ell(a)\leq\ell(b)\vee\ell(c)\\
 &\Longrightarrow \ell(a)\leq\ell(b)\text{ or }\ell(a)\leq\ell(c)\\
 &\Longrightarrow a\leq b\text{ or }a\leq c.
\end{align*}
The condition \(a\ne0\) is preserved and reflected because
\(\ell(0)=0\) and \(q\ell=\id_Q\).
Define
\[
 P_Q=\Pc(Q),\quad\tau_{P_Q}=\sigma(Q)|_{P_Q},\quad
 Z=\ell(P_Q)=Y\cap\ell(Q),\quad\tau_Z=\sigma(M)|_Z.
\]
Both \(q\) and \(\ell\) preserve directed joins, and
\(q\ell=\id_Q\). Their restrictions therefore give mutually inverse
continuous maps, and order isomorphisms,
\begin{equation}\label{eq:relative-homeomorphism}
 \ell|_{P_Q}:(P_Q,\tau_{P_Q})\longrightarrow(Z,\tau_Z),
 \qquad q|_Z:(Z,\tau_Z)\longrightarrow(P_Q,\tau_{P_Q}).
\end{equation}
For directed \(D\subseteq Z\), put \(a=\bigvee_Qq[D]\).
Equation~\eqref{eq:cop-directed-closure} gives \(a\in P_Q\), and
\(\ell q(d)=d\) for \(d\in Z\). Hence
\begin{equation}\label{eq:retained-directed-suprema}
 \bigvee_MD=\bigvee_M\ell q[D]
 =\ell\left(\bigvee_Qq[D]\right)=\ell(a)\in Z,
 \qquad \bigvee_ZD=\bigvee_YD=\bigvee_MD.
\end{equation}

Proposition~\ref{prop:canonical} gives a continuous frame
\(A=\OO(X)\cong M^{\op}\) with
\(\sigma(A^{\op})=\upsilon(A^{\op})\).
The dual of \(q\), composed with this isomorphism, is a surjective
complete lattice homomorphism \(A\to Q^{\op}\).
By Lemma~\ref{lem:quotient},
\[
 Q^{\op}\text{ is a continuous frame},\qquad
 \sigma(Q)=\upsilon(Q).
\]
The spectrum of \(Q^{\op}\) is \(P_Q\). Its hull--kernel opens are
\(P_Q\setminus\dn_Q a\), \(a\in Q\). Since \(\sigma(Q)=\upsilon(Q)\),
\[
 \tau_{\mathrm{hk}}(P_Q)=\sigma(Q)|_{P_Q}=\tau_{P_Q}.
\]
Thus \((P_Q,\tau_{P_Q})\) is countably based, locally compact, and
sober. The spectral representation gives an order isomorphism
\[
 \theta:Q\longrightarrow\Gamma(P_Q,\tau_{P_Q}),\qquad
 \theta(a)=P_Q\cap\dn_Qa.
\]
For the continuous inclusion
\(\iota:(P_Q,\tau_{P_Q})\hookrightarrow\Sigma Q\),
\[
 (\theta\circ\iota)(p)=P_Q\cap\dn_Qp
 =\dn_{P_Q}p=\eta_{P_Q}(p).
\]
An order isomorphism is Scott continuous; hence
\(\eta_{P_Q}:(P_Q,\tau_{P_Q})\to\Sigma\Gamma(P_Q,\tau_{P_Q})\)
is continuous. Theorem~\ref{thm:countable} yields
\(\tau_{P_Q}=\sigma(P_Q)\), and
\eqref{eq:relative-homeomorphism} gives \(\tau_Z=\sigma(Z)\).

For \(t\in E\), \eqref{eq:fixed-points} and \eqref{eq:coprime} give
\[
 t\in\Pc(M),\quad\ell q(t)=t
 \quad\Longrightarrow\quad q(t)\in P_Q,\quad t\in Z.
\]
Therefore \(E\subseteq Z\). Since \(\eta_X:P\to Y\) is an order
isomorphism, \(B':=\eta_X[B]\) is Scott closed in \(Y\).
The set \(B'\cap Z\) is lower in \(Z\), and for directed
\(D\subseteq B'\cap Z\),
\[
 \bigvee_ZD=\bigvee_YD\in B'\cap Z
\]
by \eqref{eq:retained-directed-suprema}. Thus \(B'\cap Z\) is
\(\sigma(Z)\)-closed, hence \(\tau_Z\)-closed.
But \eqref{eq:sequence-outside} and continuity of \(\eta_X\) give
\[
 \eta_X(x_n)\in B'\cap Z,\qquad
 \eta_X(x_n)\xrightarrow{\tau_Z}\eta_X(g)\in Z\setminus B',
\]
a contradiction. Hence \(\sigma(P)\subseteq\tau\).
Sobriety gives \(\tau\subseteq\sigma(P)\), completing the proof.
\end{proof}

\Needspace{7\baselineskip}
\begin{proof}[Alternative proof of Corollary~\ref{cor:coframe}]
Put \(X=\Spec(A)\). The Hofmann--Lawson theorem gives
\[
 X\text{ locally compact and sober},\qquad\OO(X)\cong A.
\]
Since \(A\) is a coframe, \eqref{eq:web-coframe} implies that \(X\)
is a web space. Theorem~\ref{thm:web} gives
\(\tau_A=\sigma(\Spec(A))\); the specialization order is inclusion.
\end{proof}

\subsection{Local continuity of join translations}

Let \(L\) be a continuous frame satisfying \eqref{eq:dualcondition},
and retain \(T,h,c_a\) from \eqref{eq:hostdefinition}. We use only the
basic properties in Lemma~\ref{lem:hostbase} until the frame and
coframe laws have been proved below. If \(E\subseteq L\) is upper and
closed under filtered infima, then
\(L\setminus E\in\sigma(L^{\op})=\upsilon(L^{\op})\), so \(E\)
is Lawson closed.

\begin{lemma}\label{lem:escape}
Let \(r,d_n\in L\), \(n\geq1\), satisfy
\(r\leq d_n\vee d_m\) whenever \(n\ne m\). Then
\[
 E=\up r\ \cup\ \bigcup_{n\geq1}\up d_n
\]
is Lawson closed. There are no Lawson-open sets \(W_*,W_n\)
satisfying
\begin{equation}\label{eq:excluded-sequence}
 d_n\in W_*\cap\bigcap_{i<n}W_i,
 \qquad \cl_{\lambda(L)}W_*\cap\up r=\varnothing,
 \qquad \cl_{\lambda(L)}W_n\cap\up d_n=\varnothing.
\end{equation}
\end{lemma}
\begin{proof}
Let \(D\subseteq E\) be filtered. If \(r\leq\bigwedge_LD\), then
\(\bigwedge_LD\in\up r\subseteq E\).
Otherwise choose \(e_0\in D\) with \(r\not\leq e_0\).
Since \(e_0\in E\), there is \(n\geq1\) with \(d_n\leq e_0\).
For distinct \(n,m\),
\[
 d_n,d_m\leq e_0\quad\Longrightarrow\quad
 r\leq d_n\vee d_m\leq e_0,
\]
a contradiction. Thus \(n\) is unique.
For any \(d\in D\), choose \(e\in D\) with \(e\leq e_0,d\).
Then \(r\not\leq e\), so \(d_m\leq e\) for some \(m\).
Since \(d_m\leq e_0\), uniqueness gives \(m=n\). Therefore
\[
 d_n\leq e\leq d\quad(d\in D),\qquad
 d_n\leq\bigwedge_LD,\qquad\bigwedge_LD\in E.
\]
As \(E=\up_LE\), the set \(E\) is Scott closed in \(L^{\op}\).
Condition~\eqref{eq:dualcondition} gives
\[
 L\setminus E\in\sigma(L^{\op})=\upsilon(L^{\op})\subseteq\lambda(L),
\]
so \(E\) is Lawson closed.

Suppose that \eqref{eq:excluded-sequence} holds, and put
\(C_N=\cl_{\lambda(L)}\{d_n:n\geq N\}\).
Then \(\varnothing\ne C_{N+1}\subseteq C_N\subseteq E\).
Lawson compactness gives \(z\in\bigcap_{N\geq1}C_N\). For every \(n\),
\[
 z\in C_1\subseteq E\cap\cl_{\lambda(L)}W_*\subseteq E\setminus\up r,
 \qquad
 z\in C_{n+1}\subseteq\cl_{\lambda(L)}W_n
 \subseteq L\setminus\up d_n.
\]
Hence \(z\in E\setminus(\up r\cup\bigcup_{n\geq1}\up d_n)=\varnothing\),
a contradiction.
\end{proof}

\Needspace{7\baselineskip}
\begin{lemma}\label{lem:diagonal}
For every \(b\in L\), the map \(a\mapsto b\vee a\) is Lawson
continuous at \(b\).
\end{lemma}
\begin{proof}
For \(j_b(a)=b\vee a\) and directed \(D\subseteq L\),
\[
 j_b\left(\bigvee_LD\right)
 =b\vee\bigvee_LD=\bigvee_{d\in D}(b\vee d)=\bigvee_Lj_b[D].
\]
Thus \(j_b\) is Scott continuous. Suppose \(j_b\) is not Lawson
continuous at \(b\). Since \(j_b(b)=b\) and the Lawson topology is
generated by \(\sigma(L)\) and \(\{L\setminus\up u:u\in L\}\),
there is \(u\not\leq b\) such that
\[
 b\in\cl_{\lambda(L)}G,\qquad G=\{a\in L:u\leq b\vee a\}.
\]
Continuity of \(L\) gives \(r\wb u\) with \(r\not\leq b\).
For \(a\in G\),
\[
 a\leq b\ \Longrightarrow\ u\leq b\vee a=b,
\]
so \(a\not\leq b\). The set
\[
 O_a=j_a^{-1}(\Uparrow r)=\{z:r\wb a\vee z\}
\]
is Scott open, and \(r\wb u\leq a\vee b\) implies \(b\in O_a\).
Regularity of \((L,\lambda(L))\) gives
\(W_*\in\Nbh_{\lambda(L)}(b)\) with
\(\cl_{\lambda(L)}W_*\cap\up r=\varnothing\).
Choose \(a_1\in G\cap W_*\).
Having chosen \(a_n\in G\), choose \(W_n\in\Nbh_{\lambda(L)}(b)\)
with \(\cl_{\lambda(L)}W_n\cap\up a_n=\varnothing\).
The set
\[
 V_n=W_*\cap\bigcap_{k=1}^n(W_k\cap O_{a_k})
\]
belongs to \(\Nbh_{\lambda(L)}(b)\), so \(G\cap V_n\ne\varnothing\).
Choose \(a_{n+1}\in G\cap V_n\). For \(m>n\),
\[
 a_m\in O_{a_n}\quad\Longrightarrow\quad r\leq a_n\vee a_m,
 \qquad a_n\in W_*\cap\bigcap_{i<n}W_i.
\]
With \(d_n=a_n\), all conditions in \eqref{eq:excluded-sequence}
hold, contradicting Lemma~\ref{lem:escape}.
\end{proof}

\Needspace{7\baselineskip}
\begin{lemma}\label{lem:nuclei}
Every \(t\in T\) is a nucleus on \(L\). For \(a\in L\) and
\(t\in T\),
\begin{equation}\label{eq:generator-join}
 h(a)\vee_Tt=t\circ c_a.
\end{equation}
Consequently, joining with \(h(a)\) preserves arbitrary meets in \(T\).
\end{lemma}
\begin{proof}
Let \(t\in T\), and choose a net \((a_i)_{i\in I}\subseteq L\) with
\(c_{a_i}\to t\) in the product Lawson topology.
Fix \(x\in L\), and put \(b=t(x)\), \(z_i=a_i\vee x\).
By Lemma~\ref{lem:hostbase}, \(x\leq b\). Thus
\[
 z_i\xrightarrow{\lambda(L)}b,\qquad
 z_i\vee b=a_i\vee b=c_{a_i}(b)\xrightarrow{\lambda(L)}t(b).
\]
Lemma~\ref{lem:diagonal} also gives
\(z_i\vee b\xrightarrow{\lambda(L)}b\).
Hausdorffness yields \(t(b)=b\), that is,
\[
 t(t(x))=t(x)\qquad(x\in L).
\]
Lemma~\ref{lem:hostbase} already gives monotonicity, inflationarity,
and finite-meet preservation, so \(t\) is a nucleus.

For fixed \(a\in L\), coordinate convergence gives
\[
 c_{a_i}\circ c_a=c_{a_i\vee a}\in T,
 \qquad
 c_{a_i\vee a}(x)=c_{a_i}(a\vee x)
 \xrightarrow{\lambda(L)}t(a\vee x)\quad(x\in L).
\]
Closedness of \(T\) in the product implies \(t\circ c_a\in T\).
Monotonicity and inflationarity give
\[
 t(x)\leq t(a\vee x),\qquad a\vee x\leq t(a\vee x),
\]
so \(t,h(a)\leq t\circ c_a\).
If \(v\in T\) and \(t,h(a)\leq v\), then
\[
 a\vee x\leq v(x),\qquad
 v(x)\leq v(a\vee x)\leq v(v(x))=v(x),
\]
and therefore
\[
 t(a\vee x)\leq v(a\vee x)=v(x)\qquad(x\in L).
\]
Thus \(t\circ c_a\) is the least upper bound of \(t,h(a)\) in \(T\),
proving \eqref{eq:generator-join}. For \(S\subseteq T\), pointwise
meet computation gives
\[
 h(a)\vee_T\bigwedge_TS
 =\left(\bigwedge_TS\right)\circ c_a
 =\bigwedge_{t\in S}(t\circ c_a)
 =\bigwedge_{t\in S}(h(a)\vee_Tt).
\]
For \(S=\varnothing\), both sides are the constant map \(c_1=1_T\).
\end{proof}

\Needspace{10\baselineskip}
\subsection{Approximation by join translations}

\Needspace{7\baselineskip}
\begin{lemma}\label{lem:selfapprox}
For every \(t\in T\), there is a net \((a_i)\) in \(L\) such that
\begin{equation}\label{eq:selfapprox}
 h(a_i)\longrightarrow t,
 \qquad t\circ c_{a_i}\longrightarrow t
\end{equation}
in the Lawson topology of \(T\).
\end{lemma}
\begin{proof}
Fix \(t\in T\). For \(U\in\Nbh_{\lambda(T)}(t)\) and \(x,u\in L\)
with \(u\not\leq t(x)\), we first prove
\begin{equation}\label{eq:one-coordinate}
 (\exists a\in L)\quad c_a\in U,\qquad u\not\leq t(a\vee x).
\end{equation}
Suppose, on the contrary, that
\[
 (\forall a\in L)\quad c_a\in U\Longrightarrow u\leq t(a\vee x).
\]
Put \(b=t(x)\), and choose \(r\wb u\) with \(r\not\leq b\).
For \(c_a\in U\), put \(z=a\vee x\). Then
\[
 r\wb u\leq t(z),\qquad
 z\leq b\Longrightarrow u\leq t(z)\leq t(b)=b.
\]
Thus \(r\wb t(z)\) and \(z\not\leq b\).
Choose \(W_*\in\Nbh_{\lambda(L)}(b)\) with
\(\cl_{\lambda(L)}W_*\cap\up r=\varnothing\).
Since \(t\in U\cap\operatorname{ev}_x^{-1}(W_*)\) and \(h(L)\)
is dense in \(T\), choose \(c_{a_1}\) in this intersection and put
\(z_1=a_1\vee x\).

After choosing \(c_{a_n}\in U\) and \(z_n=a_n\vee x\), choose
\(W_n\in\Nbh_{\lambda(L)}(b)\) with
\(\cl_{\lambda(L)}W_n\cap\up z_n=\varnothing\).
Define
\[
 V_n=U\cap\operatorname{ev}_x^{-1}
       \left(W_*\cap\bigcap_{k=1}^nW_k\right)
       \cap\bigcap_{k=1}^n\operatorname{ev}_{z_k}^{-1}(\Uparrow r).
\]
Each evaluation is continuous, \(t(x)=b\), and \(r\wb t(z_k)\).
Hence \(V_n\in\Nbh_{\lambda(T)}(t)\).
Choose \(c_{a_{n+1}}\in h(L)\cap V_n\), and put
\(z_{n+1}=a_{n+1}\vee x\). Then
\begin{gather*}
 z_n\in W_*\cap\bigcap_{i<n}W_i,\qquad
 \cl_{\lambda(L)}W_n\cap\up z_n=\varnothing,\\
 m>n\quad\Longrightarrow\quad
 r\leq c_{a_m}(z_n)=a_m\vee a_n\vee x=z_m\vee z_n.
\end{gather*}
Lemma~\ref{lem:escape}, with \(d_n=z_n\), gives a contradiction.
This proves \eqref{eq:one-coordinate}.

Let \(k\geq1\), and let \((x_j,u_j)\), \(1\leq j\leq k\), satisfy
\(u_j\not\leq t(x_j)\). The continuous map
\[
 \mu_k:T^k\to T,\qquad
 \mu_k(t_1,\ldots,t_k)=\bigwedge_{j=1}^kt_j
\]
satisfies \(\mu_k(t,\ldots,t)=t\in U\).
Choose \(V\in\Nbh_{\lambda(T)}(t)\) with \(\mu_k[V^k]\subseteq U\).
By \eqref{eq:one-coordinate}, choose
\[
 c_{a_j}\in V,\qquad u_j\not\leq t(a_j\vee x_j)
 \quad(1\leq j\leq k).
\]
Set \(a=\bigwedge_{j=1}^ka_j\). Then
\begin{gather*}
 c_a=\bigwedge_{j=1}^kc_{a_j}\in\mu_k[V^k]\subseteq U,\\
 t(a\vee x_j)\leq t(a_j\vee x_j)
 \quad\Longrightarrow\quad u_j\not\leq t(a\vee x_j)
 \quad(1\leq j\leq k).
\end{gather*}
For \(k=0\), choose \(c_a\in U\cap h(L)\).

Define
\begin{gather*}
 E_t=\{(x,u)\in L^2:u\not\leq t(x)\},\qquad
 I_t=\Nbh_{\lambda(T)}(t)\times\Fin(E_t),\\
 (U,F)\leq(V,G)\iff V\subseteq U\text{ and }F\subseteq G.
\end{gather*}
The pair \((U\cap V,F\cup G)\) is an upper bound of \((U,F),(V,G)\),
so \(I_t\) is directed. For \(i=(U,F)\), choose \(a_i\in L\) such that
\[
 c_{a_i}\in U,\qquad
 (\forall(x,u)\in F)\quad u\not\leq t(a_i\vee x).
\]
For \(O\in\Nbh_{\lambda(T)}(t)\),
\[
 i\geq(O,\varnothing)\quad\Longrightarrow\quad h(a_i)\in O.
\]
Thus \(h(a_i)\xrightarrow{\lambda(T)}t\).
For fixed \(x,u\) with \(u\not\leq t(x)\),
\[
 i\geq(T,\{(x,u)\})\quad\Longrightarrow\quad
 t(a_i\vee x)\in L\setminus\up u.
\]
For \(O\in\sigma(L)\) with \(t(x)\in O\), monotonicity gives
\[
 t(x)\leq t(a_i\vee x)\quad\Longrightarrow\quad
 t(a_i\vee x)\in O\qquad(i\in I_t).
\]
The Scott opens and the sets \(L\setminus\up u\) generate
\(\lambda(L)\). Therefore
\[
 t(a_i\vee x)\xrightarrow{\lambda(L)}t(x)\quad(x\in L).
\]
By \eqref{eq:generator-join}, \(t\circ c_{a_i}\in T\).
Lemma~\ref{lem:hostbase} identifies the Lawson topology of \(T\)
with the product subspace topology, so
\(t\circ c_{a_i}\xrightarrow{\lambda(T)}t\), as required.
\end{proof}

\Needspace{7\baselineskip}
\begin{lemma}\label{lem:finite-approx}
For \(u,s_1,\ldots,s_m\in T\), there is a net \((a_\gamma)\)
in \(L\) such that
\[
 h(a_\gamma)\longrightarrow u,\qquad
 h(a_\gamma)\vee_Ts_j\longrightarrow u\vee_Ts_j
 \quad(1\leq j\leq m)
\]
in the Lawson topology of \(T\).
\end{lemma}
\begin{proof}
For \(m=0\), density of \(h(L)\) gives the required net.
Suppose \(m\geq1\), and put \(w_j=u\vee_Ts_j\).
Choose a net \((d_\alpha)_{\alpha\in I}\subseteq L\) with
\(h(d_\alpha)\xrightarrow{\lambda(T)}u\).
By Lemma~\ref{lem:selfapprox} and \eqref{eq:generator-join}, for each
\(j\) there is a net \((b_{j,\beta})_{\beta\in J_j}\subseteq L\) with
\[
 h(b_{j,\beta})\xrightarrow{\lambda(T)}w_j,\qquad
 w_j\vee_Th(b_{j,\beta})\xrightarrow{\lambda(T)}w_j.
\]
Order \(I\times\prod_{j=1}^mJ_j\) coordinatewise, and define
\[
 a_\gamma=d_\alpha\wedge\bigwedge_{j=1}^mb_{j,\beta_j},
 \qquad\gamma=(\alpha,\beta_1,\ldots,\beta_m).
\]
Since \(h\) preserves finite meets and meet is Lawson continuous,
\[
 h(a_\gamma)
 =h(d_\alpha)\wedge\bigwedge_{j=1}^mh(b_{j,\beta_j})
 \xrightarrow{\lambda(T)}u\wedge\bigwedge_{j=1}^mw_j=u.
\]
For each \(j\),
\[
 w_j\leq w_j\vee_Th(a_\gamma)
 \leq w_j\vee_Th(b_{j,\beta_j})\xrightarrow{\lambda(T)}w_j.
\]
Equation~\eqref{eq:squeeze} gives
\(w_j\vee_Th(a_\gamma)\xrightarrow{\lambda(T)}w_j\).
Put \(y_\gamma=h(a_\gamma)\vee_Ts_j\).
For every convergent subnet \(y_{\gamma(\delta)}\to y\),
\[
 h(a_{\gamma(\delta)}),\ s_j
 \leq y_{\gamma(\delta)}
 \leq w_j\vee_Th(a_{\gamma(\delta)}).
\]
Closedness of the order gives \(u,s_j\leq y\leq w_j\), so
\(w_j=u\vee_Ts_j\leq y\leq w_j\). Hence \(y=w_j\).
Compactness and \eqref{eq:unique-subnet-limit} imply
\(y_\gamma\xrightarrow{\lambda(T)}w_j\), proving all required limits.
\end{proof}

\Needspace{7\baselineskip}
\begin{proposition}\label{prop:framehost}
The continuous lattice \(T\) is a frame.
\end{proposition}
\begin{proof}
Fix \(u,v,z\in T\). By Lemma~\ref{lem:finite-approx}, choose one net
\((a_\gamma)\subseteq L\) with
\[
 h(a_\gamma)\vee_Ts\xrightarrow{\lambda(T)}u\vee_Ts
 \qquad(s\in\{v,z,v\wedge z\}).
\]
Lemma~\ref{lem:nuclei} gives, for every \(\gamma\),
\[
 h(a_\gamma)\vee_T(v\wedge z)
 =(h(a_\gamma)\vee_Tv)\wedge(h(a_\gamma)\vee_Tz).
\]
The left-hand side converges to \(u\vee_T(v\wedge z)\).
Joint Lawson continuity of meet makes the right-hand side converge to
\((u\vee_Tv)\wedge(u\vee_Tz)\). Hausdorffness yields
\[
 u\vee_T(v\wedge z)=(u\vee_Tv)\wedge(u\vee_Tz).
\]
Thus \(T\) is distributive. Lemma~\ref{lem:hostbase} gives continuity
of \(T\), so the distributive continuous-lattice result recalled in
Section~\ref{sec:prelim} implies that \(T\) is a frame.
\end{proof}

\Needspace{10\baselineskip}
\subsection{The coframe law}

\Needspace{7\baselineskip}
\begin{proposition}\label{prop:coframehost}
The continuous frame \(T\) is a coframe.
\end{proposition}
\begin{proof}
Let \(u\in T\), and let \(D\subseteq T\) be nonempty filtered.
We prove
\begin{equation}\label{eq:filteredcoframe}
 u\vee_T\bigwedge_TD=\bigwedge_{d\in D}(u\vee_Td).
\end{equation}
Put
\[
 v=\bigwedge_TD,\qquad p=u\vee_Tv,\qquad
 q=\bigwedge_{d\in D}(u\vee_Td).
\]
Since \(v\leq d\) for all \(d\in D\), \(p\leq q\).
Suppose \(p\ne q\). The pointwise order gives \(x\in L\) with
\(q(x)\not\leq p(x)\). Set \(b=p(x)\), and choose
\(r\wb_Lq(x)\) with \(r\not\leq b\).
Every element of \(T\) is a nucleus, and \(u\leq p\). Hence
\begin{equation}\label{eq:fixed-b}
 b\leq u(b)\leq p(b)=p(p(x))=b,\qquad u(b)=b.
\end{equation}
Choose \(W_*\in\Nbh_{\lambda(L)}(b)\) with
\(\cl_{\lambda(L)}W_*\cap\up r=\varnothing\).

We construct \(a_n,b_n,z_n\in L\) and
\(W_n\in\Nbh_{\lambda(L)}(b)\), for \(n\geq1\), satisfying
\begin{align}
 z_n&=a_n\vee b_n\vee x\in W_*\cap\bigcap_{i<n}W_i,
       \label{eq:triangle1}\\
 r&\wb u(b_n\vee x),\label{eq:triangle2}\\
 r&\leq a_n\vee b_i\vee x\quad(1\leq i<n),\label{eq:triangle3}\\
 \cl_{\lambda(L)}W_n&\cap\up z_n=\varnothing.\label{eq:triangle4}
\end{align}
Assume that the choices are made for \(i<n\), and set
\(W=W_*\cap\bigcap_{i<n}W_i\).
By Lemma~\ref{lem:finite-approx}, applied to \(u,v\), choose
\((\alpha_j)_{j\in J}\subseteq L\) with
\[
 h(\alpha_j)\xrightarrow{\lambda(T)}u,\qquad
 h(\alpha_j)\vee_Tv\xrightarrow{\lambda(T)}p.
\]
Evaluation and \eqref{eq:generator-join} give
\begin{align*}
 v(\alpha_j\vee x)&\xrightarrow{\lambda(L)}p(x)=b\in W,\\
 \alpha_j\vee b_i\vee x&\xrightarrow{\lambda(L)}u(b_i\vee x)
 \in\Uparrow r\qquad(1\leq i<n).
\end{align*}
The second membership follows from the inductive hypothesis
\eqref{eq:triangle2}. Define
\begin{gather*}
 J_0=\{j\in J:v(\alpha_j\vee x)\in W\},\\
 J_i=\{j\in J:r\wb\alpha_j\vee b_i\vee x\}
 \quad(1\leq i<n).
\end{gather*}
Each \(J_i\) contains a final segment of \(J\). Their finite
intersection therefore contains an index \(j_n\). Set
\(a_n=\alpha_{j_n}\). Then
\begin{equation}\label{eq:choose-an}
 v(a_n\vee x)\in W,\qquad
 r\wb a_n\vee b_i\vee x\quad(1\leq i<n).
\end{equation}

Index the net \((d)\) by \((D,\geq_T)\). Equation~\eqref{eq:decreasing-net}
and continuity of evaluation give
\[
 d\xrightarrow{\lambda(T)}v,\qquad
 d(a_n\vee x)\xrightarrow{\lambda(L)}v(a_n\vee x)\in W.
\]
Choose \(d_n\in D\) with \(d_n(a_n\vee x)\in W\).
No inequality between \(d_n\) and \(d_{n-1}\) is imposed.
Apply Lemma~\ref{lem:finite-approx} to \(d_n,u\), obtaining
\((\beta_j)_{j\in I_n}\subseteq L\) with
\[
 h(\beta_j)\xrightarrow{\lambda(T)}d_n,\qquad
 h(\beta_j)\vee_Tu\xrightarrow{\lambda(T)}d_n\vee_Tu.
\]
Thus
\begin{align*}
 \beta_j\vee a_n\vee x&\xrightarrow{\lambda(L)}d_n(a_n\vee x)\in W,\\
 u(\beta_j\vee x)&\xrightarrow{\lambda(L)}(d_n\vee_Tu)(x)
 \in\Uparrow r,
\end{align*}
where \(r\wb q(x)\leq(d_n\vee_Tu)(x)\) proves the second membership.
Both index sets
\[
 I_n^0=\{j:\beta_j\vee a_n\vee x\in W\},\qquad
 I_n^1=\{j:r\wb u(\beta_j\vee x)\}
\]
contain final segments of \(I_n\). Choose \(k_n\in I_n^0\cap I_n^1\)
and set \(b_n=\beta_{k_n}\), \(z_n=a_n\vee b_n\vee x\). Then
\[
 z_n\in W,\qquad r\wb u(b_n\vee x),\qquad
 r\leq a_n\vee b_i\vee x\quad(1\leq i<n).
\]
These are \eqref{eq:triangle1}--\eqref{eq:triangle3}.
If \(z_n\leq b\), monotonicity and \eqref{eq:fixed-b} give
\[
 r\leq u(b_n\vee x)\leq u(z_n)\leq u(b)=b,
\]
a contradiction. Hence \(b\notin\up z_n\), and regularity gives
\(W_n\in\Nbh_{\lambda(L)}(b)\) satisfying \eqref{eq:triangle4}.
For \(n=1\), the index sets \(J_i\) with \(1\leq i<n\) are absent,
so the same construction starts with \(W=W_*\).

For \(1\leq i<n\), \eqref{eq:triangle3} implies
\[
 r\leq a_n\vee b_i\vee x\leq z_n\vee z_i.
\]
Thus \((z_n)\), \(W_*\), and \((W_n)\) satisfy all conditions
excluded by Lemma~\ref{lem:escape}. This contradiction proves
\eqref{eq:filteredcoframe}.

For arbitrary \(S\subseteq T\), define
\[
 D_S=\{\bigwedge_TF:F\in\Fin(S)\}.
\]
The set \(D_S\) contains \(1_T=\bigwedge_T\varnothing\), and
\[
 \bigwedge_T(F\cup G)\leq\bigwedge_TF,\bigwedge_TG
 \qquad(F,G\in\Fin(S)).
\]
Hence \(D_S\) is nonempty filtered, with
\(\bigwedge_TD_S=\bigwedge_TS\).
Equation~\eqref{eq:filteredcoframe} and finite distributivity give
\begin{align*}
 u\vee_T\bigwedge_TS
 &=\bigwedge_{F\in\Fin(S)}\left(u\vee_T\bigwedge_TF\right)\\
 &=\bigwedge_{F\in\Fin(S)}\ \bigwedge_{s\in F}(u\vee_Ts)
  =\bigwedge_{s\in S}(u\vee_Ts).
\end{align*}
For \(S=\varnothing\), both sides equal \(1_T\).
Therefore \(T\) is a coframe.
\end{proof}

\Needspace{7\baselineskip}
\begin{proof}[Alternative proof of Theorem~\ref{thm:embedding}]
Use \(T=\cl_{L^L}\{c_a:a\in L\}\) and \(h(a)=c_a\) from
\eqref{eq:hostdefinition}. Lemma~\ref{lem:hostbase} gives continuity
of \(T\) and, for \(S\subseteq L\), \(F\in\Fin(L)\), and \(a,b\in L\),
\begin{gather*}
 h(a)=h(b)\Longrightarrow a=b,\qquad
 h\left(\bigvee_LS\right)=\bigvee_Th[S],\\
 h\left(\bigwedge_LF\right)=\bigwedge_Th[F],\qquad
 a\wb_Lb\Longrightarrow h(a)\wb_Th(b).
\end{gather*}
Propositions~\ref{prop:framehost} and \ref{prop:coframehost} give the
frame and coframe laws for \(T\). Thus \(h\) has every asserted property.
\end{proof}

\subsection{Transfer through completely prime filters}\label{ap:transfer-section}

We use the separation principle \eqref{eq:filter-separation} to prove
Theorem~\ref{thm:transfer} directly in the completely prime filter
representation.

\begin{proof}[Alternative proof of Theorem~\ref{thm:transfer}]
For \(F\in\Spec(L)\), define
\[
 H_F=\{u\in A:(\exists a\in F)\ e(a)\leq u\}.
\]
Then \(H_F=\up_AH_F\) and \(1_A=e(1_L)\in H_F\).
If \(u,v\in H_F\), choose \(a,b\in F\) with
\(e(a)\leq u\) and \(e(b)\leq v\). Then
\[
 a\wedge b\in F,\qquad
 e(a\wedge b)=e(a)\wedge e(b)\leq u\wedge v,
\]
so \(u\wedge v\in H_F\). Thus \(H_F\) is a lattice filter.
For directed \(D\subseteq A\) with \(\bigvee_AD\in H_F\), choose
\(a\in F\) with \(e(a)\leq\bigvee_AD\).
Continuity of \(L\) and complete primality of \(F\) give
\[
 a=\bigvee\{b:b\wb_La\}\in F
 \quad\Longrightarrow\quad(\exists b\in F)\ b\wb_La.
\]
Hence
\[
 e(b)\wb_Ae(a)\leq\bigvee_AD
 \quad\Longrightarrow\quad(\exists d\in D)\ e(b)\leq d
 \quad\Longrightarrow\quad D\cap H_F\ne\varnothing.
\]
Therefore \(H_F\in\sigma(A)\).

For \(G\in\Spec(A)\), monotonicity of \(e\) gives
\(e^{-1}(G)=\up_Le^{-1}(G)\). Preservation of finite meets and bounds gives
\(1_L\in e^{-1}(G)\), \(0_L\notin e^{-1}(G)\), and finite-meet
closure of \(e^{-1}(G)\). For \(S\subseteq L\),
\[
 \bigvee_LS\in e^{-1}(G)
 \iff\bigvee_Ae[S]\in G
 \iff S\cap e^{-1}(G)\ne\varnothing.
\]
Thus
\[
 f:\Spec(A)\to\Spec(L),\qquad f(G)=e^{-1}(G)
\]
is well defined. Injectivity and meet preservation imply that \(e\)
reflects order:
\[
 e(a)\leq e(b)\ \Longrightarrow\ e(a\wedge b)=e(a)
 \ \Longrightarrow\ a\wedge b=a\ \Longrightarrow\ a\leq b.
\]
Fix \(F\in\Spec(L)\), and put \(p=\bigvee_L(L\setminus F)\).
By the prime-element representation, \(p\notin F\) and
\(L\setminus F=\dn_Lp\). If \(e(p)\in H_F\), then some \(a\in F\)
satisfies \(e(a)\leq e(p)\), whence \(a\leq p\) and \(p\in F\),
a contradiction. Therefore \(e(p)\notin H_F\).
By \eqref{eq:filter-separation}, choose \(G\in\Spec(A)\) with
\(H_F\subseteq G\) and \(e(p)\notin G\). For \(a\in L\),
\begin{align*}
 a\in F&\ \Longrightarrow\ e(a)\in H_F\subseteq G,\\
 a\notin F&\ \Longrightarrow\ a\leq p
 \ \Longrightarrow\ e(a)\leq e(p)
 \ \Longrightarrow\ e(a)\notin G.
\end{align*}
Hence \(f(G)=F\), proving surjectivity of \(f\).

Let \(W\subseteq\Spec(L)\) be upper and suppose
\(f^{-1}(W)=\hata v\) for \(v\in A\). Fix \(F\in W\).
If \(v\notin H_F\), \eqref{eq:filter-separation} gives
\(G\in\Spec(A)\) with \(H_F\subseteq G\) and \(v\notin G\).
Then
\[
 F\subseteq f(G),\quad F\in W=\up W
 \quad\Longrightarrow\quad f(G)\in W
 \quad\Longrightarrow\quad G\in\hata v,
\]
contrary to \(v\notin G\). Thus \(v\in H_F\), so choose \(a\in F\)
with \(e(a)\leq v\). For \(F'\in\hata a\), choose \(G'\in\Spec(A)\)
with \(f(G')=F'\). Then
\[
 a\in F'=e^{-1}(G')\ \Longrightarrow\ e(a)\in G'
 \ \Longrightarrow\ v\in G'
 \ \Longrightarrow\ F'\in W.
\]
Consequently \(F\in\hata a\subseteq W\).
Since \(F\in W\) was arbitrary,
\begin{equation}\label{eq:open-descent}
 W=\bigcup\{\hata a:a\in L,\ \hata a\subseteq W\}\in\tau_L.
\end{equation}

For directed \(\mathcal D\subseteq\Spec(A)\),
\[
 f\left(\bigvee\mathcal D\right)
 =e^{-1}\left(\bigcup_{G\in\mathcal D}G\right)
 =\bigcup_{G\in\mathcal D}e^{-1}(G)
 =\bigvee f[\mathcal D].
\]
Thus \(f\) is Scott continuous. For \(W\in\sigma(\Spec(L))\),
\[
 f^{-1}(W)\in\sigma(\Spec(A))=\tau_A,
\]
so \(f^{-1}(W)=\hata v\) for some \(v\in A\).
Equation~\eqref{eq:open-descent} gives \(W\in\tau_L\).
Therefore \(\sigma(\Spec(L))\subseteq\tau_L\).
The reverse inclusion holds for every frame, so
\(\tau_L=\sigma(\Spec(L))\).
\end{proof}

\begin{proof}[Alternative proof of Theorem~\ref{thm:main}]
The necessary implication is \cite[Theorem~4.15]{CK23}. Conversely, let \(L\) be a distributive continuous lattice satisfying \eqref{eq:dualcondition}. By the alternative proof of Theorem~\ref{thm:embedding} above, the construction \eqref{eq:hostdefinition} gives a continuous frame-coframe \(T\) and a frame embedding \(h:L\to T\) with Scott-continuous right adjoint. By the alternative proof of Corollary~\ref{cor:coframe} above, \(\Spec(T)\) is a Scott space. The alternative proof of Theorem~\ref{thm:transfer} above, applied to \(h\), completes the argument.
\end{proof}

\section{Weak hypercontinuity and dual core-compactness}\label{app:dual-compactness}

This appendix relates weak hypercontinuity to dual Scott compactness and
then to the spectrum. The main additional implication is that a spatial
frame whose dual Scott space is core-compact is necessarily continuous.

\subsection{Classical characterizations}\label{app:unified}

Put \(\beta(L)=\sigma(L)\vee\sigma(L^{\op})\). A complete lattice
is \emph{lean} if every dual Scott-closed set is Scott compact.
The following equivalences are
\cite[Theorems~4.15 and~4.19]{CK23} and
\cite[Proposition~4.4 and its proof]{CKL2022}.

\begin{proposition}\label{un:weak-equivalences}
For a continuous lattice \(L\), the following conditions are equivalent:
\begin{enumerate}[label=\textup{(\roman*)}]
\item \(L\) is weak hypercontinuous;
\item \(L\) is lean;
\item \(\beta(L)=\lambda(L)\);
\item \((L,\beta(L))\) is compact Hausdorff.
\end{enumerate}
They are also equivalent to leanness of \(L^{\op}\).
\end{proposition}

For a complete lattice, write \(a\prec b\) when
\(b\in\intr_{\upsilon(L)}\uparrow a\).
Recall that \(L\) is \emph{hypercontinuous} if
\(b=\bigvee\{a:a\prec b\}\), and \emph{hyperalgebraic} if
\(b=\bigvee\{a:a\prec a\leq b\}\), for every \(b\in L\).
A dcpo is \emph{quasicontinuous} if each Scott-open neighbourhood
\(U\) of a point \(x\) contains a finitely generated upper set
\(\uparrow F\) with \(x\in\intr_{\sigma(P)}\uparrow F\);
it is \emph{quasialgebraic} if \(\uparrow F\) can be chosen Scott open.
We use the standard Gierz--Lawson characterization
\cite[Theorem~6.4\textup{(iv)}]{GL81}; see also
\cite[Theorem~2.2]{CKL2022}.

\begin{theorem}\label{un:hyper-weak}
A continuous lattice \(L\) is hypercontinuous if and only if it is
weak hypercontinuous and \(L^{\op}\) is quasicontinuous.
\end{theorem}

\begin{corollary}\label{un:algebraic}
For an algebraic lattice, weak hypercontinuity, hypercontinuity and
hyperalgebraicity are equivalent.
\end{corollary}
\begin{proof}
If \(L\) is algebraic and weak hypercontinuous, then
\(\sigma(L^{\op})=\omega(L)\) is algebraic by
\cite[Lemma~3.1]{CKL2022}. Hence \(L^{\op}\) is quasialgebraic by
\cite[Theorem~2.5\textup{(3)}]{CKL2022}, and
Theorem~\ref{un:hyper-weak} gives hypercontinuity.
In an algebraic hypercontinuous lattice, every compact element
\(k\) satisfies \(k\prec k\): apply compactness to
\(k=\bigvee\{a:a\prec k\}\), whose joining family is directed
because it is closed under finite joins. Thus the lattice is hyperalgebraic.
The reverse implications follow from the definitions and
Theorem~\ref{un:hyper-weak}.
\end{proof}

The distinction from hypercontinuity is strict. For Lawson's
meet-continuous noncontinuous lattice \(W\) with core-compact Scott
space, the frame \(L=\OO(\Sigma W)\) is continuous and weak
hypercontinuous but not hypercontinuous
\cite[discussion following Theorem~4.15]{CK23}.

\subsection{Reflection of continuity for spatial frames}\label{cc:section}

A frame \(L\) is spatial if the canonical map
\(L\to\OO(\Spec L)\) is an isomorphism. We prove the implication
from dual Scott core-compactness to continuity without assuming
continuity in advance. The first step uses the patch topology
\(\theta_L=\operatorname{patch}\Sigma(L^{\op})\).

\begin{lemma}\label{cc:patch-meet}
Let \(L\) be a frame. If \(\Sigma(L^{\op})\) is core-compact, then
\[
 \theta_L=\operatorname{patch}\Sigma(L^{\op})
\]
is compact Hausdorff, the order relation of \(L\) is closed in
\(\theta_L\times\theta_L\), and binary meet is jointly
\(\theta_L\)-continuous.
\end{lemma}
\begin{proof}
Put \(P=L^{\op}\). By \cite[Corollary~3.3]{JJL2016}, the Scott
space of the complete lattice \(P\) is stably compact. In
particular, it is locally compact and sober. Its patch topology is
compact Hausdorff by \cite[Theorem~2.5]{JKM2001}.

For \(a\in P\), define \(c_a(x)=a\vee_Px\). Since \(P\) is a
coframe, \(c_a\) preserves arbitrary infima and has a left adjoint
\[
 \rho_a(y)=\bigwedge_P\{z\in P:y\leq_Pa\vee_Pz\}.
\]
Thus
\[
 \rho_a(y)\leq_Px\quad\Longleftrightarrow\quad y\leq_Pc_a(x).
\]
The map \(\rho_a\) preserves arbitrary suprema, and \(c_a\)
preserves nonempty suprema. Both maps are Scott continuous.
If \(K\) is Scott compact and saturated, then
\begin{equation}\label{cc:eq:adjoint-preimage}
 c_a^{-1}(K)=\uparrow_P\rho_a[K].
\end{equation}
Indeed, \(c_a(x)\in K\) implies
\(\rho_a(c_a(x))\leq_Px\). Conversely, if \(k\in K\) and
\(\rho_a(k)\leq_Px\), then \(k\leq_Pc_a(x)\), so
\(c_a(x)\in K\) because \(K\) is upper.
The set \(\rho_a[K]\) is Scott compact, and its upper closure is
Scott compact because every Scott-open set is upper. Hence
\eqref{cc:eq:adjoint-preimage} is compact saturated. The inverse
image under \(c_a\) of each generator of the patch topology is
therefore patch open.

Consequently \((P,\theta_L,\vee_P)\) is a compact Hausdorff
semitopological semilattice. By \cite[Theorem~6.6]{Lawson1973},
its operation is jointly continuous. Since \(\vee_P=\wedge_L\),
binary meet on \(L\) is jointly continuous. Finally,
\[
 \{(x,y):x\leq_Ly\}=\{(x,y):x\wedge_Ly=x\}
\]
is closed because \(\theta_L\) is Hausdorff.
\end{proof}

An indexed family \((F_i)_{i\in I}\) is point-finite if
\(\{i\in I:x\in F_i\}\) is finite for every point \(x\).
It is locally finite if every point has a neighbourhood meeting
only finitely many members, with indices counted separately.

\begin{lemma}\label{cc:point-finite}
If \(\Sigma \GG(X)\) is core-compact, then every point-finite family
of relatively closed subsets of every open subspace of \(X\) is
locally finite.
\end{lemma}
\begin{proof}
Let \(U\) be open in \(X\), put \(C_0=X\setminus U\), and write
\(F_i=A_i\cap U\), where \(A_i\) is closed in \(X\).
Suppose that the point-finite family \((F_i)_{i\in I}\) is not
locally finite at \(x\in U\). For every finite \(J\subseteq I\),
\begin{equation}\label{cc:eq:tail-closure}
 x\in\cl_X\bigcup_{i\in I\setminus J}A_i.
\end{equation}
Indeed, every neighbourhood of \(x\) meets some \(F_i\) with
\(i\notin J\).

Put \(P=\GG(X)\) and \(H=P\setminus\downarrow C_0\). Then \(H\)
is Scott open and contains \(\cl_X\{x\}\). Core-compactness gives
a Scott-open set \(V\) such that
\[
 \cl_X\{x\}\in V\ll H,
\]
where \(\ll\) is calculated in \(\OO(\Sigma P)\).
For finite \(J\subseteq I\), the finite unions of members of
\((A_i)_{i\in I\setminus J}\) form a directed family whose
supremum contains \(\cl_X\{x\}\), by
\eqref{cc:eq:tail-closure}. Scott openness of \(V\) therefore
gives a finite \(J'\subseteq I\setminus J\) with
\(\bigcup_{i\in J'}A_i\in V\). Since
\(\varnothing\notin V\subseteq H\), the set \(J'\) is nonempty.
Recursively choose pairwise disjoint nonempty finite sets
\(I_n\subseteq I\), \(n\geq1\), such that
\[
 G_n=\bigcup_{i\in I_n}A_i\in V.
\]
The family \((G_n\cap U)_{n\geq1}\) is point-finite.

For \(n\geq1\), define
\[
 D_n=\downarrow C_0\;\cup\!\bigcup_{k\geq n}\downarrow G_k.
\]
We prove that \(D_n\) is Scott closed. It is a lower set. Let
\(\mathcal E\subseteq D_n\) be directed. If
\(E\subseteq C_0\) for every \(E\in\mathcal E\), then
\(\bigvee\mathcal E\subseteq C_0\). Otherwise, choose
\(E_0\in\mathcal E\) and \(y\in E_0\setminus C_0\). The set
\[
 K_0=\{k\geq n:E_0\subseteq G_k\}
\]
is nonempty and finite, since every such \(G_k\) contains
\(y\in U\). For \(E\in\mathcal E\), put
\[
 K(E)=\{k\in K_0:E\subseteq G_k\}.
\]
If \(E_1,\ldots,E_m\in\mathcal E\), choose
\(E'\in\mathcal E\) containing \(E_0,E_1,\ldots,E_m\).
Since \(E'\nsubseteq C_0\), some \(k\geq n\) satisfies
\(E'\subseteq G_k\), and this \(k\) belongs to
\(\bigcap_{j=1}^mK(E_j)\). The subsets \(K(E)\) of the finite
set \(K_0\) thus have the finite intersection property. There is
a fixed \(k\in K_0\) with \(E\subseteq G_k\) for every
\(E\in\mathcal E\). Since \(G_k\) is closed,
\[
 \bigvee\mathcal E=\cl_X\bigcup\mathcal E\subseteq G_k.
\]
Hence \(D_n\) is Scott closed.

Point-finiteness gives \(\bigcap_nD_n=\downarrow C_0\): if
\(E\nsubseteq C_0\), any \(y\in E\setminus C_0\) belongs to
only finitely many \(G_k\), so \(E\notin D_n\) for sufficiently
large \(n\). The increasing Scott-open family
\((P\setminus D_n)_{n\geq1}\) therefore has union \(H\).
Since \(V\ll H\), some \(n\) satisfies
\(V\subseteq P\setminus D_n\). This contradicts
\(G_n\in V\cap D_n\).
\end{proof}

\begin{lemma}\label{cc:spatial-seq}
Let \(L\) be a spatial frame such that \(\Sigma(L^{\op})\)
is core-compact. For every \(\theta_L\)-convergent sequence,
\begin{equation}\label{cc:eq:sequence}
 x_n\longrightarrow v
 \quad\Longrightarrow\quad
 v=\bigvee_{N\geq1}\bigwedge_{n\geq N}x_n.
\end{equation}
\end{lemma}
\begin{proof}
Let \(X=\Spec L\), and identify its points with the proper
meet-prime elements of \(L\) by \eqref{eq:prime-spectrum}.
Put \(U_a=\{p:a\nleq p\}\). Spatiality gives
\(L\cong\OO(X)\), and therefore
\begin{equation}\label{cc:eq:spatial-operations}
 U_{\bigvee_i a_i}=\bigcup_iU_{a_i},
 \qquad
 U_{\bigwedge_i a_i}=\intr_X\bigcap_iU_{a_i}.
\end{equation}
Moreover, \(\GG(X)\cong L^{\op}\), so
Lemma~\ref{cc:point-finite} applies to \(X\).

Suppose that \(x_n\to v\) in \(\theta_L\). In the open subspace
\(U_v\), consider the relatively closed sets
\(A_n=U_v\setminus U_{x_n}\). If \(p\in U_v\), then
\(v\notin\downarrow p\). Lemma~\ref{cc:patch-meet} implies
that \(L\setminus\downarrow p\) is a \(\theta_L\)-open
neighbourhood of \(v\). Thus \(x_n\nleq p\) eventually, so
\(p\) belongs to only finitely many \(A_n\). By
Lemma~\ref{cc:point-finite}, the family \((A_n)\) is locally
finite. Hence every \(p\in U_v\) has an open neighbourhood
\(B_p\) in \(X\) and an integer \(N_p\) such that
\[
 p\in B_p\subseteq U_v,
 \qquad
 B_p\subseteq\bigcap_{n\geq N_p}U_{x_n}.
\]
Since \(B_p\) is open, \eqref{cc:eq:spatial-operations} gives
\(B_p\subseteq U_{\bigwedge_{n\geq N_p}x_n}\). Therefore
\[
 U_v\subseteq\bigcup_{N\geq1}U_{\bigwedge_{n\geq N}x_n}
 =U_{\bigvee_N\bigwedge_{n\geq N}x_n}.
\]
Spatiality yields \(v\leq\bigvee_N\bigwedge_{n\geq N}x_n\).

For the reverse inequality, fix \(N\) and put
\(c_N=\bigwedge_{n\geq N}x_n\). The \(\theta_L\)-closed set
\(\uparrow c_N\) contains the tail \((x_n)_{n\geq N}\), so
it contains \(v\). Thus \(c_N\leq v\) for every \(N\), which
proves \eqref{cc:eq:sequence}.
\end{proof}

The next result is independent of spatiality and distributivity.
As in Section~\ref{sec:prelim}, \(\Nbh_\theta(v)\) denotes
the open neighbourhoods of \(v\).

\begin{theorem}\label{cc:compact-meet}
Let \(M\) be a complete lattice equipped with a compact Hausdorff
topology \(\theta\) for which binary meet is jointly continuous.
Suppose that every \(\theta\)-convergent sequence satisfies
\[
 x_n\longrightarrow x
 \quad\Longrightarrow\quad
 x=\bigvee_{N\geq1}\bigwedge_{n\geq N}x_n.
\]
Then, for every \(v\in M\),
\begin{equation}\label{cc:eq:neighbourhood-infima}
 v=\bigvee_{U\in\Nbh_\theta(v)}\bigwedge U,
 \qquad
 \bigwedge U\ll v\quad(U\in\Nbh_\theta(v)).
\end{equation}
In particular, \(M\) is continuous.
\end{theorem}
\begin{proof}
The order relation is closed, since \(x\leq y\) is equivalent
to \(x\wedge y=x\). Equations~\eqref{eq:increasing-net} and
\eqref{eq:decreasing-net} therefore give convergence of directed
sets to their suprema and filtered sets to their infima.
In particular, for any nonempty \(A\subseteq M\),
the net
\[
 \left(\bigwedge F\right)_{F\in\mathcal P_{\mathrm{fin}}^+(A)}
 \longrightarrow\bigwedge A,
\]
where \(\mathcal P_{\mathrm{fin}}^+(A)\) consists of the nonempty
finite subsets of \(A\), ordered by inclusion.

Let \(U\in\Nbh_\theta(v)\), and let \(D\) be directed with
\(s=\bigvee D\geq v\). Since \(d\to s\), continuity of meet
gives \(d\wedge v\to s\wedge v=v\). Some \(d\wedge v\) lies
in \(U\), whence \(\bigwedge U\leq d\wedge v\leq d\).
Therefore \(\bigwedge U\ll v\).

Put
\[
 a=\bigvee_{U\in\Nbh_\theta(v)}\bigwedge U\leq v,
\]
and suppose that \(a<v\). By regularity, there is a
\(\theta\)-open neighbourhood \(O\) of \(v\) such that
\(\cl_\theta O\cap\downarrow a=\varnothing\).
The set \(\uparrow O\) is open: if \(x\geq y\in O\), then
the continuous map \(m_y(t)=t\wedge y\) satisfies
\[
 x\in m_y^{-1}(O)\subseteq\uparrow O.
\]
The set \(\uparrow\cl_\theta O\) is closed, because it is
the projection on the second coordinate of the compact set
\[
 \{(y,x):y\in\cl_\theta O,\ y\leq x\}.
\]
Define
\[
 W=M\setminus\uparrow\cl_\theta O,
 \qquad
 C=M\setminus\uparrow O.
\]
Then \(W\) is an open lower set, \(C\) is a closed lower set,
and
\begin{equation}\label{cc:eq:closed-lower}
 \downarrow a\subseteq W\subseteq C,
 \qquad v\notin C.
\end{equation}

Choose a compact neighbourhood \(K_1\) of \(v\) contained in
\(M\setminus C\). Recursively choose compact neighbourhoods
\(K_n\) of \(v\) such that
\begin{equation}\label{cc:eq:compact-neighbourhoods}
 K_{n+1}\subseteq K_n,
 \qquad K_{n+1}\wedge K_{n+1}\subseteq K_n.
\end{equation}
To perform the recursive step, choose open neighbourhoods
\(A,B\) of \(v\) with \(A\wedge B\subseteq\intr_\theta K_n\).
Choose an open neighbourhood \(T\) of \(v\) such that
\(\cl_\theta T\subseteq A\cap B\cap\intr_\theta K_n\), and
put \(K_{n+1}=\cl_\theta T\).

The compact set \(K=\bigcap_nK_n\) contains \(v\) and is closed
under binary meet, by \eqref{cc:eq:compact-neighbourhoods}.
Its finite nonempty meets belong to \(K\) and converge to
\(b=\bigwedge K\). Consequently
\begin{equation}\label{cc:eq:limit-outside}
 b\in K\subseteq K_1\subseteq M\setminus C.
\end{equation}
For each \(n\), the set \(\intr_\theta K_n\) belongs to
\(\Nbh_\theta(v)\), so
\(\bigwedge K_n\leq\bigwedge\intr_\theta K_n\leq a\),
and hence \(\bigwedge K_n\in W\).
The finite nonempty meets of \(K_n\) converge to
\(\bigwedge K_n\). Since \(W\) is open, there is a finite
nonempty \(F_n\subseteq K_n\) such that
\begin{equation}\label{cc:eq:finite-blocks}
 b_n=\bigwedge F_n\in W\subseteq C.
\end{equation}

List each \(F_n\) as one finite nonempty block and concatenate
the blocks to obtain a sequence \((y_j)_{j\geq1}\). For each
\(m\), this sequence is eventually in \(K_m\). Set
\(z_j=y_j\wedge b\). We prove that \(z_j\to b\). Let \(V\)
be an open neighbourhood of \(b\). The open set
\(m_b^{-1}(V)\) contains \(K\), since \(k\wedge b=b\) for
every \(k\in K\). Some \(K_m\) is contained in
\(m_b^{-1}(V)\). Otherwise the decreasing compact sets
\(K_m\setminus m_b^{-1}(V)\) would be nonempty and would
have nonempty intersection, contrary to
\(K\subseteq m_b^{-1}(V)\). Therefore \(z_j\in V\) eventually,
as required.

For \(N\geq1\), put \(c_N=\bigwedge_{j\geq N}z_j\).
Choose a whole block \(F_n\) occurring after the first
\(N-1\) terms. By \eqref{cc:eq:finite-blocks},
\[
 c_N\leq b_n\wedge b\leq b_n\in C.
\]
Since \(C\) is lower, \(c_N\in C\). The increasing sequence
\((c_N)\) converges to \(c=\bigvee_Nc_N\), so closedness of
\(C\) gives \(c\in C\). The assumed sequence identity, applied
at \(b\) to \(z_j\to b\), gives
\[
 c=\bigvee_N\bigwedge_{j\geq N}z_j=b,
\]
contrary to \eqref{cc:eq:limit-outside}. Thus \(a=v\).

The neighbourhood infima form a directed set, since
\(\bigwedge U,\bigwedge V\leq\bigwedge(U\cap V)\).
Each is way-below \(v\), and their supremum is \(v\).
If \(u_1,u_2\ll v\), this directed family contains an element
above both \(u_1\) and \(u_2\). Hence the entire set of
elements way-below \(v\) is directed and has supremum \(v\).
This proves continuity of \(M\).
\end{proof}

\begin{theorem}\label{cc:reflection}
For a spatial frame \(L\), the following conditions are equivalent:
\begin{enumerate}[label=\textup{(\roman*)}]
\item \(\Sigma(L^{\op})\) is core-compact;
\item \(\Sigma(L^{\op})\) is locally compact;
\item \(\Sigma(L^{\op})\) is stably compact;
\item \(L\) is weak hypercontinuous.
\end{enumerate}
Under these conditions,
\begin{equation}\label{cc:eq:topologies}
 \operatorname{patch}\Sigma(L^{\op})=\lambda(L),
 \qquad \sigma(L^{\op})=\upsilon(L^{\op}).
\end{equation}
\end{theorem}
\begin{proof}
The first three conditions are equivalent for every complete lattice
by \cite[Corollary~3.3]{JJL2016}. Assume \textup{(i)}.
Lemmas~\ref{cc:patch-meet} and~\ref{cc:spatial-seq}, followed by
Theorem~\ref{cc:compact-meet}, show that \(L\) is continuous and
satisfies \eqref{cc:eq:neighbourhood-infima} for \(\theta=\theta_L\).
If \(v\in H\in\sigma(L)\), this directed approximation gives
an open neighbourhood \(U\) of \(v\) with \(\bigwedge U\in H\),
whence \(U\subseteq H\). Thus \(\sigma(L)\subseteq\theta_L\).
Closedness of the order gives
\(\upsilon(L^{\op})\subseteq\theta_L\), so
\(\lambda(L)\subseteq\theta_L\). The former is Hausdorff and the
latter compact, hence the two topologies coincide.
Every dual Scott-open set is therefore Lawson open and lower in \(L\).
Lemma~\ref{un:lawson-lower} makes it lower open, proving
\textup{(iv)} and \eqref{cc:eq:topologies}.
Conversely, \textup{(iv)} implies \textup{(ii)} by
\cite[Proposition~3.4]{CKL2022}.
\end{proof}

\subsection{Closed-set lattices and spectra}

\begin{corollary}\label{cc:closed-sets}\label{cc:spectrum}
For a \(T_0\)-space \(X\), the following conditions are equivalent:
\begin{enumerate}[label=\textup{(\roman*)}]
\item \(\Sigma\GG(X)\) is core-compact;
\item \(X\) is core-compact and
\(\sigma(\GG(X))=\upsilon(\GG(X))\);
\item \(\operatorname{Sob}(X)\) is a locally compact Scott space,
with the Scott topology of its specialization order.
\end{enumerate}
Under these conditions the Scott topology on \(\GG(X)\) is its
lower Vietoris topology, and \(x\mapsto\cl_X\{x\}\) is continuous
from \(X\) to \(\Sigma\GG(X)\).
Equivalently, for a spatial frame \(L\), the space
\(\Sigma(L^{\op})\) is core-compact if and only if
\(\Spec L\) is locally compact and
\(\tau_L=\sigma(\Spec L)\).
\end{corollary}
\begin{proof}
Apply Theorem~\ref{cc:reflection} to \(L=\OO(X)\), using
\(\GG(X)\cong L^{\op}\) and the equivalence between
core-compactness of \(X\) and continuity of \(\OO(X)\).
The additional spectral equivalence follows from
Theorem~\ref{thm:main}, the Hofmann--Lawson representation, and
\(\Spec\OO(X)\cong\operatorname{Sob}(X)\).
The lower Vietoris and continuity assertions are
\eqref{eq:lower-vietoris} and the principal-closure embedding in
Section~\ref{sec:prelim}.
\end{proof}

The implication \textup{(i)}\(\Rightarrow\)\textup{(ii)} in
Corollary~\ref{cc:closed-sets} answers
\cite[Problem~4.9]{CK23} affirmatively.

\subsection{The lattice--spectrum correspondence}\label{app:spectrum-dictionary}

For web spaces we use the classical criterion
\begin{equation}\label{un:eq:web-open}
 \uparrow(U\cap\downarrow x)\in\OO(X)
 \qquad(U\in\OO(X),\ x\in X),
\end{equation}
equivalent to the coframe law for \(\OO(X)\)
\cite[Theorem~3]{Erne2009}.
A dcpo \(P\) is \emph{meet-continuous} if, for every directed
\(D\subseteq P\),
\begin{equation}\label{un:eq:meet-continuity}
 x\leq\bigvee D
 \quad\Longrightarrow\quad
 x\in\cl_{\Sigma P}(\downarrow x\cap\downarrow D).
\end{equation}
No binary meets are assumed. The following is the standard Scott-space
form of the web criterion \cite[p.~23:11]{Erne2019}.

\begin{lemma}\label{un:web-meet}
For a dcpo \(P\), meet-continuity of \(P\), webness of \(\Sigma P\),
and the coframe law for \(\OO(\Sigma P)\) are equivalent.
\end{lemma}

\begin{theorem}\label{un:classification}
Let \(L\) be a distributive continuous lattice, let
\(X=\Spec L\) carry its hull--kernel topology, and let \(P\) be its
specialization order. Then \(X\) is locally compact and sober.
In each subsequent row of the following table, the two entries are
equivalent.
\begin{center}
\renewcommand{\arraystretch}{1.23}
\begin{tabular}{@{}>{\raggedright\arraybackslash}p{0.36\linewidth}>{\raggedright\arraybackslash}p{0.57\linewidth}@{}}
\hline
Property of \(L\) & Property of \(X\) and \(P\)\\
\hline
No additional assumption & \(X\) is locally compact and sober\\
Weak hypercontinuous & \(X=\Sigma P\)\\
Hypercontinuous & \(X=\Sigma P\) and \(P\) is quasicontinuous\\
Coframe & \(X=\Sigma P\) and \(P\) is meet-continuous\\
Completely distributive & \(X=\Sigma P\) and \(P\) is continuous\\
\hline
\end{tabular}
\end{center}
Here \(X=\Sigma P\) asserts equality of topologies on the given set.
\end{theorem}
\begin{proof}
The Hofmann--Lawson representation \cite{HL78} gives local compactness, sobriety,
and \(L\cong\OO(X)\). The weak row is
Theorem~\ref{thm:main}.

If \(L\) is hypercontinuous, Theorem~\ref{un:hyper-weak} gives weak
hypercontinuity, so \(X=\Sigma P\). The equivalence between
quasicontinuity of \(P\) and hypercontinuity of \(\OO(\Sigma P)\)
\cite[Theorem~2.5\textup{(2)}]{CKL2022} proves both directions of the
hypercontinuous row.

If \(L\) is a coframe, Corollary~\ref{cor:coframe} gives
\(X=\Sigma P\), and Lemma~\ref{un:web-meet} gives meet-continuity of
\(P\). Conversely, that lemma and \(L\cong\OO(X)\) give the
coframe property from the corresponding spatial condition.

A completely distributive lattice is hypercontinuous
\cite[Corollary~1]{Erne2009}. Hence its spectrum is \(\Sigma P\),
and \(P\) is continuous because \(\OO(\Sigma P)\cong L\) is
completely distributive. The converse uses the same classical
equivalence \cite[Theorem~2.5\textup{(1)}]{CKL2022}.
\end{proof}

We also use the classical identity between classes
\cite[Corollary~1\textup{(f),(g)}]{Erne2009}.

\begin{corollary}\label{un:hyper-coframe}
A complete lattice is completely distributive if and only if it is
hypercontinuous and a coframe.
\end{corollary}

\begin{corollary}\label{un:meet-cover}
Every locally compact sober Scott space \(\Sigma P\) is the image of
a locally compact meet-continuous Scott space \(\Sigma Q\) under a
surjective Scott-continuous map induced by a frame embedding.
\end{corollary}
\begin{proof}
Put \(L=\OO(\Sigma P)\). By Theorem~\ref{thm:main}, this is a
distributive weak hypercontinuous lattice. Its join-translation
extension \(h:L\to T(L)\) is a continuous frame and coframe
extension with Scott-continuous right adjoint, by
Theorem~\ref{thm:embedding}.
The spectral map
\[
 \Spec T(L)\longrightarrow\Spec L,\qquad
 \mathcal F\longmapsto h^{-1}[\mathcal F]
\]
is surjective by the prime-lifting argument in the proof of
Theorem~\ref{thm:transfer}. By Theorem~\ref{un:classification},
its domain is \(\Sigma Q\) for a meet-continuous dcpo \(Q\),
and its codomain is \(\Sigma P\). Both spectral topologies are
Scott, so the map is Scott continuous.
\end{proof}

\section{Scott sobrification and sobriety criteria}\label{paper:topology}

We first apply the spectrum theorem to the sobrification of a
core-compact Scott space. For web \(d\)-spaces, we then obtain a
sobriety criterion in terms of local compactness, countable
irreducible sets and the Baire property.

\subsection{Sobrification of core-compact Scott spaces}\label{topx:reflection-section}

Write \(\operatorname{Irr}_c(X)=\Pc(\GG(X))\) for the nonempty
irreducible closed subsets of \(X\), ordered by inclusion. Its
sobrification topology consists of
\[
 U^*=\{F\in\operatorname{Irr}_c(X):F\cap U\ne\varnothing\},
 \qquad U\in\OO(X).
\]
The unit is \(\eta_X(x)=\downarrow_Xx\).
For a dcpo \(D\), put \(\mathcal S D=\operatorname{Irr}_c(\Sigma D)\).
We use the following established Scott-topological facts
\cite[Propositions~3.4 and~3.7, Lemma~3.6]{CKL2022};
part~\textup{(iii)} is also \cite[Theorem~II-4.13]{GHK03}.

\begin{lemma}\label{fn:scott-facts}
Let \(D\) be a dcpo.
\begin{enumerate}[label=\textup{(\roman*)}]
\item If \(M\) is weak hypercontinuous, then \(\Sigma M^{\op}\)
is locally compact and sober.
\item If \(\Sigma D\) is core-compact, then
\(\sigma(\GG(\Sigma D))=\upsilon(\GG(\Sigma D))\).
\item If \(\Sigma D\) is core-compact, then
\(\Sigma(D\times E)=\Sigma D\times\Sigma E\) for every dcpo \(E\).
\end{enumerate}
\end{lemma}

\begin{theorem}\label{fn:sobrification}
Let \(D\) be a dcpo with core-compact Scott space. Then
\(\OO(\Sigma D)\) is a distributive weak hypercontinuous lattice and
\[
 \operatorname{Sob}(\Sigma D)\cong\Sigma(\mathcal S D).
\]
This Scott space is locally compact and sober. Directed suprema are
given by
\begin{equation}\label{fn:eq:irr-sup}
 \sup_{i\in I}F_i=
 \cl_{\Sigma D}\Bigl(\bigcup_{i\in I}F_i\Bigr).
\end{equation}
\end{theorem}
\begin{proof}
Put \(L=\OO(\Sigma D)\). Core-compactness makes \(L\) continuous,
and Lemma~\ref{fn:scott-facts}\textup{(ii)}, together with
\(L^{\op}\cong\GG(\Sigma D)\), gives weak hypercontinuity.
Apply Theorem~\ref{thm:main} and the canonical identification
\(\Spec\OO(X)\cong\operatorname{Sob}(X)\).
The Hofmann--Lawson representation gives local compactness;
\eqref{fn:eq:irr-sup} is the usual directed-supremum formula for
irreducible closed sets \cite[Remark~2.1\textup{(1)}]{ZhaoXu2018}.
\end{proof}

\subsection{Locally compact web \texorpdfstring{\(d\)}{d}-spaces}
\label{top:web-dspaces}

For the common-lower-bound relation \(R_U\) in
\eqref{ap:eq:R}, consider the condition
\begin{equation}\label{top:eq:joint}
 R_U\text{ is open in }X\times X\quad(U\in\OO(X)).
\end{equation}
The product here is the ordinary topological product.

\begin{lemma}\label{top:pullback}
Every core-compact \(T_0\) web space satisfies
\eqref{top:eq:joint}.
\end{lemma}
\begin{proof}
The frame \(\OO(X)\) is continuous and a coframe. Thus
\(S=\operatorname{Sob}(X)\) is a locally compact sober
meet-continuous Scott space, by Corollary~\ref{cor:coframe} and
Lemma~\ref{un:web-meet}. If \(S=\Sigma P\), meet-continuity makes
\(R_V\) Scott open in \(P^2\). Indeed, let a directed family have
supremum \((s,t)\), and choose \(u\in V\) below \(s,t\).
Meet-continuity for its first projection gives \(u_1\in V\)
below \(u\) and one first coordinate; applying it to the second
projection gives \(u_2\in V\) below \(u_1\) and one second
coordinate. A common upper member of the directed family has both
coordinates above \(u_2\). By
Lemma~\ref{fn:scott-facts}\textup{(iii)}, \(R_V\) is product open.
For \(U\in\OO(X)\),
\[
 R_U^X=(\eta_X\times\eta_X)^{-1}(R_{U^*}^{S}).
\]
Indeed, a common lower bound in \(S\) is an irreducible closed set
contained in both principal closures and meeting \(U\); any point
of that intersection is a witness in \(X\). The displayed inverse
image is open.
\end{proof}

The Scott-space case of the next lemma is He's theorem
\cite[Theorem~5.17]{He2026}. We give the reduction that extends
his argument to the present topological hypotheses.

\begin{lemma}\label{top:sobriety}
Every locally compact \(T_0\) \(d\)-space satisfying
\eqref{top:eq:joint} is sober.
\end{lemma}
\begin{proof}
Write \(h(A)=\intr(\uparrow A)\).
Condition~\eqref{top:eq:joint} has two consequences used in
He's compact-selection argument. First,
\(\uparrow(U\cap\downarrow b)\) is open, being a section of
\(R_U\). Second, for \(x\in U\) and each integer \(r\geq1\),
there is an open neighbourhood \(V\subseteq U\) of \(x\) such
that every nonempty subset of \(V\) of size at most \(r\) has
a lower bound in \(U\). To see this, successively shrink
neighbourhoods \(V_j\) so that \(V_0\subseteq U\) and
\(V_{j+1}^2\subseteq R_{V_j}\); take \(2^j\geq r\) and pair
the entries of a finite tuple recursively.

These are precisely the order-topological inputs to
\cite[Lemmas~5.7--5.11 and the proof of Theorem~5.12]{He2026}:
the remaining steps use relatively closed refinements in compact
sets, the finite intersection property, and a finitely branching
tree. Thus that compact-selection argument applies unchanged and
gives
\[
 K_{n+1}\subseteq h(K_n),\quad K_n\ne\varnothing\text{ compact}
 \quad\Longrightarrow\quad
 \exists p_n\in K_n\quad p_n\leq p_{n+1}.
\]
This invocation does not assume that the original topology is Scott.

Closed subspaces inherit local compactness, the \(d\)-space
property and \eqref{top:eq:joint}. For the last assertion, lowerness
of a closed set \(F\) gives
\(R^{F}_{U\cap F}=R^X_U\cap(F\times F)\).
Let \(F\) be nonempty, irreducible and closed. Given nonempty
relative opens \(O_n\), choose compact saturated sets with
\[
 \varnothing\ne\intr_F K_n\subseteq K_n\subseteq O_n,
 \qquad K_{n+1}\subseteq\intr_F K_n.
\]
Compact selection and the \(d\)-space axiom give an increasing
sequence \(p_n\in K_n\) whose supremum belongs to every \(O_n\).

We recall the final minimal-closed-set step of
\cite[Lemma~5.16]{He2026}, checking its application here.
For nonempty relative opens \(H_n\subseteq F\), choose
\[
 \varnothing\ne W_n\subseteq Q_n\subseteq W_{n-1}\cap H_n,
 \qquad W_{-1}=F,
\]
with \(W_n\) open and \(Q_n\) compact saturated in \(F\).
By compactness and Zorn's lemma, there is a minimal closed
\(B\subseteq F\) meeting all \(Q_n\). It is irreducible: for
any open \(O\) meeting \(B\), minimality gives
\(\varnothing\ne B\cap Q_j\subseteq O\) for some \(j\), and
the \(Q_j\) are nested. The preceding countable-intersection
argument in \(B\) gives \(b\in B\cap\bigcap_nW_n\).
Minimality then implies \(B=\downarrow_F b\).
The open neighbourhoods
\(\uparrow_F(W_n\cap\downarrow_F b)\) form a countable local
base at \(b\), and \(b\in\bigcap_nH_n\).

Consequently the first-countable points of \(F\) are dense, by
taking a constant sequence \(H_n\), and directed, by interleaving
countable local bases at two such points. Their supremum \(p\)
belongs to \(F\), since \(X\) is a \(d\)-space. Density and
closedness of \(\downarrow p\) give
\(F=\downarrow p\). Thus \(X\) is sober.
\end{proof}

\begin{theorem}\label{top:web-equivalences}
For a locally compact \(T_0\) web space \(X\), the following
conditions are equivalent:
\begin{enumerate}[label=\textup{(\roman*)}]
\item \(X\) is a \(d\)-space;
\item \(X\) is sober;
\item \(X\) is the Scott space of its specialization dcpo;
\item \(X\) is the Scott space of a meet-continuous dcpo.
\end{enumerate}
\end{theorem}
\begin{proof}
Lemmas~\ref{top:pullback} and~\ref{top:sobriety} prove
\textup{(i)}\(\Rightarrow\)\textup{(ii)}.
Theorem~\ref{thm:web} proves
\textup{(ii)}\(\Rightarrow\)\textup{(iii)}, and
Lemma~\ref{un:web-meet} proves
\textup{(iii)}\(\Rightarrow\)\textup{(iv)}.
Every Scott space of a dcpo is a \(d\)-space, giving the last
implication.
\end{proof}

\subsection{Countable and Baire criteria}\label{topx:countable-section}

For \(A\in\operatorname{Irr}_c(X)\), let
\(\mathcal N_A=\{U\in\OO(X):A\cap U\ne\varnothing\}\).
This filter has a countable base precisely when the corresponding
point of \(\operatorname{Sob}(X)\) is first countable.
We use \(\omega\)-well-filteredness in the equivalent sequential
form: if decreasing nonempty compact saturated sets \(K_n\)
satisfy \(\bigcap_nK_n\subseteq U\), where \(U\) is open, then
\(K_n\subseteq U\) for some \(n\)
\cite[Proposition~3.10]{XSXZ2020}. Well-filteredness allows arbitrary
filtered families. A space is \emph{completely Baire} if all its
closed subspaces are Baire.

\begin{theorem}\label{topx:countable-test}\label{topr:baire-characterization}
Let \(X\) be a core-compact \(T_0\) web \(d\)-space. The following
conditions are equivalent:
\begin{enumerate}[label=\textup{(\roman*)}]
\item \(X\) is sober;
\item every first-countable point of \(\operatorname{Sob}(X)\)
belongs to \(\eta_X[X]\);
\item every \(A\in\operatorname{Irr}_c(X)\) for which
\(\mathcal N_A\) has a countable base is the closure of a singleton;
\item every countable irreducible subset of \(X\) has the same
closure as a singleton;
\item \(X\) is locally compact;
\item \(X\) is \(\omega\)-well-filtered;
\item \(X\) is well-filtered;
\item \(X\) is completely Baire;
\item every nonempty irreducible closed subspace of \(X\) is Baire.
\end{enumerate}
\end{theorem}
\begin{proof}
Since \(\OO(X)\) is a continuous frame and a coframe,
\(S=\operatorname{Sob}(X)\) is a locally compact sober web
Scott space by Corollary~\ref{cor:coframe}. Let \(G\) be its
first-countable points.
Lemma~\ref{lem:fc-core} gives
\[
 G\cap\downarrow_Ss\text{ nonempty and directed},
 \qquad s=\bigvee_S(G\cap\downarrow_Ss)\quad(s\in S).
\]
The image \(\eta_X[X]\) is closed under directed suprema in \(S\):
if \(D\subseteq X\) is directed and \(d=\bigvee_XD\), then the
\(d\)-space axiom gives
\[
 \bigvee_{x\in D}\eta_X(x)=\cl_XD=\downarrow_Xd=\eta_X(d).
\]
Hence \textup{(ii)} makes \(\eta_X\) surjective, proving
\textup{(i)}. The converse is immediate, and
\textup{(ii)} and \textup{(iii)} are equivalent by definition.
Sobriety implies \textup{(iv)}. If \(A\) is as in
\textup{(iii)}, choose a decreasing base \((U_n)\) for
\(\mathcal N_A\) and points \(x_n\in A\cap U_n\).
Every open set meeting \(A\) contains some \(U_n\), so
\(\cl_X\{x_n:n\in\NN\}=A\). This countable set is irreducible;
thus \textup{(iv)} implies \textup{(iii)}.

For a core-compact sober space, local compactness is classical.
The converse here is Theorem~\ref{top:web-equivalences}.
The established implications
\[
 \textup{(i)}\Longrightarrow\textup{(vii)}
 \Longrightarrow\textup{(vi)}\Longrightarrow\textup{(v)}
\]
follow from sobriety implying well-filteredness and
\cite[Theorem~3.13]{XSXZ2020}.

A locally compact sober space is completely Baire
\cite[Proposition~2.37]{Martin1999}, so
\textup{(i)} implies \textup{(viii)}, which implies \textup{(ix)}.
If \(X\) is not sober, the already proved equivalence gives a
nonprincipal \(A\in\operatorname{Irr}_c(X)\) with a decreasing
countable base \((U_n)\) for \(\mathcal N_A\).
The sets \(V_n=A\cap U_n\) form a decreasing countable
\(\pi\)-base of \(A\) and are dense in \(A\).
If \(a\in\bigcap_nV_n\), every nonempty open subset of \(A\)
contains \(a\), giving \(A=\cl_X\{a\}\), a contradiction.
Thus \(A\) is not Baire, proving
\textup{(ix)}\(\Rightarrow\)\textup{(i)}.
\end{proof}

In condition~\textup{(iv)}, the countable irreducible subset need
not be closed. No countability assumption on \(X\) or its
sobrification is imposed.

\section{Function lattices and retracts}\label{paper:functions}

We classify function lattices and compare Scott and bi-Scott
retracts.

\subsection{Function lattices}\label{app:functions}

For a dcpo \(D\) and a complete lattice \(M\), write \([D\to M]\)
for the Scott-continuous maps, ordered pointwise. This is a complete
lattice; its joins are pointwise, while its meets need not be
\cite[Theorem~2.1]{Erne2013}. For distributive continuous \(M\), the
standard function-space representation gives
\begin{equation}\label{fn:eq:open-representation}
 [D\to M]\cong\OO(\Sigma D\times\Spec M);
\end{equation}
see \cite[Theorem~II-4.10\textup{(2),(6)}]{GHK03} and
\(M\cong\OO(\Spec M)\).

\Needspace{14\baselineskip}
\begin{theorem}\label{fn:classification}
Let \(D\) be a nonempty dcpo and \(M\) a nontrivial distributive
continuous lattice. Put \(A=[D\to M]\). The entries in each row
are equivalent.
\begin{center}
\renewcommand{\arraystretch}{1.18}
\begin{tabular}{@{}l>{\raggedright\arraybackslash}p{0.27\linewidth}>{\raggedright\arraybackslash}p{0.60\linewidth}@{}}
\hline
 & Property of \(A\) & Necessary and sufficient condition\\
\hline
\textup{(i)} & Continuous & \(\Sigma D\) is core-compact\\
\textup{(ii)} & Weak hypercontinuous & \(\Sigma D\) is core-compact and \(M\) is weak hypercontinuous\\
\textup{(iii)} & Hypercontinuous & \(D\) is quasicontinuous and \(M\) is hypercontinuous\\
\textup{(iv)} & Coframe & \(D\) is meet-continuous and \(M\) is a coframe\\
\textup{(v)} & Completely distributive & \(D\) is continuous and \(M\) is completely distributive\\
\hline
\end{tabular}
\end{center}
In \textup{(iv)}, continuity of \(A\) is not assumed; a continuous
coframe is obtained by combining \textup{(i)} and \textup{(iv)}.
\end{theorem}
\begin{proof}
Part \textup{(i)} is the classical continuity theorem
\cite[Theorem~II-4.7]{GHK03}. Part \textup{(v)} is
\cite[Theorem~3.1]{Erne2013}, together with
\(\OO(\Sigma D)\) completely distributive if and only if \(D\)
is continuous \cite[Theorem~2.5\textup{(1)}]{CKL2022}.

For \textup{(ii)}, suppose first that \(\Sigma D\) is
core-compact and \(M\) is weak hypercontinuous. Write
\(\Spec M=\Sigma P\), using Theorem~\ref{thm:main}.
This Scott space is locally compact. Hence
\(\Sigma(D\times P)=\Sigma D\times\Sigma P\) is core-compact,
by Lemma~\ref{fn:scott-facts}\textup{(iii)}. Now
\eqref{fn:eq:open-representation} and
Theorem~\ref{fn:sobrification} make \(A\) weak hypercontinuous.
Conversely, weak hypercontinuity of \(A\) implies core-compactness
of \(\Sigma D\) by \textup{(i)}. Define
\[
 c:M\to A,\quad c(a)(d)=a,
 \qquad q:A\to M,\quad q(f)=\bigwedge_{d\in D}f(d).
\]
Then \(c\dashv q\) and \(qc=\id_M\). Thus \(q\) is
dual-Scott continuous, whereas \(c\) is continuous for the lower
topologies, since
\[
 c^{-1}(A\setminus\uparrow f)
 =M\setminus\uparrow\bigvee_{d\in D}f(d).
\]
For \(U\in\sigma(M^{\op})\), it follows that
\(q^{-1}(U)\in\sigma(A^{\op})=\upsilon(A^{\op})\), and then
\(U=c^{-1}q^{-1}(U)\in\upsilon(M^{\op})\), as required.

For \textup{(iii)}, constant maps and evaluation exhibit \(M\)
as a retract of \(A\) by maps preserving arbitrary joins.
The maps \(j:2\to M\), \(j(0)=0\), \(j(1)=1\), and
\(p:M\to2\), \(p(a)=0\) if and only if \(a=0\), also
preserve arbitrary joins and satisfy \(pj=\id_2\).
Postcomposition therefore exhibits \(\OO(\Sigma D)=[D\to2]\)
as another such retract. Hypercontinuity passes to these retracts
\cite[Lemma~3.16]{CKL2022}, and
\cite[Theorem~2.5\textup{(2)}]{CKL2022} gives the necessary
quasicontinuity of \(D\). Sufficiency is the function-lattice
theorem of Yuan and Kou \cite[Corollary~3.2]{YuanKou2010},
which does not require distributivity of the target.

Finally, \eqref{fn:eq:open-representation} and the web-space
characterization \eqref{eq:web-coframe} reduce \textup{(iv)} to
the fact that a nonempty finite product is a web space exactly
when both factors are web spaces
\cite[Proposition~1.4]{Erne2019}. Apply
Lemma~\ref{un:web-meet} to \(\Sigma D\), and use
\(M\cong\OO(\Spec M)\) for the other factor.
\end{proof}

The Cartesian closedness of weak hypercontinuous lattices and of
their full distributive subcategory, with Scott-continuous maps,
has been established separately.

\subsection{Scott and bi-Scott retracts}\label{app:categories}

A retraction consists of maps \(i:M\to A\) and \(r:A\to M\)
with \(ri=\id_M\). A map is \emph{bi-Scott continuous} if it
preserves directed suprema and filtered infima. Unless a stronger
property is explicitly stated, neither map is required to preserve
finite lattice operations. A monotone retract of a complete lattice
is complete, with
\begin{equation}\label{cat:eq:retract-operations}
 \bigvee_M H=r\Bigl(\bigvee_A i[H]\Bigr),\qquad
 \bigwedge_M H=r\Bigl(\bigwedge_A i[H]\Bigr).
\end{equation}
Scott retracts of continuous lattices are continuous
\cite[Corollary~I-2.3]{GHK03}.

The following retraction will be used in each row of the
classification. For distributive \(M\), its adjunction occurs in
\cite[Proposition~4.5]{CKL2022}.

\begin{proposition}\label{fn:envelope}
Let \(M\) be weak hypercontinuous and put
\(N=\OO(\Sigma M^{\op})\). Then \(N\) is a distributive weak
hypercontinuous lattice, and the maps
\begin{equation}\label{fn:eq:envelope}
 I:M\to N,\quad I(a)=M\setminus\uparrow a,
 \qquad R:N\to M,\quad R(U)=\bigwedge(M\setminus U)
\end{equation}
are bi-Scott continuous, with \(I\dashv R\) and \(RI=\id_M\).
\end{proposition}
\begin{proof}
Lemma~\ref{fn:scott-facts}\textup{(i)} and
Theorem~\ref{fn:sobrification} give the assertion about \(N\).
For \(a\in M\) and \(U\in N\),
\[
 I(a)\subseteq U
 \quad\Longleftrightarrow\quad M\setminus U\subseteq\uparrow a
 \quad\Longleftrightarrow\quad a\leq R(U).
\]
Also \(RI(a)=\bigwedge\uparrow a=a\). Thus \(I\dashv R\),
\(I\) preserves arbitrary joins, and \(R\) preserves arbitrary
meets.
Under complementation \(N^{\op}\cong\GG(\Sigma M^{\op})\),
the opposite-order map \(I\) becomes the Scott-continuous
principal-ideal map \cite[Proposition~3.10\textup{(1)}]{CKL2022}.
Hence \(I\) is bi-Scott continuous.

Let \(a\ll_M b\) and put
\[
 K=M\setminus\intr_{\sigma(M)}\uparrow a.
\]
Since \(b\in\intr_{\sigma(M)}\uparrow a\) and this interior is
upper, \(I(a)\subseteq K\subseteq I(b)\). The set \(K\) is
closed in the compact Lawson space of \(M\). Since
\(\sigma(M^{\op})=\upsilon(M^{\op})\subseteq\lambda(M)\),
it is compact in \(\Sigma M^{\op}\). Consequently
\(I(a)\ll_N I(b)\). Lemma~\ref{ap:adjoint-waybelow} now gives
Scott continuity of \(R\), so \(R\) is bi-Scott continuous.
\end{proof}

\begin{theorem}\label{cat:karoubi}\label{cat:retract-table}
For each class in the first column, its retracts are exactly the
lattices in the corresponding second or third column. The Scott
column uses Scott-continuous maps; the bi-Scott column uses
bi-Scott-continuous maps.
\begin{center}
\renewcommand{\arraystretch}{1.2}
\begin{tabular}{@{}>{\raggedright\arraybackslash}p{0.31\linewidth}>{\raggedright\arraybackslash}p{0.20\linewidth}>{\raggedright\arraybackslash}p{0.41\linewidth}@{}}
\hline
Class of ambient lattices & Scott retracts & Bi-Scott retracts \(M\)\\
\hline
Distributive weak hypercontinuous & All continuous lattices & Weak hypercontinuous lattices\\
Distributive hypercontinuous & All continuous lattices & Hypercontinuous lattices\\
Distributive hyperalgebraic & All continuous lattices & Hyperalgebraic lattices\\
Continuous frames that are coframes & All continuous lattices & Weak hypercontinuous lattices with \(M^{\op}\) meet-continuous\\
Completely distributive & All continuous lattices & Weak hypercontinuous lattices with \(M^{\op}\) continuous\\
\hline
\end{tabular}
\end{center}
Every row also identifies the idempotent completion of the full
category in the first column with the full category in the
corresponding column, using the indicated maps.
\end{theorem}
\begin{proof}
The Scott column follows from the classical characterization of
continuous lattices as Scott retracts of powerset lattices
\cite[Corollary~I-4.18\textup{(4)}]{GHK03}: powerset lattices
belong to all five ambient classes, and all five classes consist
of continuous lattices. The same observation applies to algebraic
lattices, distributive algebraic lattices, and powerset lattices
themselves.

Suppose that \(i:M\to A\) and \(r:A\to M\) are bi-Scott
continuous, \(ri=\id_M\), and \(A\) is weak hypercontinuous.
Then \(M\) is continuous. The surjection \(r\) sends the compact
space \((A,\beta(A))\) continuously onto \((M,\beta(M))\), so
\(\beta(M)\) is compact. Since
\(\lambda(M)\subseteq\beta(M)\) and \(\lambda(M)\) is
Hausdorff, the two topologies coincide. By
Proposition~\ref{un:weak-equivalences}, \(M\) is weak
hypercontinuous. Conversely,
Proposition~\ref{fn:envelope} realizes every weak hypercontinuous
\(M\) as a bi-Scott retract of a distributive one.

For the hypercontinuous row, use
Theorem~\ref{un:hyper-weak}: it adds quasicontinuity of the
opposite order to weak hypercontinuity. Quasicontinuity is
preserved by Scott retracts
\cite[Proposition~2.7.8]{Jia2018Thesis}. Conversely, when
\(M^{\op}\) is quasicontinuous,
\(N=\OO(\Sigma M^{\op})\) in
Proposition~\ref{fn:envelope} is hypercontinuous
\cite[Theorem~2.5\textup{(2)}]{CKL2022}.

For the hyperalgebraic row, suppose \(M\) is a bi-Scott retract
of a hyperalgebraic lattice \(A\). Both are weak hypercontinuous,
so their bi-Scott and Lawson topologies agree. The section embeds
the compact Lawson space of \(M\) as a closed subspace of that
of \(A\). Algebraicity of \(A\) makes its Lawson topology
zero-dimensional, and hence so is that of \(M\). The classical
zero-dimensionality characterization
\cite[Theorem~2.11]{Hofmann1977} applies because \(M\), being
continuous, also has Lawson-continuous meet. Thus \(M\) is algebraic.
Corollary~\ref{un:algebraic} gives hyperalgebraicity. Conversely,
if \(M\) is hyperalgebraic, then
\(N=\sigma(M^{\op})=\upsilon(M^{\op})\) is algebraic by
\cite[Lemma~3.1]{CKL2022}, since an algebraic lattice is
quasialgebraic. The envelope \(N\) is therefore distributive and
hyperalgebraic, again by Corollary~\ref{un:algebraic}.

For the coframe row, every ambient lattice is weak hypercontinuous
by Theorem~\ref{un:classification}, and its opposite order is
meet-continuous. Meet-continuity is inherited by Scott retracts
\cite[Proposition~3.1.7]{Jia2018Thesis}. Conversely, if \(M\)
is weak hypercontinuous and \(M^{\op}\) is meet-continuous,
then \(N=\OO(\Sigma M^{\op})\) is also a coframe by
Lemma~\ref{un:web-meet}.

For the last row, complete distributivity gives continuity in
both orders, so necessity follows from the first bi-Scott row and
the continuous-retract theorem in the opposite order. Conversely,
continuity of \(M^{\op}\) makes
\(\OO(\Sigma M^{\op})\) completely distributive
\cite[Theorem~2.5\textup{(1)}]{CKL2022}, and the same
retraction applies.

Finally, an idempotent of the indicated kind splits through its
fixed-point poset, with operations as in
\eqref{cat:eq:retract-operations}. Each target class is closed
under the relevant retracts, by the same arguments, and every
object in it is a retract of a source object. This is the
idempotent-completion criterion \cite[Section~5.8]{Kelly1982}.
\end{proof}

In particular, let \(\mathbf{ScSpec}\) denote the distributive
continuous lattices with intrinsic Scott spectra, so that
\(\mathbf{ScSpec}\) and the distributive weak hypercontinuous
lattices have the same objects. With \(\mathbf{ContLat}\) and
\(\mathbf{WHyp}\) denoting continuous and weak hypercontinuous
lattices respectively, the first row gives
\begin{equation}\label{cat:eq:karoubi}
 \operatorname{Kar}(\mathbf{ScSpec}_{Scott})
 \simeq\mathbf{ContLat}_{Scott},\qquad
 \operatorname{Kar}(\mathbf{ScSpec}_{bi})
 \simeq\mathbf{WHyp}_{bi}.
\end{equation}
Here \(\operatorname{Kar}\) denotes the usual category obtained
by splitting idempotents \cite[Section~5.8]{Kelly1982}.

For finite ambient lattices the distinction disappears: every
monotone map is bi-Scott continuous, and the retracts of finite
distributive lattices, or of finite powerset lattices, are exactly
all finite lattices. For the nontrivial direction, a finite lattice
\(M\) is a retract of \(\mathcal P(M)\) via
\(i(x)=\downarrow x\) and \(r(H)=\bigvee H\).

\section{The join-translation extension and lattice structure}\label{paper:auxiliary}

For a distributive continuous lattice \(L\), retain the construction
\(T(L)\), the embedding \(h_L(a)=c_a\), which preserves arbitrary
joins and finite meets, and its Scott-continuous right adjoint
\(r_L(t)=t(0)\) from Section~\ref{sec:host}.
Lemma~\ref{lem:hostbase} gives continuity of \(T(L)\), pointwise
meets and directed joins, and Lawson density of \(h_L[L]\), without
assuming weak hypercontinuity. We give a universal property of this
specific extension and describe the lattice properties that it detects.

\Needspace{17\baselineskip}
\subsection{A universal property of dense extensions}

\begin{proposition}\label{ap:canonical-quotient}
Let \(L\) be a distributive weak hypercontinuous lattice, let \(M\)
be a continuous frame-coframe, and let \(e:L\to M\) be a frame
embedding with Scott-continuous right adjoint \(r_e\). Put
\[
 S=\cl_{\lambda(M)}e[L].
\]
Then \(S\) is a continuous frame-coframe, and there is a unique
surjective complete lattice homomorphism \(q:S\to T(L)\) with
\(qe=h_L\). It is given by
\begin{equation}\label{ap:eq:canonical-quotient}
 q(s)(x)=r_e\bigl(s\vee_M e(x)\bigr)
 \qquad(s\in S,\ x\in L).
\end{equation}
\end{proposition}

\begin{proof}
Both finite lattice operations on \(M\) are Lawson continuous
(Lemma~\ref{ap:joint}). Thus the closure of the bounded sublattice
\(e[L]\) is again a bounded sublattice. Increasing nets of finite
joins and decreasing nets of finite meets show that \(S\) is a
complete sublattice of \(M\). It inherits the frame and coframe
laws. Its inclusion \(j:S\hookrightarrow M\) has a left adjoint
and preserves directed joins, so \(S\) is a Scott retract of \(M\)
and hence continuous. Moreover, \(\lambda(S)\) is the subspace
topology from \(\lambda(M)\)
\cite[Theorem~III-1.8]{GHK03}.

Since \(r_e\) preserves arbitrary meets and directed joins, it is
Lawson continuous; also \(r_e e=\id_L\). Formula
\eqref{ap:eq:canonical-quotient} therefore defines a Lawson-continuous
map \(S\to L^L\) with \(q(e(a))=c_a\). Density and compactness give
\(q[S]=T(L)\). The coframe law in \(M\), meet preservation by
\(r_e\), and the pointwise meet formula in \(T(L)\) give
\[
 q\left(\bigwedge_S A\right)(x)
 =r_e\left(\bigwedge_{s\in A}(s\vee e(x))\right)
 =\bigwedge_{s\in A}q(s)(x)
 \qquad(A\subseteq S).
\]
For \(s\in S\) and \(a,x\in L\),
\[
 q(s\vee e(a))(x)=r_e(s\vee e(a\vee x))=q(s)(a\vee x).
\]
Consequently \(q(s\vee e(a))=q(s)\vee_{T(L)}c_a\), by
Lemma~\ref{ap:join}. Taking a net \(e(a_i)\to t\) and using
Lawson continuity of joins in \(S\) and \(T(L)\) gives
\(q(s\vee t)=q(s)\vee q(t)\). Also \(q(0_S)=c_0=0_{T(L)}\).
Increasing nets of finite joins now give preservation of all joins.
Finally, every complete lattice homomorphism between continuous
lattices is Lawson continuous \cite[Theorem~III-1.8]{GHK03}; agreement
on the dense set \(e[L]\) therefore proves uniqueness.
\end{proof}

Equivalently, \(h_L:L\to T(L)\) is terminal among frame embeddings
\(e:L\to S\) into continuous frame-coframes whose right adjoints
are Scott continuous and whose images are Lawson dense, with
commuting complete lattice homomorphisms as morphisms. This
universal property concerns the specified extensions of a fixed
\(L\).

Canonical extensions of frames, as introduced in \cite{Jakl2020},
satisfy the compactness condition
\[
 \bigwedge e[F]\leq e(a)\quad\Longrightarrow\quad a\in F
 \qquad(F\text{ a Scott-open filter});
\]
see \cite[Section~3.1.1]{JaklSuarez2025}. The join-translation
extension need not satisfy this condition. For \(L=[0,1]\) in
its usual order, the map \(a\mapsto c_a\) is Lawson continuous,
so compactness gives \(T(L)=h_L[L]\). If \(0<r<1\), then
\(F=(r,1]\) is a Scott-open filter with
\(\bigwedge h_L[F]=h_L(r)\) and \(r\notin F\).

\subsection{Structural correspondences}

We use hypercontinuity as in Section~\ref{app:unified}. A lattice
is superalgebraic if it is generated under joins by completely
join-prime elements; equivalently, it is algebraic and completely
distributive \cite[Corollary~1]{Erne2009}.

\begin{theorem}\label{auxs:structure}
Let \(L\) be a distributive continuous lattice. In each row of the
following table, the two entries are equivalent.
\begin{center}
\renewcommand{\arraystretch}{1.18}
\begin{tabular}{@{}>{\raggedright\arraybackslash}p{0.38\linewidth}>{\raggedright\arraybackslash}p{0.55\linewidth}@{}}
\hline
Property of \(L\) & Property of its extension\\
\hline
Weak hypercontinuous & \(T(L)\) is a coframe\\
Hypercontinuous & \(T(L)\) is completely distributive\\
Coframe & \(T(L)=h_L[L]\)\\
Boolean & \(T(L)\) is Boolean\\
\hline
\end{tabular}
\end{center}
The coframe row is also equivalent to Lawson continuity of
\(h_L\), and to bi-Scott continuity of \(h_L\). If \(L\) is weak
hypercontinuous, the hypercontinuous row is also equivalent to
continuity of \(T(L)^{\op}\) and to hypercontinuity of \(T(L)\).
\end{theorem}

\begin{proof}
The forward implication of the weak row is
Theorem~\ref{thm:embedding}. Conversely, if \(T(L)\) is a coframe,
it is finitely distributive and continuous, hence a frame. Its spectrum
is a Scott space by Corollary~\ref{cor:coframe}. Apply
Theorem~\ref{thm:transfer} to \(h_L\), and then
Theorem~\ref{thm:main}, to obtain weak hypercontinuity of \(L\).

For the coframe row, if \(T(L)=h_L[L]\), write
\(\bigwedge_{a\in A}c_a=c_b\). Evaluating at \(0\) and at \(x\)
gives \(b=\bigwedge A\) and
\[
 x\vee\bigwedge A=\bigwedge_{a\in A}(x\vee a).
\]
Conversely, if \(L\) is a coframe, its join translations are Lawson
continuous by Lemma~\ref{ap:joint}. Thus \(h_L[L]\) is compact and
Lawson closed in \(L^L\), so it equals \(T(L)\). The same
compactness argument shows that Lawson continuity of \(h_L\)
implies surjectivity; surjectivity makes \(h_L\) an order
isomorphism and hence Lawson and bi-Scott continuous. Conversely,
if \(h_L\) is bi-Scott continuous, finite-meet preservation and
preservation of filtered infima give, for every \(A\subseteq L\),
\[
 h_L\left(\bigwedge A\right)
 =\bigwedge_{F\in\Fin(A)}^{T(L)}h_L\left(\bigwedge F\right)
 =\bigwedge_{a\in A}^{T(L)}h_L(a).
\]
Evaluating at \(x\) gives the coframe identity above. This argument
includes \(A=\varnothing\), since \(\Fin(A)\) contains the empty set.

We prove the hypercontinuous row by explicit completely-below
approximation. Write \(u\triangleleft v\) when every set \(A\) with
\(v\leq\bigvee A\) contains an \(a\) with \(u\leq a\). Complete
distributivity is equivalent to
\(v=\bigvee\{u:u\triangleleft v\}\)
\cite[Corollary~1\textup{(b),(g)}]{Erne2009}. Suppose that \(L\) is
hypercontinuous. Then it is weak hypercontinuous, and \(T=T(L)\)
is a coframe. For \(p\in L\), put
\[
 d_p=\bigwedge_{a\nleq p}^{T}c_a.
\]
If \(u\nleq c_p\), Lawson density supplies a net
\(c_{a_i}\to u\) eventually outside \(\downarrow c_p\), hence
eventually above \(d_p\). The closedness of \(\uparrow d_p\) gives
\(d_p\leq u\). It follows that
\begin{equation}\label{auxs:eq:completely-below}
 t\nleq c_p\quad\Longrightarrow\quad d_p\triangleleft t,
\end{equation}
because \(t\leq\bigvee A\) and \(t\nleq c_p\) force some
\(u\in A\) to lie outside \(\downarrow c_p\). For finite
\(F\subseteq L\), the coframe law and preservation of joins by
\(h_L\) give
\begin{equation}\label{auxs:eq:finite-separation}
 \bigvee_{p\in F}^{T}d_p
 =\bigwedge_{a\in L\setminus\downarrow F}^{T}c_a.
\end{equation}
Indeed, distribute the finite join across the meets: the resulting
translations are \(c_{\bigvee_{p\in F}a_p}\), with \(a_p\nleq p\).
Their parameters run through \(L\setminus\downarrow F\), since
every such parameter is obtained by the constant choice.
For \(F=\varnothing\), both sides are \(c_0\).

Each principal filter \(\uparrow x\) is hypercontinuous: if
\(a\prec y\) and \(y\geq x\), restriction of the upper-open
neighbourhood witnessing \(a\prec y\) gives
\(a\vee x\prec_{\uparrow x}y\). Its upper topology is the
subspace upper topology. The finite-complement neighbourhood base
therefore gives, for \(y\geq x\),
\begin{equation}\label{auxs:eq:relative-hyper}
 y=\bigvee_{\substack{F\subseteq\uparrow x\text{ finite}\\
                       y\notin\downarrow F}}
       \bigwedge\bigl((\uparrow x)\setminus\downarrow F\bigr).
\end{equation}
Every term is below \(y\); conversely, each approximant
\(v\prec_{\uparrow x}y\) is below the infimum of a basic
upper-open neighbourhood contained in \(\uparrow v\).

For \(t\in T\), let \(q=\bigvee\{d_p:t\nleq c_p\}\). Then
\(q\leq t\) by \eqref{auxs:eq:completely-below}. Fix \(x\in L\).
If \(F\subseteq\uparrow x\) is finite and
\(t(x)\notin\downarrow F\), then \(t\nleq c_p\) for every
\(p\in F\), and \eqref{auxs:eq:finite-separation} gives
\[
 q(x)\geq\bigwedge_{a\in L\setminus\downarrow F}(a\vee x)
       =\bigwedge\bigl((\uparrow x)\setminus\downarrow F\bigr).
\]
Apply \eqref{auxs:eq:relative-hyper} to \(y=t(x)\) to obtain
\(q(x)\geq t(x)\). Thus
\[
 t=\bigvee_{t\nleq c_p}^{T}d_p,
 \qquad d_p\triangleleft t\text{ for every term},
\]
which proves complete distributivity of \(T\). Conversely,
\(r_L:T(L)\to L\) preserves arbitrary meets and directed joins
and is surjective. If \(T(L)\) is completely distributive, the
CL-quotient characterization of Gierz and Lawson
\cite[Definition~3.5 and Theorem~6.4\textup{(i)}]{GL81} therefore
makes \(L\) hypercontinuous.

For the Boolean row, the forward implication follows from the
coframe row. Conversely, suppose \(T(L)\) is Boolean. A continuous
Boolean lattice is a powerset lattice, by the classical spatial
representation of continuous frames \cite{HL78}; identify \(T(L)\)
with \(\mathcal P(S)\). The image \(H=h_L[L]\) is a Lawson-dense
subframe. The map \(b=h_Lr_L:\mathcal P(S)\to H\) is Scott
continuous, satisfies \(b(A)\subseteq A\), and has \(b(S)=S\).
Hence
\[
 S=\bigcup_{F\subseteq S\text{ finite}}b(F).
\]
Every \(s\in S\) therefore belongs to some finite member of \(H\).
Lawson density in \(\mathcal P(S)\) implies that for \(s\ne t\)
there is a member of \(H\) containing \(s\) and omitting \(t\).
Finite intersections inside a finite neighbourhood now give
\(\{s\}\in H\). Since \(H\) is closed under arbitrary unions,
\(H=\mathcal P(S)\), and \(L\) is Boolean.

Finally, for a continuous frame-coframe, complete distributivity
is equivalent both to continuity of the dual and to
hypercontinuity \cite[Corollary~1\textup{(f),(g),(i)}]{Erne2009}.
This proves the two further formulations of the hypercontinuous row.
\end{proof}

There is also the elementary equivalence
\(L\text{ finite}\Longleftrightarrow T(L)\text{ finite}\):
use injectivity of \(h_L\) in one direction and the coframe row
in the other. For algebraicity we record only the forward
implications. If \(L\) is algebraic, then \(T(L)\) is algebraic;
if \(L\) is algebraic and weak hypercontinuous, then \(T(L)\)
is superalgebraic. Indeed, \(L^L\) is algebraic, and the inclusion
\(T(L)\hookrightarrow L^L\) preserves arbitrary meets and directed
joins. Its left adjoint therefore preserves compact elements
\cite[Theorem~IV-1.4]{GHK03}; applying it to compact
approximations gives algebraicity of \(T(L)\). The second
assertion now follows from the weak row and the classical
algebraic-coframe theorem \cite[Lemma~1\textup{(7)}]{Erne2009}.

\section*{Supplementary materials}
The supplementary materials and the Lean formalization code will be available at \url{https://github.com/ChanYuxu/Recent-Progress-on--Domain-Theory}.

\end{document}